\documentclass[10pt,leqno]{amsart}
\usepackage{epsfig}
\usepackage{amssymb}
\usepackage{amscd}
\usepackage[matrix,arrow]{xy}

\newtheorem{theo}{Theorem}[section]
\newtheorem{lemma}[theo]{Lemma}
\newtheorem{defi}[theo]{Definition}
\newtheorem{prop}[theo]{Proposition}
\newtheorem{conj}[theo]{Conjecture}
\newtheorem{cor}[theo]{Corollary}
\newtheorem{remark}[theo]{Remark}
\newtheorem{property}[theo]{Property}
\newtheorem{example}[theo]{Example}
\numberwithin{equation}{section}

\def\CP{\mathbb{CP}}
\def\R{\mathbb{R}}
\def\C{\mathbb{C}}
\def\Z{\mathbb{Z}}
\def\FS{\mathrm{Lag_{vc}}}
\def\db#1{ \bD^b({#1})}
\def\coh{\operatorname{coh}}
\def\Qcoh{\operatorname{Qcoh}}
\def\lto{\longrightarrow}
\def\wt{\widetilde}

\def\A{{\mathcal A}}
\def\B{{\mathcal B}}
\def\cC{{\mathcal C}}
\def\D{{\mathcal D}}
\def\E{{\mathcal E}}
\def\F{{\mathcal F}}
\def\G{{\mathcal G}}
\def\O{{\mathcal O}}
\def\P{{\mathcal{P}}}

\def\bH{{\mathbf H}}
\def\bC{{\mathbf C}}
\def\bD{{\mathbf D}}
\def\bR{{\mathbf R}}
\def\bL{{\mathbf L}}

\def\PP{{\mathbb P}}
\def\FF{{\mathbb F}}

\def\Hom{\operatorname{Hom}}
\def\End{\operatorname{End}}
\def\Ext{\operatorname{Ext}}
\def\hom{{{\mathcal H}om}}

\def\kk{{\mathbf k}}
\def\bP{{\mathbf P}}
\def\bA{{\mathbf A}}
\def\Proj{{\mathbf{Proj}}\,}
\def\PROJ{{\mathbb{P}\mathrm{roj}}\,}
\def\bz{{\mathbf 0}}
\def\Spec{{\mathbf{Spec}}\,}

\def\bdot{{\raise1pt\hbox{$\scriptscriptstyle\bullet$}}}
\def\op{{\oplus}}
\def\gr{\operatorname{gr}}
\def\qgr{\operatorname{qgr}}
\def\tors{\operatorname{tors}}
\def\Gr{\operatorname{Gr}}
\def\QGr{\operatorname{QGr}}
\def\Tors{\operatorname{Tors}}
\def\Gm{{\mathbf G}_m}
\def\mod{\operatorname{mod}}
\def\Mod{\operatorname{Mod}}
\def\rad{\operatorname{rad}}
\def\Com{\operatorname{Com}^{\bdot}}
\def\Be{{B_{\theta}}}
\def\Ce{{C_{\theta}^{\bdot}}}
\def\Cc{\mathfrak{C}_{\theta}}
\def\Bc{\mathfrak{B}_{\theta}}
\def\Xe{{X}^{\bdot}}
\def\fun{{MK_n}}
\def\De{{F(n)}}
\def\Te{{\mathcal T}_{\theta}}
\def\cX{\mathcal{X}}

\def\Sch{\operatorname{Sch}}

\def\aa{{\mathbf{\overline{a}}}}
\def\Mat{{\operatorname{M}(n+1, \kk^*)}}

\title[Mirror symmetry for weighted projective planes]
{Mirror symmetry for weighted projective planes
and their noncommutative deformations }

\author{Denis Auroux}
\address{Department of Mathematics, M.I.T., Cambridge MA 02139, USA
\newline\indent
Centre de Math\'ematiques, Ecole Polytechnique, 91128 Palaiseau, France}
\email{auroux@math.mit.edu, auroux@math.polytechnique.fr}
\author{Ludmil Katzarkov}
\address{Department of Mathematics, University of California, Irvine, CA
92697, USA}
\email{lkatzark@math.uci.edu}
\author{Dmitri Orlov}
\address{
Algebra Section, Steklov Mathematical Institute, Russian Academy
of Sciences, \newline\indent 8 Gubkin str., Moscow 119991, Russia}

\email{orlov@mi.ras.ru}

\begin{document}
\maketitle

\tableofcontents

\section{Introduction}

The phenomenon of Mirror Symmetry, in its ``classical'' version, was first
observed  for Calabi-Yau manifolds, and mathematicians were introduced to
it through a series of remarkable papers \cite{GP,CDGP,Vafa,Wit,Dub,Man,...}.
Some very strong conjectures have been made about its topological interpretation
-- e.g.\ the Strominger-Yau-Zaslow conjecture. In a different direction, the framework of
mirror symmetry was extended by Batyrev, Givental, Hori, Vafa, etc.\ to the
case of Fano manifolds.

In this paper, we approach mirror symmetry for Fano manifolds from the
point of view suggested by the work of Kontsevich and his
remarkable Homological Mirror Symmetry (HMS)
conjecture \cite{KICM}. We extend the previous investigations in the
following two directions:

\begin{itemize}
\item  Building on recent works by Seidel \cite{Se1}, Hori and Vafa
\cite{HV} (see also an earlier paper by Witten \cite{Wit2}), we prove HMS for some Fano manifolds, namely weighted projective
lines and planes, and Hirzebruch surfaces. This extends, at a greater level
of generality, a result of Seidel \cite{Se2} concerning the case of the
usual $\CP^2$.

\item  We obtain the first explicit description of the extension
of HMS to noncommutative deformations of Fano algebraic varieties.
\end{itemize}

In the long run, the goal is to explore in greater depth the fascinating
ties brought forth by HMS between complex algebraic geometry and symplectic
geometry, hoping that the currently more developed algebro-geometric
methods will open a fine opportunity for obtaining new interesting results in
symplectic geometry. We first describe the results of this paper in some more detail.
\medskip

Most of the classical works on string theory deal with the case of $N=2$
superconformal sigma models with a Calabi-Yau target space. In this
situation the corresponding field theory has two topologically twisted
versions, the A- and B-models, with D-branes of types A and B respectively.
Mirror symmetry interchanges these two classes of D-branes.
In mathematical terms, the category of B-branes on a Calabi-Yau
manifold $X$ is the derived category of coherent sheaves on $X$,
$\mathbf{D}^b(\coh(X))$. The so-called (derived) Fukaya category $D\mathcal{F}(Y)$ has been proposed
as a candidate for the category of A-branes on a Calabi-Yau manifold $Y$;
in short this is a category whose objects are Lagrangian submanifolds equipped with flat
vector bundles.
The HMS Conjecture claims that if two Calabi Yau manifolds $X$ and $Y$ are
mirrors to each other then $\mathbf{D}^b(\coh(X))$ is equivalent to $D\mathcal{F}(Y)$.

Physicists also consider more general $N=2$ supersymmetric field theories
and the corresponding D-branes; among these, two families of theories
are of particular interest to us: on one hand, sigma models with a Fano
variety as target space, and on the other hand, $N=2$ Landau-Ginzburg models.
Mirror symmetry puts the former
in correspondence with a certain subclass of the latter.
In particular, B-branes on a
Fano variety are described by the derived category of coherent
sheaves, and under mirror symmetry they correspond to the A-branes of
a mirror Landau-Ginzburg model. These A-branes are
described by a suitable analogue
of the Fukaya category, namely the derived category of
Lagrangian vanishing cycles.

In order to demonstrate this feature of mirror symmetry,
we use a procedure introduced by
Batyrev \cite{Bat}, Givental \cite{Giv},
Hori and Vafa \cite{HV}, which we will call the {\it toric mirror ansatz}.
Starting from a complete intersection $Y$ in a toric variety, this
procedure yields a description of an affine subset of its
mirror Landau-Ginzburg model (to obtain a full description of
the mirror it is usually necessary to consider a partial (fiberwise)
compactification) -- an open symplectic
manifold  $(X,\omega)$ and a symplectic fibration $W:X\to \C$
(see e.g.\ \cite{KALI}).

Following ideas of Kontsevich \cite{Ko} and Hori-Iqbal-Vafa \cite{HIV},
Seidel rigorously defined (in the case of non-degenerate critical points)
a derived category of
Lagrangian vanishing cycles $D(\FS(W))$ \cite{Se1},
whose objects represent A-branes on $W:X\to \C$.

In the  case of Fano manifolds the statement of the HMS conjecture is the
following:

\begin{conj}
The category of A-branes $D(\FS(W))$ is  equivalent
to the derived category of coherent sheaves (B-branes) on $Y$.
\end{conj}

We will  prove this conjecture  for various examples.

There is also a parallel statement of HMS relating the derived category of
B-branes on $W:X\to \C$, whose definition was suggested by Kontsevich and
carried out algebraically in \cite{O}, and the derived Fukaya category of
$Y$. Since very little is known about these Fukaya categories, we will not
discuss the details of this statement in the present paper. Our hope in
this direction is that algebro-geometric methods will allow us
to look at Fukaya categories from a different perspective.

The case we will be mainly concerned with in this paper is that of
the weighted projective plane $\CP^2(a,b,c)$ (where $a,b,c$ are three
mutually prime positive integers).
Its mirror is the affine hypersurface
$X=\{x^ay^bz^c=1\}\subset
(\C^*)^3$, equipped with an exact symplectic form $\omega$ and
the superpotential $W=x+y+z$.
Our main theorem is:

\begin{theo}\label{thm:main}
\bfseries{HMS holds for  $\CP^2(a,b,c)$ and its noncommutative
deformations.}
\end{theo}

Namely, we show that the derived category of coherent sheaves (B-branes)
on the weighted projective plane $\CP^2(a,b,c)$ is
equivalent to the derived category of vanishing cycles (A-branes) on
the affine hypersurface $X\subset(\C^*)^3$. Moreover, we also show that
this mirror correspondence between derived categories can be extended to
toric noncommutative deformations of $\CP^2(a,b,c)$ where B-branes are concerned,
and their
mirror counterparts, non-exact deformations of the symplectic structure
of $X$ where A-branes are concerned.

Observe that weighted projective planes are rigid in terms of commutative
deformations, but have a one-dimesional moduli space of toric noncommutative
deformations ($\CP^2$ also has some other noncommutative
deformations, see \S 6.2). We expect a similar phenomenon to hold in many
cases where the toric mirror ansatz applies. An interesting question will
be to extend this correspondence to the case of general noncommutative
toric vareties.

We will also consider some other examples besides weighted projective planes,
in order to demonstrate the ubiquity of HMS:

\begin{itemize}

\item as a warm-up example, we give a proof of HMS for weighted projective lines
(a result also announced by D.\ van Straten in \cite{DVS}).

\item we also discuss HMS for Hirzebruch surfaces $\mathbb{F}_n$.
For  $n\ge 3$, the canonical class is no longer
negative ($\mathbb{F}_n$ is not Fano), and HMS does not hold directly,
because some modifications of the toric mirror ansatz are needed, as
already noticed in \cite{HIV}. The direct application of the
ansatz produces a Landau-Ginzburg model whose derived category of
vanishing cycles is identical to that on the mirror of
the weighted projective plane $\CP^2(1,1,n)$.
In order to make the HMS conjecture work
we need to restrict ourselves to an open subset in
the target space $X$ of this Landau-Ginzburg model.

\item we will also outline an idea of the
proof of HMS (missing only some Floer-theoretic arguments about certain
moduli spaces of pseudo-holomorphic discs)
for some higher-dimensional Fano manifolds,
e.g.\ $\CP^3$.

\end{itemize}

A word of warning is in order here. We do not describe completely and do
not make use of the full potential of the toric mirror ansatz in this paper.
Indeed we do not compactify and desingularize the open manifold $X$.
Compactification and desingularization procedures will be
addressed in full detail in future papers \cite{FAKO} dealing with
the cases of more general Fano manifolds and manifolds of general type,
where these extra steps are needed in order to exhibit the whole category
of D-branes of the Landau-Ginzburg model.
In this paper we work with specific examples for which
compactification and desingularization are not needed (conjecturally
this is the case for all toric varieties). However there are two principles
which are readily apparent from these specific examples:

\begin{itemize}

\item
noncommutative deformations of Fano manifolds are related to variations
of the cohomology class of the symplectic form on the mirror Landau-Ginzburg
models;

\item
even in the toric case, a fiberwise compactification of the Landau-Ginzburg model
is required in order to obtain general noncommutative deformations. The noncompact case
then arises as a limit where the symplectic
form on the compactified fiber acquires poles along the compactification
divisor.

\end{itemize}

Moreover there are two features of HMS for toric varieties,
which become apparent in this paper and which
we would like to emphasize:

\begin{itemize}

\item  it is important to think of singular toric varieties as smooth
quotient stacks. As a consequence of the work of Cox \cite{COX} this
characterization is possible in many cases;

\item  as suggested by our specific examples, we would like to conjecture that
the derived category of coherent sheaves over a smooth toric quotient stack
is always generated by an exceptional collection of line bundles.

\end{itemize}

The paper is organized as follows. In Chapter 2 we give a detailed
description of derived categories of coherent sheaves over weighted
projective spaces and some of their noncommutative deformations. After
recalling the definition of the weighted projective space
$\mathbb{P}(\overline{\mathbf{a}})$ as a quotient
stack, we describe the category of coherent sheaves over
$\mathbb{P}(\overline{\mathbf{a}})$ and its noncommutative deformations
$\mathbb{P}_\theta(\overline{\mathbf{a}})$, and describe explicitly
generating exceptional collections for
$\db{\mathrm{coh}(\mathbb{P}_\theta(\overline{\mathbf{a}}))}$
(Theorem \ref{thm:excoll} and Corollary \ref{cor:excoll}). This is a novel
result, and we believe that it suggests a procedure that applies
to many other examples of noncommutative toric varieties.
We also discuss derived categories of
coherent sheaves over Hirzebruch surfaces.

In Chapter 3 we introduce the category of Lagrangian vanishing cycles
associated to a Lefschetz fibration, and outline the main steps involved
in its determination; to illustrate the definitions, we treat the case of
the mirror of a weighted projective line. After this warm-up, in Chapter 4
we turn to our main examples, namely the Landau-Ginzburg
models mirror to weighted projective planes and their non-exact symplectic
deformations. More precisely we start by studying the vanishing cycles
and their intersection properties, which allows us to determine all the morphisms
in $\mathrm{Lag_{vc}}$ (Lemma \ref{l:isects}). Next we study moduli spaces of
pseudo-holomorphic discs in the fiber in order to determine Floer products
(Lemmas \ref{l:m3}--\ref{l:m2}); this
gives formulas for compositions of morphisms and higher products
in $\mathrm{Lag_{vc}}$ (the latter turn out to be identically
zero). Finally, after a discussion of Maslov index and grading,
we establish an explicit correspondence between deformation parameters
on both sides (noncommutative deformation of the weighted projective plane,
and complexified K\"ahler class on the mirror)
and complete the proof of Theorem \ref{thm:main}.

Chapter 5 deals with the case of mirrors to Hirzebruch surfaces, showing
how their categories of Lagrangian vanishing cycles relate to those of
mirrors to weighted projective planes $\CP^2(n,1,1)$. In particular we prove HMS for $\mathbb{F}_n$
when $n\in\{0,1,2\}$, and show how for $n\ge 3$ a certain degenerate limit
of the Landau-Ginzburg model singles out a full subcategory of
$\mathrm{Lag_{vc}}$ whose derived category is equivalent to
that of coherent sheaves on the Hirzebruch surface.

Finally, in Chapter 6 we make various observations and concluding remarks,
related to the following directions for future research:

\begin{itemize}

\item HMS for Del Pezzo surfaces, and for higher-dimensional
weighted projective spaces (cf.\ \S \ref{ss:higherdim}
for a discussion of the case of $\CP^3$);

\item HMS for general (non toric) noncommutative deformations
(cf.\ \S \ref{ss:ncp2} for a discussion of the case of
$\CP^2$);

\item the ``other side'' of HMS -- relating derived Fukaya categories
to derived categories of B-branes on the mirror Landau-Ginzburg model.
\end{itemize}
Another topic that will be investigated in a forthcoming paper \cite{AKOS}
is HMS for products: our considerations
for $\mathbb{F}_0=\CP^1\times\CP^1$ suggest a certain
product formula on both sides of HMS~-- if we consider two manifolds
$Y_1$, $Y_2$ with mirror Landau-Ginzburg models $(X_1,W_1)$ and $(X_2,W_2)$,
then the mirror of $Y_1\times Y_2$ is simply $(X_1\times X_2, W_1+W_2)$,
and we have the following general conjecture:

\begin{conj}\label{conj:product}
$D(\FS(W_1+W_2))$ is equivalent to the product $D(\FS(W_1)\otimes \FS(W_2))$.
\end{conj}

More precisely, the vanishing cycles of $W_1+W_2$ are in one-to-one
correspondence with pairs of vanishing cycles of $W_1$ and $W_2$, and it
can be checked (cf.\ \S \ref{ss:product}) that
$$\mathrm{Hom}_{\FS(W_1+W_2)}((A_1,A_2),(B_1,B_2))\simeq
\mathrm{Hom}_{\FS(W_1)}(A_1,B_1)\otimes \mathrm{Hom}_{\FS(W_2)}(A_2,B_2).$$
The conjecture asserts that Floer products behave in the expected manner
with respect to these isomorphisms.
\smallskip

{\bf  Acknowledgements:}
We are thankful to P.\ Seidel for many helpful discussions and explanations
concerning categories of Lagrangian vanishing cycles, and to A.\ Kapustin
for explaining some features of HMS for Hirzebruch surfaces and pointing
out some references. We have also benefitted from discussions with A. Bondal,
F. Bogomolov, S. Donaldson, M.~Douglas, V. Golyshev,
M. Gromov, K. Hori, M. Kontsevich, Yu. Manin, T. Pantev, Y. Soibelman, C. Vafa, E. Witten.

Finally, we are grateful
to IPAM (especially to M. Green and H.D. Cao) for the wonderful
working conditions during the IPAM program ``Symplectic Geometry and
Physics'', where a big part of this work was done.

DA was partially supported by NSF grant DMS-0244844. LK was partially
supported by NSF grant DMS-9878353 and NSA grant H98230-04-1-0038.
DO was partially supported by the Russian Foundation for
Basic Research (grant No.~02-01-00468), Russian Presidential grant
for young scientists No.~MD-2731.2004.1,
CRDF Award No.~RM1-2405-MO-02, and the Russian Science
Support Foundation.

\section{Weighted projective spaces}\label{sec:algebra}

\subsection{Weighted projective spaces as stacks}

We start by reviewing definitions from the theory of
weighted projective spaces.

Let $\kk$ be a base field.
Let $a_0, \ldots,a_n$ be positive integers.
Define the graded algebra $S=S(a_0, \ldots,a_n)$ to be the
polynomial algebra $\kk[x_0, \ldots,x_n]$ graded by $\deg x_i
=a_i.$
Classically the projective variety $\Proj S$
is called  the weighted projective space with weights $a_0,\ldots,a_n$
and is denoted by $\bP(a_0, \ldots,
a_n).$
Consider the action of the algebraic group $\Gm=\kk^*$ on the affine space
$\bA^{n+1}$ given in some affine coordinates $x_0, \ldots,x_n$ by the formula
\begin{equation}\label{acti}
\lambda(x_0, \ldots,x_n)=(\lambda^{a_0} x_0, \ldots,\lambda^{a_n}
x_n).
\end{equation}
In
geometric terms, the weighted projective space $\bP(a_0, \ldots,
a_n)$ is the quotient variety $(\bA^{n+1}\backslash \bz)\big/\Gm$
under the induced action of the group $\Gm.$

The variety $\bP(a_0, \ldots,a_n)$ is a rational
$n$-dimensional projective variety, singular in general, whose affine
pieces $x_i\ne 0$ are isomorphic to $\bA^n \big/
\Z_{a_i}.$
For example, the variety $\bP(1, 1, n)$ is the projective cone
over a twisted rational curve of degree $n$  in $\bP^n.$

Denote by $\aa$ the vector $(a_0,\ldots,a_n)$ and write
$\bP(\aa)$ instead $\bP(a_0,\ldots,a_n)$ for brevity.

There is also another way to define the quotient of the action above:
in the category of stacks.
The quotient stack
$$
\left[(\bA^{n+1}\backslash\bz)\big/\Gm\right]
$$
will be denoted by $\PP(\aa)$ and will also be called the weighted projective space.
The stack
$\PP(\aa)$ is smooth, and from many points of view it is a more
natural object than $\bP(\aa).$

We now review the notion of stack as needed to understand
 our main example -- weighted projective spaces.
A detailed treatment of algebraic stacks can be found
in \cite{LM} and \cite{WF}.

There are two  ways of thinking about an algebraic stack:
\begin{itemize}
\item[a)]  as a category $\cX,$ with additional properties;
\item[b)]  from an atlas $R\rightrightarrows U,$ with $R$ and $U$ schemes,
$R$ determining
an equivalence relation on $U.$
\end{itemize}

 From the categorical point of view a {\sf stack} is a category $\cX$ fibered in groupoids
$p:\cX\to \Sch$ over the category
$\Sch$ of $\kk$-schemes, satisfying two descent (sheafy) properties
in \'etale topology. An {\sf algebraic stack} has to satisfy
some additional representability
conditions.
For the precise definition see \cite{LM,WF}.

Any scheme $X\in \Sch$ defines a category $\Sch/X$:
its objects are pairs $(S,\phi)$ with
$\lbrace S\stackrel{\phi}{\to} X\rbrace$ a map in $\Sch,$
and a morphism from $(S,\phi)$ to $(T,\psi)$ is a morphism $f:T\to S$ such that
$\phi f=\psi.$ The category $\Sch/X$ comes with a natural functor to $\Sch$.
Thus, any scheme is an algebraic stack.

Another example, the most important one for us,
comes from an action of an algebraic group $G$ on a scheme $X.$
The {\sf quotient stack} $\left[X/G\right]$ is defined to be the category
whose objects are those
$G$-torsors (principal homogeneous right $G$-schemes) $\G\to S$ which
are locally trivial in the \'etale topology, together with a
$G$-equivariant map from $\G$ to $X.$

In order to work with coherent sheaves on a stack it is  convenient
to use an atlas for the stack.
We describe very briefly groupoid presentations (or atlases) of algebraic stacks.
A pair of schemes $R$ and $U$ with morphisms $s,t,e,m,i,$ satisfying certain
group-like properties, is called a groupoid in $\Sch$ or an {\sf algebraic groupoid}.
For any scheme $S$ the morphisms $s,t:R\to U$ (``source'' and ``target'')
determine two maps from the set $\Hom(S,R)$
to the set $\Hom(S,U).$ A quick way to state all relations between $s,t,e,m,i$ is to
say that the induced morphisms make the ``objects'' $\Hom(S,U)$ and ``morphisms''
$\Hom(S,R)$ into a category in which all arrows are invertible.
We will denote an algebraic groupoid by $R\rightrightarrows U,$
omitting the notations for $e,m,$ and $i.$

Any scheme $X$ determines a groupoid
$X\rightrightarrows X,$ whose morphisms are identity maps.
The main example for us is the transformation groupoid associated to an
algebraic group action $X\times G\to X$, which provides an atlas
for the quotient stack $\left[ X/G\right].$
The {\sf transformation groupoid} $X\times G\rightrightarrows X$ is  defined by
$$
s(x,g)=x,
\quad
t(x,g)=x\cdot g,
\quad
m((x,g),(x\cdot g,h))=(x, g\cdot h),
\quad
e(x)=(x, e_{G}),
\quad
i(x,g)=(x\cdot g, g^{-1}).
$$

If $R\rightrightarrows U$ is an atlas for a stack $\cX,$
giving a coherent sheaf on $\cX$
is equivalent to giving a coherent sheaf $\F$ on $U,$ together with an isomorphism
$s^*\F\stackrel{\sim}{\to} t^*\F$ on $R$ satisfying a cocycle condition on
$R_t\mathop\times\limits_{U}{}_s R.$
In particular, for a quotient stack $\left[X/G\right]$
the category of coherent sheaves is equivalent to the category of
$G$-equivariant sheaves on $X$ due to effective descent for strictly flat
morphisms of algebraic stacks (see, e.g., \cite{LM}, Thm.\ 13.5.5).
Applying this fact to weighted projective spaces, we obtain that
\begin{equation}\label{eqis}
\coh(\PP(\aa))\cong \coh^{\Gm}_{\aa}(\bA^{n+1}\backslash \bz),
\end{equation}
where $\coh^{\Gm}_{\aa}(\bA^{n+1}\backslash \bz)$ is the category
of $\Gm$-equivariant coherent sheaves on $(\bA^{n+1}\backslash \bz)$
with respect
to the action given by rule (\ref{acti}).






\subsection{Coherent sheaves on weighted projective spaces}

Let $A=\bigoplus\limits_{i\ge 0} A_{i}$ be a
finitely generated graded algebra. Denote by
$\operatorname{mod}(A)$ the category of finitely generated right
$A$-modules and by $\gr(A)$ the category of finitely generated
graded right $A$-modules in which morphisms are the
homomorphisms of degree zero. Both are abelian categories.

Denote by $\tors(A)$ the full subcategory of $\gr(A)$ which
consists of those graded $A$-modules which have finite dimension
over $\kk$.

\begin{defi}
Define  the category $\qgr(A)$ to be the quotient category $
\gr(A)/\tors(A).$ The objects of $\qgr(A)$ are the objects
of the category $\gr(A)$ (we denote by $\widetilde{M}$ the object
in $\qgr(A)$ which corresponds to a module $M$). The morphisms in
$\qgr(A)$ are defined to be
$$
\Hom_{\qgr}(\widetilde{M},
\widetilde{N})=\lim_{\stackrel{\lto}{M'}} \Hom_{\gr}(M', N),
$$
where $M'$ runs over all submodules of $M$ such that $M/M'$ is
finite dimensional over $\kk.$
\end{defi}

The category $\qgr(A)$ is an abelian category and  there is a shift
functor on it: for a given graded module
$M=\mathop{\bigoplus}\limits_{i\ge 0} M_{i}$ the shifted module
$M(p)$ is defined by $M(p)_{i}= M_{p+i},$ and the induced shift
functor on the quotient category $\qgr(A)$ sends $\widetilde{M}$
to $\widetilde{M}(p)=\smash{\widetilde{M(p)}}.$

Similarly, we can consider the category $\Gr(A)$ of all graded
right $A$-modules. It contains the subcategory $\Tors(A)$ of
torsion modules. Recall that a module $M$ is called torsion if for
any element $x\in M$ one has $x A_{\ge s}=0$ for some $s,$ where
$A_{\ge s}= \mathop{\bigoplus}\limits_{i\ge s} A_{i}.$ We denote
by $\QGr(A)$ the quotient category $\Gr(A)/\Tors(A).$
 It is clear
that the intersection of the categories $\qgr(A)$ and $\Tors(A)$
in the category $\QGr(A)$ coincides with $\tors(A).$ In
particular, the category $\QGr(A)$ contains $\qgr(A)$ as a full
subcategory. Sometimes it is convenient to work with $\QGr(A)$
instead of $\qgr(A).$

In the case when the algebra  $A=\mathop{\bigoplus}\limits_{i\ge
0} A_{i}$ is a commutative graded algebra generated over $\kk$ by
its first component (which is assumed to be finite dimensional)
J.-P.~Serre \cite{Se} proved that the category of coherent sheaves $\coh(X)$
on the projective variety $X=\Proj A$ is equivalent to the
category $\qgr(A)$. Such an equivalence also holds for
the category of quasicoherent sheaves on $X$ and the category
$\QGr(A)=\Gr(A)/\Tors(A).$

This theorem can be extended to general finitely generated
commutative algebras if we work at the level of quotient stacks.

Let $ S=\mathop\bigoplus\limits_{p=0}^{\infty} S_p $ be a commutative
graded $\kk$-algebra which is connected, i.e. $S_0=\kk.$ The
grading on $S$ induces an action of the group $\Gm$ on the affine
scheme $\Spec S.$ Let $\bz$ be the closed point of $\Spec S$ that
corresponds to the ideal $S_+=S_{\ge 1}\subset S.$ This point is invariant
under the action.

\begin{defi}
 Denote by $\PROJ S$ the quotient stack $\left[(\Spec S\backslash
\bz)\big/\Gm\right].$
\end{defi}

There is a natural map $\PROJ S \to \Proj S,$ which is an isomorphism when
 the
algebra $S$ is generated by its first component $S_1.$

\begin{prop}
 Let $S=\mathop{\oplus}\limits_{i\ge 0} S_{i}$ be a
graded finitely generated algebra. Then the category of
(quasi)coherent sheaves on the quotient stack $\PROJ(S)$ is
equivalent to the quotient category $\qgr(S)$ $($resp.\ $\QGr(S))$.
\end{prop}
\begin{proof}
Let $\bz$ be the closed point on the affine scheme $\Spec S$ which
corresponds to the maximal ideal $S_{+}\subset S.$ Denote by $U$
the scheme $(\Spec S\backslash\bz).$ We know that the category of
(quasi)coherent sheaves on the stack $\PROJ S$ is equivalent to
the category of $\Gm$-equivariant (quasi)coherent sheaves on
$U.$ The category of (quasi)coherent sheaves on $U$ is equivalent
to the quotient of the category of (quasi)coherent sheaves on
$\Spec S$ by the subcategory of (quasi)coherent sheaves with
support on $\bz.$ This is also true for the categories of
$\Gm$-equivariant sheaves. But the category of (quasi)coherent
$\Gm$-equivariant sheaves on $\Spec S$ is just the
category $\gr(S)$ (resp.\ $\Gr(S)$) of graded modules over $S,$ and the
subcategory of (quasi)coherent sheaves with support on $\bz$
coincides with the subcategory $\tors(S)$ (resp.\ $\Tors(S)$). Thus, we
obtain that $\coh(\PROJ S)$ is equivalent to the quotient category
$\qgr(S)=\gr(S)/\tors(S)$ (and $\Qcoh(\PROJ S)$ is equivalent to
$\QGr(S)=\Gr(S)/\Tors(S)$).
\end{proof}

\begin{cor}\label{coreq}
The category of (quasi)coherent sheaves on the weighted projective
space $\PP(\aa)$ is equivalent to the category
$\qgr(S(a_0,\ldots,a_n))$ $($resp.\ $\QGr(S(a_0,\ldots,a_n)))$.
\end{cor}

We conclude this section by giving the definition of noncommutative
weighted projective
spaces and the categories of coherent sheaves on them.
Consider a  matrix
 $\theta=(\theta_{ij})$ of dimension $(n+1)\times (n+1)$ with entries
$\theta_{ij}\in\kk^*$ for all $i,j.$ The set of all such matrices
will be denoted by $\Mat.$
Consider the graded algebra
$S_{\theta}=S_{\theta}(a_0,\ldots,a_n)$  generated
by elements $x_{i}, i=0,\ldots,n$ of degree $a_i$ and with  relations
$$
\theta_{ij}x_i x_j= \theta_{ji}x_j x_i
$$
for all $i$ and $j.$
This algebra is a noncommutative deformation of the algebra
$S(a_0,\ldots,a_n).$
It can be easily checked that the algebra $S_{\theta}$ depends
only on the matrix $\theta^\mathrm{an},$ with entries
\begin{equation}\label{antis}
\theta^\mathrm{an}_{ij}:=\theta_{ij}\theta_{ji}^{-1}
\qquad
\text{for all}
\quad 0\le i,j\le n.
\end{equation}
Thus, if $(\theta')^\mathrm{an}=\theta^\mathrm{an}$
for two matrices $\theta'$ and $\theta,$
then $S_{\theta'}\cong S_{\theta}.$

As before, denote by $\qgr(S_{\theta})$ the
quotient category  $\gr(S_{\theta})/\tors(S_{\theta}),$ where
$\gr(S_{\theta})$ is the category of finitely generated graded
right $S_{\theta}$-modules and $\tors(A)$ is the full
subcategory of $\gr(S_{\theta})$ consisting of graded modules of
finite dimension over $\kk$.

Corollary \ref{coreq} suggests that the category
$\qgr(S_{\theta})$ can be considered as  the category of coherent
sheaves on a noncommutative weighted projective space.
We  will
denote this space by $\PP_{\theta}(\aa)$ and will write
$\coh(\PP_{\theta})$ instead $\qgr(S_{\theta}).$
Similarly, the category of quasi-coherent sheaves
$\Qcoh(\PP_{\theta})$ is defined as the quotient
$\QGr(S_{\theta})=\Gr(S_{\theta})/\Tors(S_{\theta}).$

\subsection{Cohomological properties of coherent sheaves on
$\PP_{\theta}(\aa)$}

In this section we discuss properties of categories of coherent
sheaves on the noncommutative weighted projective spaces
$\PP_{\theta}(\aa).$ Note that the usual commutative
weighted projective space is a particular case of the noncommutative
one, when $\theta$ is the matrix with all entries equal to 1.

All algebras $S_{\theta}(a_0,\ldots,a_n)$ are noetherian.
This follows from the fact that they are Ore extensions of
commutative polynomial algebras (see for example \cite{MR}). For
the same reason the algebras $S_{\theta}(a_0,\ldots,a_n)$ have
finite right (and left) global dimension, which is equal to
$(n+1)$ (see \cite{MR}, p.~273).
Recall that the global dimension of a ring $A$ is the minimal
number $d$ (if it exists) such that for any two modules $M$ and $N$
we have $\Ext^{d+1}_{A}(M, N)=0.$

The notion of a regular algebra was
introduced in \cite{AS}. As we will see below,  regular algebras
have many good properties. More details can be found in \cite{Ar}.

\begin{defi}\label{reg}
 A graded algebra $A$ is called regular of dimension $d$
if it satisfies the following conditions:

\begin{tabular}{ll}
{\rm(1)}& $A$ has global dimension $d,$\\
{\rm(2)}& $A$ has polynomial growth, i.e. $\dim A_p\le
cp^{\delta}$ for some
$c, \delta\in \R,$\\
{\rm(3)}& $A$ is Gorenstein, meaning that $\Ext^{i}_{A}(\kk, A)=0$
if $i\ne d,$
and $\Ext^{d}_{A}(\kk, A)=\kk(l)$\\
& for some $l.$ The number $l$ is called
the Gorenstein parameter.
\end{tabular}
\end{defi}
Here $\Ext_{A}$ stands for the $\Ext$ functor in the category
of right modules $\operatorname{mod}(A).$

\begin{prop} The algebra $S_{\theta}(a_0,\ldots,a_n)$
is a noetherian regular
algebra of global dimension $n+1.$ The Gorenstein parameter $l$
of this algebra is equal to the sum $\mathop\sum\limits_{i=0}^{n}
a_i.$
\end{prop}
\begin{proof}
Property (1) holds, as for all Ore extensions of commutative
polynomial algebras. Property (2) holds because our algebras have the
same growth
as ordinary polynomial algebras. Property (3) follows from the
following Koszul resolution of the right module $\kk_{S_{\theta}}$
\begin{multline}\label{Kos}
0\to S_{\theta}(-\mathop\sum\limits_{i=0}^{n} a_i)
\to\bigoplus_{i_0<\ldots<i_{n-1}}S_{\theta}(-\mathop\sum\limits_{j=0}^{n-1}
a_{i_j}) \to
\cdots\\
\cdots\to
\bigoplus_{i_0<i_1}S_{\theta}(-a_{i_0}-a_{i_1})\to
\bigoplus_{i=0}^n S_{\theta}(-a_i)\to S_{\theta}\to
\kk_{S_{\theta}}\to 0,
\end{multline}
and the fact that the transposed complex is a resolution of the left module
${}_{S_{\theta}}\kk,$ shifted to the degree $l=\sum a_i.$
The explicit formula for the differentials in the complex (\ref{Kos})
will be given later (see (\ref{Kos2})).
\end{proof}

Denote by $\O(i)$  the object $\widetilde{S_{\theta}(i)}$ in
the category $\coh(\PP_{\theta})=
\qgr(S_{\theta}).$
Consider the sequence $\lbrace \O(i) \rbrace_{i\in \Z}.$
It can be checked that  the following properties hold true:

\begin{itemize}
\item[(a)] For any coherent sheaf
$\F$ there are integers $k_{1},\ldots,k_{s}$ and
an epimorphism
$$\mathop{\op}\limits_{i=1}^{s} \O(-k_{i})\to \F.$$

\item[(b)] For every epimorphism $\F\to \G$ the induced map
$\Hom(\O(-n), \F)\to \Hom(\O(-n), \G)$ is surjective for $n\gg 0.$
\end{itemize}

A sequence which satisfies such conditions will be called {\sf ample}.
It is proved in \cite{Ar} that the sequence $\{ \O(i) \}$ is ample
in $\qgr(A)$ for any graded right noetherian $\kk$-algebra $A$ if it
satisfies the extra condition:
$$
(\chi_{1})  : \quad \dim_{\kk}\Ext^{1}_{A}(\kk, M) < \infty
$$
for any finitely generated graded $A$-module $M.$

This condition can be verified for all noetherian regular algebras
(see \cite{Ar}, Theorem 8.1). In particular, the sequence
$\lbrace \O(i) \rbrace_{i\in \Z}$
in the category
$\coh(\PP_{\theta})$ is  ample.

For any sheaf $\F\in \qgr(A)$ we can define a graded module
$\Gamma(\F)$ by the rule:
$$
\Gamma(\F):= \mathop{\oplus}\limits_{i\ge 0} \Hom(\O(-i), \F)
$$
It is proved in \cite{Ar} that for any noetherian algebra $A$ that
satisfies the condition $(\chi_1)$ the correspondence $\Gamma$ is a
functor from $\qgr(A)$ to $\gr(A)$ and the composition of $\Gamma$
with the natural projection $\pi : \gr(A)\lto \qgr(A)$ is isomorphic
to the identity functor (see \cite{Ar}, \S~3,4).

We formulate next a result about the cohomology of sheaves on
noncommutative weighted projective spaces. This result is proved in
\cite{Ar} (Theorem 8.1) for a general regular algebra and parallels the
commutative case.

\begin{prop}\label{cohom}
 Let $S_{\theta}=S_{\theta}(a_0,\ldots,a_n)$ be the algebra of the noncommutative
weighted projective space $\PP_{\theta}=\PP_{\theta}(\aa).$
 Then
\begin{itemize}
\item[{\rm 1)}] The cohomological dimension of the category
$\coh(\PP_{\theta}(\aa))$ is equal to
$n,$ i.e. for any two coherent sheaves $\F, \G\in \coh(\PP_{\theta})$ the space
$\Ext^{i}(\F, \G)$ vanishes if $i>n.$

\item[{\rm 2)}] There are isomorphisms
\begin{equation}\label{cohinv}
H^{p}(\PP_{\theta}, \O(k))=
    \begin{cases}
    (S_{\theta})_{k}& \text{for $p=0,\; k\ge 0$}\\
    (S_{\theta})_{-k-l}^*& \text{ for $p=n,\; k\le -l$}\\
    0& \text{otherwise}
    \end{cases}
\end{equation}
\end{itemize}
\end{prop}

This proposition and the ampleness of the sequence $\{ \O(i) \}$
imply the following corollary.

\begin{cor}\label{zero}
For any sheaf $\F\in
\coh(\PP_{\theta})$ and for all sufficiently large $i\gg 0$ we have
$
H^{k}(\PP_{\theta}, \F(i))=0
$
for all $k>0.$
\end{cor}
\begin{proof}
The group $H^k(\PP_{\theta}, \F(i))$ coincides with $\Ext^{k}(\O(-i), \F).$
Let $k$ be the maximal integer (it exists because the global
dimension is finite) such that for some $\F$ there exists
arbitrarily large $i$ such that $\Ext^{k}(\O(-i), \F)\ne0.$ Assume
that $k\ge 1.$ Choose an epimorphism
$\mathop{\op}\limits_{j=1}^{s}\O(-k_{j}) \to \F.$ Let $\F_1$
denote its kernel. Then for $i>\max\{k_j\}$ we have
$\Ext^{>0}(\O(-i),\mathop{\op}\limits_{j=1}^{s}\O(-k_{j}))=0,$
hence $\Ext^{k}(\O(-i), \F)\ne0$ implies $\Ext^{k+1}(\O(-i),
\F_1)\ne0.$ This contradicts the assumption of the maximality of
$k.$
\end{proof}


One of the  useful properties of commutative smooth projective varieties
is the existence of the dualizing sheaf. Recall that a sheaf
$\omega_X$ is called dualizing if for any $\F\in \coh(X)$ there are
natural isomorphisms of $\kk$-vector spaces
$$
H^i(X, \F)\cong \Ext^{n-i}(\F, \omega_X)^*,
$$
where $*$ denotes the $\kk$--dual space. The Serre duality theorem
asserts the existence of a dualizing sheaf for smooth projective
varieties.  In this case the dualizing sheaf is a line bundle and
coincides with the sheaf of differential forms $\Omega_{X}^n$ of
top degree.

Since the definition of $\omega_X$ is given in abstract categorical
terms, it can be extended to the noncommutative case as well. More
precisely, we will say that $\qgr(A)$ satisfies classical {\it Serre
duality} if there is an object $\omega\in \qgr(A)$ together with
natural isomorphisms
$$
\Ext^{i}(\O,-)\cong \Ext^{n-i}(-, \omega)^*.
$$

Our noncommutative varieties $\PP_{\theta}(\aa)$ satisfy classical
Serre duality, with dualizing sheaves being $\O(-l),$
where $l=\sum a_i$ is the Gorenstein
parameter for $S_{\theta}(a_0,\ldots,a_n).$ This follows from the paper
\cite{YZ}, where the existence of a dualizing sheaf in $\qgr(A)$
has been proved for a  class of algebras which includes all
noetherian regular algebras. In addition, the authors of \cite{YZ}
showed that the dualizing sheaf coincides with $\wt{A}(-l),$ where
$l$ is the Gorenstein parameter for $A.$

There is a reformulation of Serre duality
in terms of bounded derived categories \cite{BK}.
A {\it Serre
functor} in the bounded derived category $\bD^b(\coh(\PP_{\theta}))$
is by definition an exact auto\-equi\-valence $S$ of $\bD^b(\coh(\PP_{\theta}))$
such that for any objects $X,Y\in \bD^b(\coh(\PP_{\theta}))$ there is
a bifunctorial isomorphism
$$
\Hom(X, Y)\stackrel{\sim}{\lto}\Hom(Y,SX)^*.
$$
Serre duality can be reinterpreted as the existence
of a Serre functor in the bounded derived category.


\subsection{Exceptional collection on $\PP_{\theta}(\aa)$}\label{ss:excoll}

For many reasons it is more natural to work not with the abelian category
of coherent sheaves
but with its bounded derived category $\db{\coh(\PP_{\theta})}.$
The purpose of this section is to describe the bounded derived category
of coherent sheaves
on the noncommutative weighted projective spaces
in the terms of exceptional collections.

First, we briefly recall  the definition of the bounded derived category for
 an abelian category $\A.$
We start with
 the category $\bC^b(\A)$ of bounded differential complexes
$$
M^{\bdot}=(0\lto\cdots\lto M^p\stackrel{d^p}{\lto}
M^{p+1}\stackrel{d^{p+1}}{\lto} M^{p+2}\lto\cdots\lto 0),
\quad M^p\in\A, \quad p\in\Z, \quad
d^2=0.
$$
A morphism of complexes $f: M^{\bdot}\lto N^{\bdot}$ is called
null-homotopic if $f^p=d_N h^p + h^{p+1} d_M$ for all $p\in\Z$
and some family of morphisms $h^{p}: M^{p}\lto N^{p-1}.$ Now the
{\it homotopy category} $\bH^b(\A)$ is defined as a category with the same
objects as $\bC^b(\A),$ whereas morphisms in $\bH^b(\A)$
are equivalence classes $\overline{f}$
of morphisms of complexes modulo null-homotopic
morphisms. A morphism of complexes $s: N^{\bdot}\to
M^{\bdot}$ is called a {\it quasi-isomorphism} if the induced morphisms $H^p s:
H^p(N^{\bdot})\to H^p(M^{\bdot})$ are isomorphisms for all
$p\in\Z.$ Denote by $\Sigma$ the class of all
quasi-isomorphisms.
The bounded  derived category $\bD^b(\A)$ is now defined as the
{\sf localization} of $\bH^b(\A)$ with respect to the class $\Sigma$ of all
quasi-isomorphisms. This means that the derived category has the
same objects as the homotopy category $\bH^b(\A),$ and that
morphisms in the derived category
are given by left fractions $s^{-1}\circ f$ with $s\in\Sigma.$

\begin{remark}\label{rem1}
For any full subcategory
$\E\subset\A$ one can construct the homotopy category $\bH^b(\E)$ and a functor
$\bH^b(\E)\to\bD^b(\A).$ In some cases, for example when $\A$ is the abelian
category of modules over an algebra $A$ of finite global dimension
and $\E$ is the subcategory of projective modules, this functor
$\bH^b(\E)\to\bD^b(\A)$ is an equivalence of triangulated categories.
\end{remark}

Second, we recall the notion of an exceptional collection.
\begin{defi}
An object $E$ of a $\kk$-linear triangulated category $\D$ is
said to be {\sf exceptional} if
$\Hom(E, E[k])=0$ for all $k\ne 0,$ and $\Hom(E, E)=\kk.$

An ordered set of exceptional objects $\sigma=\left(E_0,\ldots E_n\right)$
is called {\sf an exceptional collection} if
$\Hom(E_j, E_i[k])=0$ for $j>i$ and all $k.$
The exceptional collection $\sigma$ is called {\sf strong} if  it satisfies
the additional condition $\Hom(E_j, E_i[k])=0$ for all $i,j$ and for $k\ne 0.$
\end{defi}
\begin{defi} An exceptional collection
$\left(E_0,\ldots,E_n\right)$ in a category $\D$ is called {\sf full}
if it generates the category
$\D,$ i.e. the minimal triangulated subcategory of $\D$
containing all objects $E_i$ coincides
with $\D.$ We write in this case
$$
\D=\left\langle E_0,\ldots,E_n\right\rangle.
$$
\end{defi}

Consider the bounded derived category of coherent sheaves
$\db{\coh(\PP_{\theta})}.$
We prove that this category has an exceptional collection which is strong
and full.
In this case we will say that the noncommutative weighted projective space
$\PP_{\theta}$ possesses
a full strong exceptional collection.

\begin{theo}\label{thm:excoll}
For any noncommutative weighted projective space
$\PP_{\theta}(\aa)$ and for any $k\in\Z$ the ordered set
$\sigma(k)=\left(\O(k),\ldots,O(k+l-1)\right),$
where $l=\sum a_i$ is the Gorenstein parameter of $S_{\theta},$
forms a full strong exceptional collection
in the category $\db{\coh(\PP_{\theta})}.$
\end{theo}
\begin{proof}
It  follows directly from Proposition \ref{cohom}
that the collection $\sigma(k)$ is exceptional and strong.
To prove that the collection is full let us consider the
triangulated subcategory
$\D\subset \db{\coh(\PP_{\theta})}$ generated by the collection $\sigma(k).$
The exact sequence (\ref{Kos}) induces the exact sequence
\begin{multline}\label{Kos1}
0\to \O(-\mathop\sum\limits_{i=0}^{n} a_i)
\to\bigoplus_{i_0<\ldots<i_{n-1}}\O(-\mathop\sum\limits_{j=0}^{n-1}
a_{i_j}) \to
\cdots\\
\cdots\to
\bigoplus_{i_0<i_1}\O(-a_{i_0}-a_{i_1})\to
\bigoplus_{i=0}^n \O(-a_i)\to \O
\to 0.
\end{multline}
Shifting it by $k+l$ one obtains that the object $\O(k+l)$ also belongs to
$\D$ and repeating
this procedure deduce that  $\O(i)$ for all $i$ belongs to $\D.$
Assume that $\D$ does not coincide with $\db{\coh(\PP_{\theta})}$
and take an object
$U$ which  does not belong to $\D.$
It is proved in \cite{Bo} (Theorem 3.2) that the subcategory $\D$ is
{\it admissible}, i.e.\
the natural embedding functor
$\D\hookrightarrow\db{\coh(\PP_{\theta})}$
has right and left adjoint functors. Denote by $j$ the right adjoint
and complete
the canonical map $jU\lto U$ to
a distinguished triangle
$$
jU\lto U\lto C\lto jU[1].
$$
It follows from adjointness that for any object $V\in \D$ the space
$\Hom(V,C)$ vanishes. The object $C$ is a bounded complex of coherent sheaves.
Denote by $H^k(C)$ the leftmost nontrivial cohomology of the complex $C.$
The ampleness of the sequence $\lbrace \O(i) \rbrace_{i\in \Z}$
guarantees that for sufficiently large $i$ the space
$\Hom(\O(-i),H^k(C))$ is nontrivial.
This implies that $\Hom(\O(-i)[-k], C)$ is nontrivial, which contradicts
to the fact
that the object $\O(-i)[-k]$ belongs to $\D.$
\end{proof}

The strong exceptional collection on the ordinary projective space
$\PP^n$ was constructed by Beilinson in \cite{Be}.
This question for the weighted projective spaces was considered
in \cite{Ba}.

\begin{defi}
The {\sf algebra of the strong exceptional collection}
$\left(E_0,\ldots,E_n\right)$
is the algebra of endomorphisms of the object
$\mathop\oplus\limits_{i=0}^{n}E_i.$
 Denote by $\Te$ the
sheaf $\mathop\oplus\limits_{i=0}^{l-1}\O(i)$
and by $\Be$ the algebra of the collection
$\left(\O,\ldots,\O(l-1)\right)$
on the noncommutative weighted projective space $\PP_{\theta},$ i.e.\ $\Be=\End(\Te).$
\end{defi}
The algebra $\Be$ is a finite dimensional algebra over $\kk.$
Denote by $\mbox{mod--}\Be$ the category of finitely generated right modules
over $\Be.$ For any coherent sheaf $\F\in \coh(\PP_{\theta})$
the space $\Hom(\Te,\F)$ has a structure of right
$\Be$-module. Denote by $P_i$ the modules
$\Hom(\Te, \O(i))$ for $i=0,\ldots,(l-1).$
All these are projective $\Be$-modules and
$\Be=\mathop\oplus\limits_{i=0}^{l-1} P_i.$
The algebra $\Be$ has $l$ primitive idempotents $e_i, i=0,\ldots,l-1$ such that
$1_\Be=e_0+\cdots+e_{l-1}$ and $e_i e_j=0$ if $i\ne j.$ The right projective modules
$P_i$ coincide with $e_i \Be.$ The morphisms between them
can be easily described since
$$
\Hom(P_i, P_j)=\Hom(e_i \Be, e_j \Be)\cong e_j \Be e_i\cong
\Hom(\O(i), \O(j))=(S_{\theta})_{j-i}.
$$

Moreover, the algebra $\Be$ has  finite global dimension.
This follows from the fact that any
right (and left) module $M$ has a finite projective resolution
consisting of the projective
modules $P_i.$
Indeed the map
$$
\bigoplus_{i=0}^{l-1}\Hom(P_i,M)\otimes P_i\lto M
$$
is surjective  and there are no non-trivial homomorphisms from $P_{l-1}$ to
the kernel of this map.
Iterating this procedure  we get a finite resolution of $M.$

Sometimes it is useful to represent the algebra $\Be$ as a category $\Bc$
which has
$l$ objects, say $v_0,\ldots,v_{l-1},$ and morphisms  defined by
$$
\Hom(v_i, v_j)\cong\Hom(\O(i), \O(j))\cong (S_{\theta})_{j-i}
$$
with the natural composition law.
Thus $\Be=\mathop\bigoplus\limits_{0\le i,j\le l-1}\Hom(v_i, v_j).$

The algebra $\Be$ is a basis algebra.
This means that the quotient of $\Be$ by the radical $\rad(\Be)$
is isomorphic to the direct sum of $l$ copies of the field $\kk.$
The category $\mbox{mod--}\Be$ has $l$ irreducible modules which will be denoted
$Q_i, i=0,\ldots,l-1$, and $\mathop\oplus\limits_{i=0}^{l-1}
Q_i=\Be/\rad(\Be).$
 The modules $Q_i$ are chosen so that $\Hom(P_i, Q_j)\cong \delta_{i,j}\, \kk.$

\bigskip

Our next topic is the notion of {\sf mutation} in an exceptional collection.
Let $\sigma=(E_0,\dots,E_n)$ be an exceptional collection in a
triangulated category $\D.$ Consider a pair $(E_i, E_{i+1})$ and the canonical
maps
$$
\Hom^{\bdot}(E_i, E_{i+1})\otimes E_i\lto E_{i+1}
\qquad
\text{and}
\quad
E_i\lto \Hom^{\bdot}(E_i, E_{i+1})^*\otimes  E_{i+1},
$$
where by definition
\begin{align*}
&\Hom^{\bdot}(E_i, E_{i+1})\otimes E_i=\bigoplus_{k\in \Z}\Hom^{k}(E_i, E_{i+1})\otimes E_i[-k],\\
&\Hom^{\bdot}(E_i, E_{i+1})^*\otimes E_{i+1}=\bigoplus_{k\in \Z}\Hom^{-k}(E_i, E_{i+1})\otimes E_{i+1}[-k]
\end{align*}
(recall that the tensor product of a vector space $V$ with an object $X$ may be considered as the direct sum
of $\dim V$ copies of the object $X$).

We define objects $LE_{i+1}$ and  $RE_i$
as the objects obtained from the distinguished triangles
\begin{align*}
&LE_{i+1}\lto \Hom^{\bdot}(E_i, E_{i+1})\otimes E_i\lto E_{i+1},\\
&E_i\lto \Hom^{\bdot}(E_i, E_{i+1})^*\otimes  E_{i+1}\lto RE_i .
\end{align*}
The object $LE_{i+1}$ (resp.\ $RE_i$) is called by {\sf left (right) mutation} of
$E_{i+1}$ (resp.\ $E_i$)
in the collection $\sigma.$
It can be checked that the objects $LE_{i+1}$ and  $RE_i$ are exceptional and, moreover,
the two  collections
\begin{align*}
L_i\sigma=\left( E_0,\dots, E_{i-1}, LE_{i+1}, E_i, E_{i+2},\dots, E_n \right)\\
R_i\sigma=\left( E_0,\dots, E_{i-1}, E_{i+1}, RE_{i}, E_{i+2},\dots, E_n \right)
\end{align*}
are exceptional as well. These collections are called left and right mutations
of the collection $\sigma$ in the pair $(E_i,E_{i+1}).$
Consider $R_i$ and $L_i$ as operations on the set of all exceptional collections in the category $\D$.
It is easy to see that they are mutually inverse, i.e. $R_i L_i=1.$
Moreover, $L_i$ (resp. $R_i$) satisfy the Artin braid group relations:
$$
L_i L_{i+1} L_i=L_{i+1} L_i L_{i+1},
\qquad
R_i R_{i+1}R_i=R_{i+1}R_i R_{i+1}
$$
(see \cite{Bo,Go}).

Denote by $L^{(k)} E_{i}$ with  $k\le i$ the result of multiple left
mutations of the object $E_i$
in the collection $\sigma.$ Analogously for right mutations.

\begin{defi}
The exceptional collection $(L^{(n)} E_n, L^{(n-1)} E_{n-1},\dots E_0)$ is called the
{\sf left dual collection} for the collection $(E_0,\dots, E_n).$ Analogously, the
{\sf right dual collection} is defined as $(E_n,RE_{n-1},\dots, R^{(n)} E_0).$
\end{defi}

\begin{example}\label{ex:dual}
For example, let us consider the full exceptional collection
$(P_0,\dots, P_{l-1})$ in the category $\db{\mod-\Be},$ consisting
of the projective $\Be$-modules $P_i.$ It can be shown (e.g.
\cite{Bo}, Lemma 5.6) that the irreducible modules $Q_i, 0\le i<l$
can be expressed as
$$
Q_i\cong L^{(i)} P_i[i].
$$
Thus, the left dual for the exceptional collection
$(P_0,\dots, P_{l-1})$ coincides with the collection
$(\nobreak{Q_{l-1}[1-l]}, \dots, Q_0).$
\end{example}

\subsection{A description of the derived categories of coherent sheaves
on $\PP_{\theta}(\aa)$}

The natural isomorphisms
$\Hom(P_i, P_j)\cong \Hom(\O(i), \O(j)),$ which are direct consequences of
the construction of the algebra $\Be,$
allow us to construct a functor
$\bar{F}:\bH^b(\P)\lto \db{\coh(\PP_{\theta})},$
where $\P$ is the full subcategory
of the category of right modules $\mbox{mod--}\Be$ consisting of
finite direct sums of the projective modules
$P_i, i=0,\ldots,l-1.$
The functor $\bar{F}$ sends $P_i$ to $\O(i)$ and any bounded complex
of projective modules to the corresponding complex of $\O(i), i=0,\ldots,l-1.$
It  follows from Remark \ref{rem1} that the functor $\bar{F}$ induces
a functor
$$
F:\db{\mbox{mod--}\Be}\lto\db{\coh(\PP_{\theta})}.
$$

\begin{theo}
The functor $F:\db{\mbox{\rm mod--}\Be}\lto\db{\coh(\PP_{\theta})}$
is an equivalence of the derived categories.
\end{theo}
Since the exceptional collection $\left(\O,\ldots,\O(l-1)\right)$
generates the category $\db{\coh(\PP_{\theta})}$ it is sufficient to
check that the functor $F$ is fully faithful.
We know that for any $0\le i,j\le l-1$ and any $k$ there are isomorphisms
$$
\Hom(P_i, P_j[k])\stackrel{\sim}{\lto}\Hom(FP_i,FP_j[k])=\Hom(\O(i), \O(j)[k]).
$$
Since $P_i, i=0,\ldots,l-1$ generate $\db{\mbox{mod--}\Be},$
the proof of the theorem is a consequence of the following lemma.
\begin{lemma}\label{ff}
Let $\A$ be abelian category and $\D$ be a triangulated category.
Let $F: \db{\A}\lto\D$ be an exact functor and let
$\lbrace E_i \rbrace_{i\in I}$ be a set
of objects of $\db{\A}$ which generates $\db{\A}$
(i.e. the minimal full triangulated
subcategory of $\db{\A}$ containing all $E_i$ coincides with $\db{\A}$).
Assume that the maps
$$
\Hom(E_i, E_j[k])\lto \Hom(FE_i, FE_j[k])
$$
are isomorphisms for all $i,j\in I$ and any $k\in \Z.$
Then the functor $F$ is fully faithful.
\end{lemma}
\begin{proof}
This lemma is known and results from d\'evissage (e.g. \cite{Ha},10.10, \cite{Ke}4.2).
We first consider the full subcategory $\cC\in\db{\A}$
which consists of all objects $X$ such that the maps
$$
\Hom(X,E_i[k])\stackrel{\sim}{\lto}\Hom(FX, FE_i[k])
$$
are isomorphisms for all $i\in I$ and all $k\in \Z.$
The category $\cC$ is a triangulated subcategory,
because it is closed with respect to the translation functor and,
for any distinguished triangle
$$
X\lto Y\lto Z\lto X[1],
$$
if $X$ and $Y$ belong to $\cC,$ then $Z$  belongs too.
The last statement is a consequence of the five lemma, i.e., since the morphisms
$f_1,f_2,f_4,f_5$ in the diagram
$$
\begin{CD}
\Hom(Y[1], E_i)@>>>&\Hom(X[1], E_i)@>>>&\Hom(Z, E_i)@>>>\\
@Vf_1VV&@Vf_2VV&@VVf_3V\\
\Hom(FY[1], FE_i)@>>>&\Hom(FX[1], FE_i)@>>>&\Hom(FZ, FE_i)@>>>\\
&&&@>>>&\Hom(Y, E_i)@>>>&\Hom(X, E_i)\\
&&&&&&@VVf_4V&@VVf_5V\\
&&&@>>>&\Hom(FY, FE_i)@>>>&\Hom(FX, FE_i)
\end{CD}
$$
are isomorphisms, the morphism $f_3$ is an isomorphism too.
The subcategory $\cC$ contains the objects $E_i$ and, hence,
coincides with $\db{\A}.$
Now consider the full
subcategory $\B\subset\db{\A}$  consisting of all objects $X$ such that the map
$$
\Hom(Y,X[k])\stackrel{\sim}{\lto}\Hom(FY, FX[k])
$$
is an isomorphism for every object $Y\in\db{\A}$ and all $k\in \Z.$
By the same argument as above the subcategory $\B$ is triangulated and
contains all $E_i.$
Therefore, it coincides with $\db{\A}.$ This proves the lemma and completes
the proof of the Theorem.
\end{proof}

There is also a right adjoint to $F,$ namely a functor
$G:\db{\coh(\PP_{\theta})}\lto \db{\mod-\Be}.$
To construct it we have to consider the functor
$$
\Hom(\Te, -):\Qcoh(\PP_{\theta})\lto \Mod-\Be
$$
where $\Mod-\Be$ is the category of all right modules over $\Be.$
Since $\Qcoh(\PP_{\theta})$ has enough injectives and has finite global dimension
there is a right derived functor
$$
\bR\Hom(\Te, -):\db{\Qcoh(\PP_{\theta})}\lto \db{\Mod-\Be}.
$$
$\db{\coh(\PP_{\theta})}$ is equivalent
to the full subcategory $\bD^b_{\coh}(\Qcoh(\PP_{\theta}))$ of
$\db{\Qcoh(\PP_{\theta})}$
whose objects are complexes with cohomologies in $\coh(\PP_{\theta}).$
Moreover, the functor $\bR\Hom(\Te, -)$ sends an object of
$\bD^b_{\coh}(\Qcoh(\PP_{\theta}))$ to an object of the subcategory
$\bD^b_{\mod}(\Mod-\Be),$ which is also equivalent to $\db{\mod-\Be}.$
This gives us a functor
$$
G=\bR\Hom(\Te, -):\db{\coh(\PP_{\theta})}\lto \db{\mod-\Be}.
$$
The functor $G$ is right adjoint to $F,$ and it is an equivalence of
categories as well.

In the end of this paragraph we describe an equivalence relation
$\theta\sim\theta'$
 on the space of
all matrices $\theta$ with $\theta_{ij}\in\kk^*$ for all $i,j$
under which the noncommutative weighted projective spaces
$\PP_{\theta}$ and $\PP_{\theta'}$ have equivalent
abelian categories of coherent sheaves.
It was mentioned above that the graded algebras
$S_{\theta}$ depend only on the matrix
$\theta^{an}$ defined by the rule (\ref{antis}).
However, it can also happen that two different
algebras $S_{\theta}$ and $S_{\theta'}$ produce isomorphic algebras
$B_{\theta}$ and $B_{\theta'}.$

\begin{prop}\label{equiv}
Let $(m_0,\ldots,m_{n})\in (\kk^*)^{(n+1)}$ be any vector with non-zero
entries. Suppose that two matrices $\theta,\theta' \in \Mat$
are related by the formula
\begin{equation}\label{relth}
\theta'_{ij}=\theta_{ij}\cdot m_i^{a_j}.
\end{equation}
Then the algebras $B_{\theta'}$ and $B_{\theta}$ are isomorphic.
\end{prop}
\begin{proof}
Consider the category ${\mathfrak B}_{\theta'}$
and its autoequivalence $\tau$
which acts by identity on the objects and acts on the spaces
$\Hom(v_i, v_j)$ as the multiplication by $(m_i)^{(j-i)}.$
There is a natural basis of the spaces $\Hom(v_i, v_j)$ which is induced
by the monomial basis $x_{i_0}\cdots x_{i_k},$ $0\le i_0\le\cdots\le i_k\le n$
of $S_{\theta'}.$
The transformation of this basis under the equivalence
$\tau$ gives us a new basis
in which the category ${\mathfrak B}_{\theta'}$ coincides with the category
${\mathfrak B}_{\theta}$ equipped with its natural basis coming from
the monomial basis of $S_{\theta}.$
The equivalence of the categories ${\mathfrak B}_{\theta'}$
and ${\mathfrak B}_{\theta}$ implies an isomorphism of the algebras
$B_{\theta'}$
and $B_{\theta}.$
\end{proof}

If now the algebras $B_{\theta'}$ and $B_{\theta}$ are isomorphic,
then the composition of the functors
$$
\bD^b(\coh(\PP_{\theta'}))\stackrel{G_{\theta'}}{\lto}
\bD^b(\mod-B_{\theta'})\cong\bD^b(\mod-B_{\theta})\stackrel{F_{\theta}}{\lto}
\bD^b(\coh(\PP_{\theta}))
$$
is an equivalence of derived categories.
This equivalence evidently takes a sheaf $\O(i),$ $0\le i \le l-1$ on
$\PP_{\theta'}$ to the sheaf $\O(i)$ on $\PP_{\theta}.$
Using the resolution (\ref{Kos1}) it can be easily checked that
this functor takes $\O(i)$ to $\O(i)$ for all $i\in \Z.$
Now, it follows from the ampleness condition on $\lbrace \O(i)\rbrace$
and Corollary \ref{zero}
that the functor  sends
the subcategory $\coh(\PP_{\theta'})$
to  $\coh(\PP_{\theta})$
and induces an equivalence
$\coh(\PP_{\theta'})\cong\coh(\PP_{\theta}).$
We just proved:

\begin{cor}\label{cor:equivtheta}
If the matrices $\theta'$ and $\theta$ are connected by the relation
(\ref{relth}) then the
noncommutative weighted projective spaces
$\PP_{\theta'}(\aa)$ and $\PP_{\theta}(\aa)$
have equivalent abelian categories of coherent sheaves
$\coh(\PP_{\theta'})$ and $\coh(\PP_{\theta}).$
\end{cor}

In the case $n=1,$ it follows immediately that
for any $\theta,\theta'\in M(2,\kk^*)$ the categories
$\coh(\PP_\theta(a_0,a_1))$ and $\coh(\PP_{\theta'}(a_0,a_1))$ are
equivalent.

Next consider the case $n=2.$
For any matrix $\theta\in M(3,\kk^*)$
 denote  the  expression
$$
(\theta^{an}_{01})^{a_2} (\theta^{an}_{12})^{a_0} (\theta^{an}_{20})^{a_1}=
(\theta_{01})^{a_2} (\theta_{12})^{a_0} (\theta_{20})^{a_1}
(\theta_{10})^{-a_2} (\theta_{21})^{-a_0} (\theta_{02})^{-a_1}
$$
by $q(\theta).$
Now, the result of Proposition
\ref{equiv} can be written in the following form.
\begin{cor}\label{cor:ncwp2}
Let $n=2$ and let $\theta'$ and $\theta$ be two matrices from
${\rm M}(3,\kk^*).$
If $q(\theta')=q(\theta)$ then the abelian categories
$\coh(\PP_{\theta'}(a_0,a_1,a_2))$ and
$\coh(\PP_{\theta}(a_0,a_1,a_2))$ are equivalent.
\end{cor}

\subsection{DG algebras and Koszul duality.}\label{ss:koszul}
The aim of this section is to give another description of the derived category
$\db{\coh(\PP_{\theta})}.$ It was shown above that
this category is equivalent to the derived category
$\db{\mod-\Be}.$ We introduce a finite dimensional
differential $\Z$-graded algebra (DG algebra)
$\Ce$ and prove that the category $\db{\coh(\PP_{\theta})}$ is equivalent
to the derived category  of $\Ce.$

This new description of the derived category in terms of the
DG-algebra $\Ce$ naturally yields an exceptional collection
(Corollary \ref{cor:excoll}), which is essentially the
(left) dual of the collection
described in Theorem \ref{thm:excoll}, cf.\ the discussion at the
end of \S \ref{ss:excoll}.

We recall here that a DG
 algebra over  $\kk$ is a graded associative $\kk$--algebra
$$
R=\bigoplus_{p\in \Z} R^p
$$
with a differential $d$ of degree $+1$
such that
$$
d(rs)=(dr)s+(-1)^p r(ds)
$$
for all $r\in R^p, s\in R.$
We will suppose that $R$ is noetherian as a graded algebra.

A right DG module over a  DG algebra is a graded right $R$--module
$M=\bigoplus_{p\in\Z} M^p$ with a differential $\nabla$ of degree 1 such that
$$
\nabla(mr)=(\nabla m)r+(-1)^p mdr
$$
for all $m\in M^p$ and $r\in R.$

A morphism of DG $R$-modules $f: M\lto N$ is called
null-homotopic if $f=d_N h + h d_M,$ where
$h: M\lto N$ is a morphism of the underlying graded
$R$-modules which is homogeneous of degree
$-1.$ The
homotopy category $\bH^b(R)$ is defined as a category which has all
finitely generated DG $R$-modules as objects,
and whose morphisms  are the equivalence classes
$\overline{f}$ of morphisms of DG $R$-modules modulo null-homotopic
morphisms.
A morphism of DG $R$-modules $s: M\to
N$ is called a quasi-isomorphism if the induced morphism $H^* s:
H^*(M)\to H^*(N)$ is an isomorphism of graded vector spaces.
Now, by definition, the  derived category $\bD^b(R)$ is  the
localization
$$
\bD^b(R):=\bH^b(R)\left[\Sigma^{-1}\right],
$$
where $\Sigma$ is the class of all quasi-isomorphisms.
It can be checked that there are canonical isomorphisms
$$
\Hom_{\bD^b(R)}(R, M)\stackrel{\sim}{\lto}
\Hom_{\bH^b(R)}(R, M)\stackrel{\sim}{\lto} H^0 M
$$
for each DG $R$-module $M.$

Any ordinary $\kk$-algebra $A$ can be considered as the DG algebra
$A^{\bdot}$ with
$A^0=A$ and $A^p=0$ for all $p\ne 0.$ In this case the derived category
of the DG algebra $\bD^b(A^{\bdot})$ identifies with
the bounded derived category of finitely generated
right $A$-modules,
i.e. $\bD^b(A^{\bdot})\cong \bD^b(\mod-A).$
For a detailed exposition of the facts about derived categories of DG algebras,
see \cite{Ke,Ke2}.

Now denote by $\Be_s$ the algebra $\Be/\rad(\Be)$ and consider it as
a right $\Be$-module, isomorphic
to the sum $\mathop\oplus\limits_{i=0}^{l-1}Q_i$ of all irreducibles.
Introduce the finite dimensional DG algebra
$$
\Ext^{\bdot}_{\Be}(\Be_s, \Be_s)=
\mathop\oplus\limits_{p\in\Z} \Ext^p_{\Be}(\Be_s, \Be_s)
$$
 with the natural composition law and trivial differential.
In what follows we  give a precise description of
 this DG algebra  and  prove the existence
of an equivalence
$$
\db{\coh(\PP_{\theta})}\cong\bD^b(\Ext^{\bdot}_{\Be}(\Be_s, \Be_s)),
$$
which gives the promised description of the category
$\bD^b(\coh(\PP_{\theta})).$

Let us introduce a graded DG algebra
$\varLambda^{\bdot}=\varLambda^{\bdot}(a_0,\ldots,a_n).$
As a DG algebra it is the skew-symmetric algebra with trivial differential
which is generated by skew-commutative elements
$y_i,$ $i=0,\,\ldots,\,l-1$ of degree $1,$ i.e.
$$
\varLambda^{\bdot}=\bigoplus_{p=0}^{n+1} \varLambda^{p},
$$
where $y_i\in \varLambda^{1}$ with the relations $y_iy_j=-y_jy_i$ for all
$0\le i,j\le n.$

The additional grading on the DG algebra $\varLambda^{\bdot}(a_0,\ldots,a_n)$
is defined by putting
$y_i\in \varLambda^{\bdot}_{-a_i}.$ Thus $\varLambda^{\bdot}(a_0,\ldots,a_n)$
is just a bigraded skew-symmetric algebra
$$
\varLambda^{\bdot}(a_0,\ldots,a_n)=\bigoplus_{p,i\in \Z}\varLambda^{p}_{i}
$$
with generators $y_i\in \varLambda^1_{-a_i}.$
For any $(n+1)\times(n+1)$-matrix $\theta$ we
also can define a graded DG algebra
$\varLambda^{\bdot}_{\theta}(a_0,\ldots,a_n)$ as the DG algebra with trivial
differential and generated
by elements $y_i\in (\varLambda_{\theta})^{1}_{-a_i}, i=0,\ldots,n$
with the relations
$$
\theta_{ij}y_i y_j +\theta_{ji}y_j y_i=0
$$
for all $0\le i,j\le n.$

Consider  the following complex $\Com$ of right $S_{\theta}$-modules
\begin{multline}\label{Kos2}
\Com:=0\to S_{\theta}(-\mathop\sum\limits_{i=0}^{n} a_i)
\to\bigoplus_{i_0<\ldots<i_{n-1}}S_{\theta}(-\mathop\sum\limits_{j=0}^{n-1}
a_{i_j}) \to
\cdots\\
\cdots\to
\bigoplus_{i_0<i_1}S_{\theta}(-a_{i_0}-a_{i_1})\to
\bigoplus_{i=0}^n S_{\theta}(-a_i)\to S_{\theta}\to
 0,
\end{multline}
in which the differentials are defined componentwise as follows:
for any set $I=\lbrace i_0,\ldots i_k\rbrace$
the differential sends the generator
of $S_{\theta}(-\mathop\sum\limits_{i\in I} a_i)$ to
the sum of the elements
$$
(-1)^s\left(\prod_{i\in I}\theta_{i i_s}\right) x_{i_s}$$
of $S_{\theta}\bigl(-\mathop\sum\limits_{i\in (I\backslash i_s)} a_i\bigr),$
for $0\le s\le k.$
With this  we see that the complex $\Com$
is a free resolution of the right $S_{\theta}$-module
$\kk_{S_{\theta}}.$

Now we define a structure of left DG module over the DG algebra
$\varLambda^{\bdot}_{\theta}$ on the complex
$\Com$, such that the element $y_j$ takes the generator
of $S_{\theta}(-\mathop\sum\limits_{i\in I} a_i)$ to the generator
of $S_{\theta}(-\mathop\sum\limits_{i\in (I\backslash i_s)} a_i)$
with coefficient
$$
(-1)^s\prod_{i\in I}\theta_{i_s i}
$$
if $j=i_s\in I=\{i_0,\dots,i_k\},$ and takes it to zero if $j\not\in I.$
It can be checked that this action is well defined and makes the complex
$\Com$ a DG $\varLambda^{\bdot}_{\theta}$-$S_{\theta}$-bimodule.
\begin{remark}\label{asbi} It is not difficult to see that the complex
$\Com$ as a graded
$\varLambda^{\bdot}_{\theta}$-$S_{\theta}$-bimodule
(i.e. without differential) is isomorphic to
$(\varLambda^{\bdot}_{\theta})^*\mathop\otimes\limits_{\kk} S_{\theta},$
where $(\varLambda^{\bdot}_{\theta})^*$ is $\Hom_{\kk}(\varLambda^{\bdot}_{\theta}, \kk).$
\end{remark}

\begin{defi}
Define  a DG category $\Cc$ (actually graded category,
because all differentials
are trivial) as a DG category with $l$ objects, say  $w_0,\ldots,w_{l-1},$
and the spaces of morphisms between which are the complexes
$$
\Hom^{\bdot}(w_j, w_i)\cong (\varLambda^{\bdot}_{\theta})_{i-j}
$$
with the natural composition law induced by that
of the DG algebra
$\varLambda^{\bdot}_{\theta}.$
\end{defi}
 It follows from the definition
of the DG algebra $\varLambda^{\bdot}_{\theta}$ that
$$
\Hom^{\bdot}(w_j,w_i)=0 \qquad\text{when}\quad j< i.
$$
\begin{defi}
Define the DG algebra $\Ce$ as the DG algebra of the
DG category $\Cc,$ i.e.
$$
\Ce:=\bigoplus_{0\le i,j\le l-1}\Hom^{\bdot}(w_j,w_i).
$$
\end{defi}
The quotient of this DG algebra by its radical is isomorphic to $\kk^{\oplus l}.$
In particular the DG algebra $\Ce,$ similarly to the algebra $\Be,$
has $l$ irreducible DG modules
in degree 0.
Moreover, as a right DG $\Ce$-module the algebra $\Ce$ is  a direct sum
$$
\Ce=\bigoplus_{i=0}^{l-1} H_{i},
\qquad
\text{where}
\quad
H_i=\bigoplus_{0\le j\le l-1}\Hom^{\bdot}(w_j,w_i),
$$
and $H_i$ are homotopically projective right DG $\Ce$-modules.

Let us construct a DG $\Ce\mbox{-}\Be$-bimodule $\Xe,$ obtained from
the DG $\varLambda^{\bdot}_{\theta}\mbox{-}S_{\theta}$-bimodule
$\Com$ by the formula
$$
\Xe=\bigoplus_{0\le i, j\le l-1} \Xe(i, j),
\qquad\text{with}
\quad
\Xe(i,j)\cong \Com_{j-i}
$$
where $\Com_{j-i}$ is the degree $(j-i)$
component of the graded complex $\Com.$
In particular, $\Xe(i,j)=0$ when $i>j$ and $\Xe(i,i)\cong\kk$ for all $i.$
The structure of DG $\Ce\mbox{-}\Be$-bimodule on $\Xe$ comes
from the structure of DG
$\varLambda^{\bdot}_{\theta}\mbox{-}S_{\theta}$-bimodule
on $\Com.$
The bimodule $\Xe$ is quasi-isomorphic to $\kk^{\oplus l},$
and it is quasi-isomorphic to $\Be/\rad(\Be)$ as a right $\Be$-module
and to $\Ce/\rad(\Ce)$ as a left DG $\Ce$-module.
This fact allows us to say that the DG algebra $\Ce$ is the
{\sf Koszul dual} to
the algebra $\Be.$

\begin{remark}\label{asbi2} It follows from Remark \ref{asbi}
that $\Xe$ as a graded
$\Ce\mbox{-}\Be$-bimodule (i.e. without differential) is isomorphic to
$$
\bigoplus_{i=0}^{l-1} H_i^*\otimes P_i,
$$
where $H_i^*$ are the left DG $\Ce$-modules $\Hom_{\kk}(H_i, \kk).$
In other words, as a graded
$\Ce\mbox{-}\Be$-bimodule $\Xe$ is isomorphic to
$\Ce^*\otimes_{\kk^{\oplus l}}\Be.$
\end{remark}

For any right DG $\Ce$-module $N$, the tensor product
$N\otimes_{\kk}\Xe$ is naturally a complex of right $\Be$-modules,
in which the module structure is given by the action of $\Be$ on $\Xe$,
and the grading and differential are given by
$$
\begin{array}{ll}
(N\otimes_{\kk}{\Xe})^k=\mathop\bigoplus\limits_{p+q=k}N^p\otimes_{\kk} X^q,
\qquad
&d(n\otimes x)=(dn)\otimes x +(-1)^p n\otimes dx
\end{array}
$$
for all $n\in N^p, x\in \Xe.$
The $\kk$-submodule generated by all differences
$nc\otimes x- m\otimes cx$ is closed under the differential and
under multiplication by any element of $\Be.$
So the quotient by this submodule, which we denote by
$N\otimes_{\Ce}\Xe,$ is a well-defined complex of right $\Be$-modules.

For any complex $M$ of right $\Be$-modules we define a right DG
$\Ce$-module
$$
\begin{array}{ll}
\hom_{\Be}(\Xe,M)^k=\mathop\prod\limits_{p-q=k}\Hom_{\Be}(X^q, M^p),
\qquad
&
(df)(x)=d(f(x))-(-1)^nf(dx).
\end{array}
$$
In this way we get a pair of adjoint functors
$(-)\otimes_{\Ce}\Xe$ and $\hom_{\Be}(\Xe,-)$
between homotopy categories,
which induce a pair of adjoint functors on the level of
derived categories as well:
$$
\stackrel{\bL}{\otimes}_{\Ce}\Xe: \bD^{b}(\Ce)\lto \bD^b(\mod\mbox{-}\Be),
\qquad
\bR\Hom_{\Be}(\Xe, -): \bD^b(\mod\mbox{-}\Be)\lto \bD^b(\Ce).
$$
Moreover, since $\Xe$ is a projective finitely generated  right
$\Be$-module and a flat
left $\Ce$-module,
both functors $(-)\otimes_{\Ce}\Xe$ and $\hom_{\Be}(\Xe,-)$
 between homotopy categories preserve acyclicity.
Hence, the derived functors in this case are defined by the same formulas.
For more information about derived functors see e.g. \cite{Ke}.

\begin{theo}
The functors
$
\stackrel{\bL}{\otimes}_{\Ce}\Xe$ and $
\bR\Hom_{\Be}(\Xe, -)$
are equivalences of triangulated categories.
\end{theo}
\begin{proof}
It is evident that the first functor $\stackrel{\bL}{\otimes}_{\Ce}\Xe$
takes
 $\Ce$ as a right DG $\Ce$-module to $\Xe$ as a right $\Be$-module
which is isomorphic to $\Be_s=\mathop\oplus\limits_{i=0}^{l-1} Q_i$
in the derived category
$\bD^b(\mod\mbox{-}\Be).$
On the other hand, it follows from Remark \ref{asbi2} and the equalities
$\Hom_{\Be}(P_i, Q_j)=\delta^i_j \kk$ that
the latter functor,
$\bR\Hom_{\Be}(\Xe, -),$
takes the module $\Be_s=\mathop\oplus\limits_{i=0}^{l-1} Q_i$
to the free DG module $\Ce=\bigoplus_{i=0}^{l-1} H_{i}$ and takes
$Q_i$ to $H_i$ for any
$0\le i \le l-1.$
Thus, the composition functor
$\bR\Hom_{\Be}(\Xe, -)\stackrel{\bL}{\otimes}_{\Ce}\Xe$
sends $B_s$ to itself and it also sends
all direct summands $Q_i$ to $Q_i.$ The adjunction maps
$$
\bR\Hom_{\Be}(\Xe, Q_i)\stackrel{\bL}{\otimes}_{\Ce}\Xe \lto Q_i
$$
cannot be trivial, hence they are isomorphisms for all $i.$
Therefore, we obtain isomorphisms
$$
\Hom_{\Be}(Q_i, Q_j[k])\stackrel{\sim}{\lto}
\Hom_{\Ce}(\bR\Hom_{\Be}(\Xe, Q_i), \bR\Hom_{\Be}(\Xe, Q_j)[k])\cong
\Hom(H_i, H_j[k])
$$
for any $0\le i,j \le l-1$ and all $k\in\Z.$

Since $Q_i, i=0,\ldots,l-1$ generate the derived category
$\bD^b(\mod\mbox{-}\Be),$  Lemma \ref{ff} implies that the functor
$$
\bR\Hom_{\Be}(\Xe, -):
\bD^b(\mod\mbox{-}\Be)\lto \bD^b(\Ce)
$$
is fully faithful.

Consider the triangulated subcategory $\D$ of $\bD^b(\Ce)$ generated by
$H_i,$ $i=0,\ldots,l-1.$ By Remark \ref{asbi2} $\Xe$ as a graded
$\Ce\mbox{-}\Be$-bimodule is isomorphic to
$\bigoplus_{i=0}^{l-1} H_i^*\otimes P_i,$ and hence, the dual to $\Xe$ over
$\kk$ gives a resolution of $\Ce/\rad(\Ce)$ in terms of $H_i.$
Therefore, the subcategory $\D$
contains  all
irreducible DG modules and coincides with
the whole $\bD^b(\Ce).$
Thus, $H_i,$ $i=0,\ldots,l-1$ generate the category  $\bD^b(\Ce),$ and
the functor $\bR\Hom_{\Be}(\Xe, -)$
is an equivalence of the derived categories.
\end{proof}

\begin{cor}
There is an isomorphism of DG algebras
$$
\Ce\cong\bigoplus_{0 \le i,j\le l-1}\Ext^{\bdot}( Q_i, Q_j).
$$
\end{cor}
The assertion of the Corollary  is clear now,
because the functor $\stackrel{\bL}{\otimes},$
which is an equivalence, sends $\Ce$ to $\Be_s=\mathop\oplus\limits_{i=0}^{l-1} Q_i.$

\begin{cor} \label{cor:excoll}
The derived category of coherent sheaves $\bD^b(\coh(\PP_{\theta}))$
on the noncommutative weighted space $\PP_{\theta}$ is equivalent
to the derived category $\bD^b(\Ce).$
\end{cor}

\subsection{Hirzebruch surfaces $\FF_n$}\label{ss:dbfn}
The  surfaces $\FF_n$ are minimal rational surfaces  defined as
the projectivizations $\Proj(\O\oplus \O(-n))$ of the vector bundles $\O\oplus\O(-n)$
over $\PP^1.$
The surface $\FF_n$ has a $(-n)$-section that will be denoted by $s.$
There is a simple connection between $\FF_n$ and the weighted projective plane
$\bP(1,1,n),$ namely the latter can be obtained from $\FF_n$ by contracting
 the $(-n)$\-section $s.$ In this way $\FF_n$ is a resolution of the singularity
of $\bP(1,1,n).$ Thus,  we have two different resolutions of the singularity
of $\bP(1,1,n)$:
$$
\xymatrix{
\FF_n \ar[dr]&&\PP(1,1,n)\ar[dl]\\
&\bP(1,1,n)}
$$
For this reason the derived categories of coherent sheaves on
$\FF_n$ and on $\PP(1,1,n)$
are closely related to each other.
We will show that for $n\ge 2$ there is a fully faithful functor
$$
\fun: \bD^b(\coh(\FF_n))\lto\bD^b(\coh(\PP(1,1,n)))
$$
and will give its description.

Denote by $f$ the class of the fiber of $\FF_n$ in the Picard group.
Since $\FF_n$ is a $\PP^1$-bundle over $\PP^1$ the derived category
of coherent sheaves on $\FF_n$ has an exceptional collection
of length 4 (see \cite{O1}). More precisely, we have

\begin{prop}
The collection $\sigma=\left( \O, \O(f), \O(s+nf), \O(s+(n+1)f)\right)$
is a full strong exceptional collection on $\FF_n.$
The derived category $\bD^b(\coh(\FF_n))$ is equivalent to
the derived category $\bD^b(\mod\mbox{-}\De),$ where $\De$ is the algebra
of the exceptional collection~$\sigma.$
\end{prop}

Denote by $U$ the two dimensional vector space $H^0(\FF_n, \O(f)).$
 From the exact sequence
$$
0\lto \O\lto\O(s+nf)\lto \O_s\lto 0
$$
we find that $H^0(\FF_n, \O(s+nf))$ is the direct sum of the space $S^nU$ and
a one-dimensional space. Analogously, we can check that $H^0(\FF_n, \O(s+(n+1)f))$
is isomorphic to $S^n U \oplus U.$

On the other hand, we know that the weighted projective plane
$\PP(1,1,n)$ has an exceptional collection
$$
\left(\O,\O(1),\ldots,\O(n),\O(n+1)\right).
$$
Denote the algebra of this exceptional collection by $B(1,1,n).$
It follows from Proposition \ref{cohom} that the space
$H^0(\PP(1,1,n), \O(1))$ is isomorphic to $U,$ $H^0(\PP(1,1,n),
\O(n))$ is isomorphic to the direct sum of $S^n U$ and a
one-dimensional space, and $H^0(\PP(1,1,n), \O(n+1))$ is
isomorphic to $S^nU\oplus U.$ This implies that the algebra of the
exceptional collection $\left( \O, \O(f), \O(s+nf),
\O(s+(n+1)f)\right)$ on $\FF_n$ is isomorphic to the algebra of
the exceptional collection $\left( \O, \O(1), \O(n),
\O(n+1)\right)$ on $\PP(1,1,n).$

Thus, the algebra of endomorphisms of the projective $B(1,1,n)$-module
$$
M=P_0\oplus P_1\oplus P_n\oplus P_{n+1}
$$
coincides with $\De,$ which makes $M$ a $\De\mbox{-}B(1,1,n)$-bimodule.
The natural functor
$$
(-)\stackrel{\bL}{\otimes}_{\De}M: \bD^b(\mod\mbox{-}\De)\lto \bD^b(\mod\mbox{-}B(1,1,n))
$$
takes the free module $\De$ to $M,$ and there are isomorphisms
$$
\Hom_{\De}(\De,\De[k])\stackrel{\sim}{\lto}\Hom_{B(1,1,n)}(M,M[k]).
$$
Since the direct summands of $\De$ generate the derived category
$\bD^b(\mod\mbox{-}\De),$ Lemma \ref{ff} guarantees that the functor
$(-)\stackrel{\bL}{\otimes}_{\De}M$ is fully faithful.
Using the descriptions of the derived categories of coherent
sheaves on $\FF_n$ and $\PP(1,1,n)$ in terms of the exceptional collections,
we obtain the following theorem.

\begin{theo}
The functor
$$
\fun: \bD^b(\coh(\FF_n))\lto\bD^b(\coh(\PP(1,1,n)))
$$
induced by $(-)\stackrel{\bL}{\otimes}_{\De}M$ is fully faithful.
\end{theo}

\section{Categories of Lagrangian vanishing cycles}
\label{sec:fsdef}

\subsection{The category of vanishing cycles of an affine Lefschetz fibration}
\label{ss:fs}
We begin this section by briefly reviewing Seidel's construction of
a Fukaya-type $A_\infty$-category associated to a symplectic Lefschetz
fibration \cite{Se1,Se2,Se3}, following a proposal of Kontsevich \cite{Ko}.
For an account of the underlying physics, the reader is referred to
the work of Hori et al \cite{HIV}.

Let $(X,\omega)$ be an open symplectic manifold, and let $f:X\to \C$ be a
symplectic Lefschetz fibration, i.e.\ a $C^\infty$ complex-valued function
with isolated non-degenerate critical points $p_1,\dots,p_r$ near which
$f$ is given in local coordinates by $f(z_1,\dots,z_n)=f(p_i)+z_1^2+
\dots+z_n^2$, and whose fibers are symplectic submanifolds of $X$.
Fix a regular value $\lambda_0$ of $f$, and consider an arc $\gamma\subset\C$
joining $\lambda_0$ to a critical value $\lambda_i=f(p_i)$.
Using the horizontal distribution given by the symplectic orthogonal to
the fibers of $f$, we can transport the vanishing cycle at $p_i$ along the
arc $\gamma$ to obtain a Lagrangian disc $D_\gamma\subset X$ fibered above
$\gamma$, whose boundary is an embedded Lagrangian sphere $L_\gamma$ in the
fiber $\Sigma_0=f^{-1}(\lambda_0)$. When the fibers of $f$ are non-compact,
parallel transport along the horizontal distribution is not always
well-defined; we will always assume that the symplectic form $\omega$
satisfies the conditions required to make the construction valid.
The Lagrangian disc $D_\gamma$ is
called the {\it Lefschetz thimble} over $\gamma$, and its boundary
$L_\gamma$ is the vanishing cycle associated to the critical point
$p_i$ and to the arc $\gamma$.

Let $\gamma_1,\dots,\gamma_r$ be a collection of arcs in $\C$ joining the
reference point $\lambda_0$ to the various critical values of $f$,
intersecting each other only at $\lambda_0$, and ordered in the clockwise
direction around $p_0$. Each arc $\gamma_i$ gives rise to a Lefschetz
thimble $D_i\subset X$, whose boundary is a Lagrangian sphere $L_i\subset
\Sigma_0$. After a small perturbation we can always assume that these
spheres intersect each other transversely inside $\Sigma_0$.

\begin{defi}[Seidel]\label{def:fs}
The directed category of vanishing cycles $\FS(f,\{\gamma_i\})$ is an
$A_\infty$-category (over a coefficient ring $R$) with $r$ objects
$L_1,\dots,L_r$ corresponding to the vanishing cycles (or more accurately
to the thimbles); the morphisms between
the objects are given by
$$\mathrm{Hom}(L_i,L_j)=\begin{cases}
CF^*(L_i,L_j;R)=R^{|L_i\cap L_j|} & \mathrm{if}\ i<j\\
R\cdot id & \mathrm{if}\ i=j\\
0 & \mathrm{if}\ i>j;
\end{cases}$$
and the differential $m_1$, composition $m_2$ and higher order
products $m_k$ are defined in terms of Lagrangian Floer homology inside
$\Sigma_0$. More precisely,
$$m_k:\mathrm{Hom}(L_{i_0},L_{i_1})\otimes \dots\otimes
\mathrm{Hom}(L_{i_{k-1}},L_{i_k}) \to \mathrm{Hom}(L_{i_0},L_{i_k})[2-k]$$
is trivial when the inequality $i_0<i_1<\dots<i_k$ fails to hold (i.e.\ it
is always zero in this case, except for $m_2$ where composition with an
identity morphism is given by the obvious formula).
When $i_0<\dots<i_k$, $m_k$ is defined by fixing a generic
$\omega$-compatible almost-complex structure on $\Sigma_0$ and counting
pseudo-holomorphic maps from a disc with $k+1$ cyclically ordered
marked points on its boundary to $\Sigma_0$, mapping the marked points
to the given intersection points between vanishing cycles, and the portions
of boundary between them to $L_{i_0},\dots,L_{i_k}$ respectively.
\end{defi}

While the general definition of Lagrangian Floer homology is a very
delicate task \cite{FO3}, we will only consider cases where most of the
technical considerations can be skipped. For example, Seidel considers
the case where the symplectic form $\omega$ is exact ($\omega=d\theta$ for
some $1$-form $\theta$) and the $L_i$ are exact Lagrangian submanifolds in
$\Sigma_0$ (i.e.\ $\theta_{|L_i}=dg_i$ is also exact). Here, we assume
instead that
the restricted symplectic form $\omega_{|\Sigma_0}$ is exact and that the
homotopy groups $\pi_2(\Sigma_0)$ and $\pi_2(\Sigma_0,L_i)$ are trivial.
The first condition
prevents the bubbling of pseudo-holomorphic spheres, while the second one
prevents the bubbling of pseudo-holomorphic discs in the definition of
Lagrangian Floer homology. Therefore, the moduli spaces of
pseudo-holomorphic maps involved in the definition of $\FS(f,\{\gamma_i\})$
have well-defined fundamental classes.

Another assumption that we will make concerns the Maslov class, which we
will assume to vanish over $L_i$. In fact, we will restrict ourselves to
the case where $X$ and $\Sigma_0$ are affine Calabi-Yau manifolds, and the
spheres $L_i$ can be lifted to graded Lagrangian submanifolds of $\Sigma_0$,
e.g.\ by fixing a holomorphic volume form on $\Sigma_0$ and choosing a real
lift of the phase $\exp(i\phi)=\Omega_{|L_i}/vol_{L_i}:L_i\to S^1$.
This makes it possible to define a $\Z$-grading (by Maslov index) on the
Floer complexes $CF^*(L_i,L_j;R)$, as will be discussed later (see also
\cite{Se1}).

For simplicity, Seidel uses $R=\Z/2$ as coefficient ring in his
definition; however the moduli spaces considered below are orientable, so it
is possible to assign a sign $\pm 1$ to each pseudo-holomorphic curve
and hence define Floer homology over $\Z$. We will further extend the
coefficient ring to $R=\C$, and count the contribution of each
pseudo-holomorphic disc $u:(D^2,\partial D^2)\to (\Sigma_0,\bigcup L_i)$
in the moduli space with a coefficient of the form $\pm \exp(-\int_{D^2}
u^*\omega)$. Weighting by area is irrelevant in the case of exact
Lagrangian vanishing cycles considered by Seidel, where it does not affect
at all the structure of the category: indeed, the symplectic
areas can then be expressed in terms of the primitives $g_i$ of
$\theta$ over $L_i$, and can be eliminated from the description simply
by a rescaling of the chosen bases of the Floer complexes (considering the
basis $\{\exp{(g_i(p)-g_j(p))}\,p,\ p\in L_i\cap L_j\}$ of $CF^*(L_i,L_j)$).
On the contrary, in the non-exact case it is important to incorporate this
weighting by area into the definition.

Hence, given two intersection points $p\in L_i\cap L_j$, $q\in L_j\cap
L_k$ ($i<j<k$), we have by definition
$$m_2(p,q)=\sum_{\substack{r\in L_i\cap L_k\\\deg r=\deg p+\deg q}}
\Biggl(\sum_{[u]\in \mathcal{M}(p,q,r)}\!\!\pm \exp(-\int_{D^2} u^*\omega)
\Biggr)\,r$$
where $\mathcal{M}(p,q,r)$ is the moduli space of pseudo-holomorphic maps
$u$ from the unit disc to $M$ (equipped with a generic $\omega$-compatible
almost-complex structure) such that $u(1)=p$, $u(\mathrm{j})=q$,
$u(\mathrm{j}^2)=r$ (where $\mathrm{j}=\exp(\frac{2i\pi}{3})$), and mapping the
portions of unit circle $[1,\mathrm{j}]$, $[\mathrm{j},\mathrm{j}^2]$,
$[\mathrm{j}^2,1]$ to $L_i$, $L_j$ and $L_k$ respectively.
The other products are defined similarly.

It is worth mentioning that this definition of Floer homology over complex
numbers is in fact
essentially equivalent to the use of coefficients in a Novikov ring, since
in both cases the main goal is to keep track of (relative) homology classes.

Although the category $\FS(f,\{\gamma_i\})$ depends on the chosen ordered
collection of arcs $\{\gamma_i\}$, Seidel has obtained the following result
\cite{Se1}:

\begin{theo}[Seidel]
If the ordered collection $\{\gamma_i\}$ is
replaced by another one $\{\gamma'_i\}$, then the categories
$\FS(f,\{\gamma_i\})$ and $\FS(f,\{\gamma'_i\})$ differ by a sequence of
mutations.
\end{theo}

Hence, the category naturally associated to the
Lefschetz fibration $f$ is not the finite directed category defined above,
but rather a (bounded) {\it derived} category, obtained from $\FS(f,\{\gamma_i\})$
by considering twisted complexes of formal direct sums of Lagrangian
vanishing cycles, and adding idempotent splittings and formal inverses
of quasi-isomorphisms. It is a classical result that, if two
categories differ by mutations, then their derived categories are
equivalent; hence the derived category $D(\FS(f))$ only depends on the
Lefschetz fibration $f$ rather than on the choice of an ordered system
of arcs \cite{Se1}.

We finish this overview with a couple of remarks. In ``usual'' Fukaya
categories, objects are pairs consisting of a compact Lagrangian submanifold
and a flat connection on some complex vector bundle defined over it. In the
case of the category associated to a Lefschetz fibration, the objects are
vanishing cycles, or perhaps more accurately, the Lefschetz thimbles bounded
by the vanishing cycles. Since the thimbles are contractible, all flat
vector bundles over them are trivial, which eliminates the need to consider
Floer homology with twisted coefficients. This ceases to be true in
presence of a non-trivial B-field, but even then the equivalence class of
the connection is entirely determined by the thimble. Another related issue
is the choice of a spin structure on the vanishing cycles in order to fix
the orientation on the moduli spaces: in the one-dimensional case that will
be of interest to us, each vanishing cycle admits two distinct spin
structures ($H^1(S^1,\Z/2)=\Z/2$). However we must always consider the spin
structure which extends to the thimble, i.e.\ the non-trivial one.

The reader is referred to Seidel's papers \cite{Se1,Se2} for various
examples -- we will focus specifically on
the Landau-Ginzburg models mirror to weighted projective spaces and
Hirzebruch surfaces.

\subsection{Structure of the proof of Theorem \ref{thm:main}}
\label{ss:outline}
Derived categories of coherent sheaves on
the weighted projective planes $\PP^2(a,b,c)$ and their noncommutative
deformations $\PP^2_\theta(a,b,c)$ have been described in Chapter \ref{sec:algebra}.
Hence, to prove Theorem \ref{thm:main}, we need to find a similar
description of the derived categories of Lagrangian vanishing cycles on
the mirror Landau-Ginzburg models.

Recall that the mirror to $\PP^2_\theta(a,b,c)$ is $(X,W)$, where $X$ is the
affine hypersurface $\{x^ay^bz^c=1\}\subset(\C^*)^3$, equipped with an exact
(for the commutative case) or non-exact (for the noncommutative case)
symplectic form, and the superpotential $W=x+y+z.$

By construction, categories of Lagrangian vanishing cycles for Lefschetz
fibrations always admit full exceptional collections. Indeed, for any choice
of arcs $\{\gamma_i\}$ the objects $L_i$ of $\FS(W,\{\gamma_i\})$ form
a generating exceptional collection of the derived category. Hence, in view
of Theorem \ref{thm:excoll} and Corollary \ref{cor:excoll}, all we need to
do is exhibit a set of arcs $\{\gamma_i\}$ for which $\FS(W,\{\gamma_i\})$
is equivalent to one of the categories $\Bc$ or $\Cc$ introduced
in \S \ref{sec:algebra} (it turns out that the latter choice is slightly
easier to achieve).

Recall from Corollary \ref{cor:excoll} that $\db{\coh(\PP^2_\theta(a,b,c))}$
is equivalent to the derived category of the DG-algebra $\Ce$
associated to the finite DG-category $\Cc$ which has $l=a+b+c$ objects
$w_0,\ldots,w_{l-1},$ with morphisms between them given by the complexes
$$
\Hom^{\bdot}(w_j, w_i)\cong (\varLambda^{\bdot}_{\theta})_{i-j}
$$
with the natural composition law induced by that of the deformed exterior
algebra $\varLambda^{\bdot}_{\theta}$ on three generators of degrees
$-a,-b,-c$, with relations of the form $\theta_{ij}y_iy_j+\theta_{ji}y_jy_i$
where $\theta\in M(3,\C^*)$ (see \S \ref{ss:koszul}).
Moreover, by Corollary \ref{cor:ncwp2}, this category depends only on
the quantity $$q(\theta)=(\theta_{01})^{c} (\theta_{12})^{a} (\theta_{20})^{b}
(\theta_{10})^{-c} (\theta_{21})^{-a} (\theta_{02})^{-b}.$$

 From a practical viewpoint,
the cyclic group $\Z/(a+b+c)$ acts by diagonal multiplication
on $X$, and the superpotential $W=x+y+z$ is equivariant with respect to
this action. The $(a+b+c)$ critical values of $W$ form a
single orbit under this action (see \S \ref{ss:vcs}).
In order to exploit this symmetry, it is therefore natural to choose
the smooth fiber $\Sigma_0=W^{-1}(0)$ as our reference fiber, and an ordered
system of arcs $\gamma_i\subset\C$ ($i=0,\dots,a+b+c-1$) consisting of
straight line segments from the origin to the various critical values
$\lambda_i$.

With this understood, Theorem \ref{thm:main} follows immediately from
Corollary \ref{cor:excoll} and the following statement:

\begin{theo}\label{thm:fswp2}
$\FS(W,\{\gamma_i\})$ is a DG category, and it is equivalent to $\Cc$
for any $\theta\in \mathrm{M}(3,\kk^*)$ such that $q(\theta)=\exp(i[B+i\omega]\cdot[T])$,
where $[B+i\omega]\in H^2(X,\C)$ is the complexified K\"ahler class, and $[T]$ is
the generator of $H_2(X,\Z)$.
\end{theo}

The proof of Theorem \ref{thm:fswp2} consists of several steps, carried
out in the various subsections of \S \ref{sec:fsp2}. First, as a prerequisite to the determination of
the vanishing cycles, one needs a convenient description of the reference
fiber $\Sigma_0$. This is done by considering the projection to the first
coordinate axis,
$\pi_x:\Sigma_0\to\C^*$, which makes $\Sigma_0$ a $(b+c)$-fold covering
of $\C^*$ branched in $(a+b+c)$ points (Lemma \ref{l:sigma0}).
With this understood, it becomes fairly easy to identify the vanishing
cycles associated to the arcs $\gamma_j$, at least in the special case where the symplectic form is
anti-invariant under complex conjugation (which implies its exactness).
Indeed, this assumption implies that the vanishing cycles $L_j$ are {\it
Hamiltonian isotopic} (and hence equivalent from the point of view of
Floer theory) to the double lifts via $\pi_x$ of certain arcs
$\delta_j\subset\C^*$ (Lemma \ref{l:vcs}) which can be described explicitly
(Figure \ref{fig:vcs}).

With an explicit description of the vanishing cycles at hand, it becomes
possible to understand the Floer complexes $CF^*(L_i,L_j)$,
by studying the intersections between $L_i$ and $L_j$ for all $0\le
i<j<a+b+c$. Using the projection to the first coordinate, these correspond
to certain specific intersections between the arcs $\delta_i$ and $\delta_j$
in $\C^*,$ as dictated by
the combinatorics of the branched covering $\pi_x$.
Such a description is given by Lemma \ref{l:isects}, from
which it follows readily that $CF^*(L_i,L_j)\simeq
(\varLambda^{\bdot}_\theta)_{i-j}$ for all $i,j.$

The next step is to study the Floer differentials and products in
$\FS(W,\{\gamma_i\})$ by counting pseudo-holomorphic maps from
$(D^2,\partial D^2)$ to $(\Sigma_0,\bigcup L_i)$. This is done by
searching for immersed polygonal regions in $\Sigma_0$ with boundary
contained in $\bigcup L_i$, or equivalently, images of such regions
under the projection $\pi_x$ (see \S \ref{ss:products}). In our case,
it turns out that the only possible contributions come from triangular
regions in $\Sigma_0$; hence, the Floer differential $m_1$ and the
higher compositions $(m_k)_{k\ge 3}$ are identically zero (Lemmas
\ref{l:isects} and \ref{l:m3}) for purely topological reasons, while
the Floer product $m_2$ has a particularly simple structure (Lemma
\ref{l:m2}). In
particular, the $A_\infty$-category $\FS(W,\{\gamma_i\})$ is actually
a DG category with trivial differential.

The grading in $\FS(W,\{\gamma_i\})$ is determined by the Maslov indices
of intersection points. Since the Maslov class vanishes, each $L_i$ can be
lifted to a {\it graded Lagrangian} submanifold of $\Sigma_0$ by choosing
a real lift of its phase function (see \S \ref{ss:grading}). The degree of
a given intersection point $p\in L_i\cap L_j$ is then determined by
the difference between the phases of $L_i$ and $L_j$ at $p.$ Although the
determination of phases is the most technical part of the argument, it
actually presents little conceptual difficulty, and after some calculations
one readily checks that the grading of morphisms in $\FS(W,\{\gamma_i\})$
is the expected one. Namely, the ``generating'' morphisms corresponding to
the generators of the deformed exterior algebra $\varLambda^{\bdot}_\theta$
have degree 1, and their pairwise products have degree 2 (cf.\ Lemma
\ref{l:grading}).

The argument is then completed by determining more precisely the structure
coefficients for the Floer product $m_2$, which depend on the
symplectic areas of the various pseudo-holomorphic discs and on the choice
of consistent orientations of the moduli spaces (see \S \ref{ss:antisym}).
In the case where the
symplectic form is anti-invariant under complex conjugation, the argument
is greatly simplified by symmetry considerations, and the Floer products
obey the anticommutation rules of an (undeformed) exterior algebra (Lemma
\ref{l:antisym}) -- recall that complex conjugation anti-invariance implies
exactness of the symplectic form. In the non-exact case or in presence of a
non-zero B-field, there is no simple method for determining the symplectic
areas of the various pseudo-holomorphic discs involved in the definition of
$m_2$. However the deformation of the category $\FS(W,\{\gamma_i\})$ is
governed by a single parameter (analogous to the quantity $q(\theta)$
introduced in Corollary \ref{cor:ncwp2}), for which
a simple topological interpretation can be found, involving only the
evaluation of $[B+i\omega]$ on the generator of $H_2(X,\Z)$
(Lemmas \ref{l:deform} and \ref{l:deform2}).

This provides the desired characterization of the category of Lagrangian
vanishing cycles, and Theorem \ref{thm:fswp2} becomes an easy
corollary of Lemmas \ref{l:isects}--\ref{l:deform2}. The
only subtle point is that the objects of the category $\Cc$ are numbered
``backwards'' (because the generators of
$\varLambda^{\bdot}_\theta$ are assigned {\it negative} degrees), so the
equivalence of categories actually takes the objects $L_0,\dots,
L_{a+b+c-1}$ of $\FS(W,\{\gamma_i\})$ to the
objects $w_{a+b+c-1},\dots,w_0$ of~$\Cc.$

\subsection{Mirrors of weighted projective lines}
As a warm-up example, we prove HMS for the weighted projective lines
$\CP^1(a,b),$ where $a,b$ are mutually prime positive integers (see also
\cite{Se2} and \cite{DVS}). The
argument is an extremely simplified version of that outlined in \S
\ref{ss:outline}.
Indeed, the mirror Landau-Ginzburg model is the curve
$X=\{x^ay^b=1\}\subset(\C^*)^2$ equipped with the superpotential $W=x+y$,
whose generic fiber is just a finite set of $a+b$ points; so
most of the considerations that arise
in the case of weighted projective planes are irrelevant here (in particular
the symplectic structure on $X$ plays no role whatsoever, which is
consistent with the fact that the category $\coh(\PP_\theta(a,b))$ does not
depend on $\theta$).

More precisely, the fiber of $W$ above a point $\lambda\in\C$ is
$$W^{-1}(\lambda)=\{(x,\lambda-x)\in(\C^*)^2,\ x^a(\lambda-x)^b=1\},$$
which consists of $a+b$ distinct points, unless $P(x)=x^a(\lambda-x)^b-1$
has a double root. Since
$$P'(x)=\Bigl(\frac{a}{x}-\frac{b}{\lambda-x}\Bigr)(P(x)+1),$$ a root of $P$ is a
double root if and only if $x=\frac{a}{a+b}\lambda$; hence a double root
exists if and only if $P(\frac{a}{a+b}\lambda)=0$, i.e.\
\begin{equation}\label{eq:p1crit}
\lambda^{a+b}=\frac{(a+b)^{a+b}}{a^ab^b}.\end{equation}

Let $\lambda_0$ be the positive real root of this equation, and let
$\lambda_j=\lambda_0\zeta^{-j}$ where $\zeta=\exp(\frac{2\pi i}{a+b})$:
then the critical values of $W$ are exactly
$\lambda_0,\dots,\lambda_{a+b-1}$.
We choose $\Sigma_0=W^{-1}(0)$ as our reference fiber, and consider the
ordered system of arcs $\gamma_0,\dots,\gamma_{a+b-1}$, where
$\gamma_j\subset\C$ is a straight line segment joining the origin to
$\lambda_j$. With this understood, we have the following result, which
implies that HMS holds for $\CP^1(a,b)$:

\begin{theo}\label{thm:fswp1}
$\FS(W,\{\gamma_i\})$ is a DG category, equivalent to $\Cc$
for any $\theta\in \mathrm{M}(2,\kk^*).$
\end{theo}

In order to prove Theorem \ref{thm:fswp1}, we study the vanishing cycles
of the superpotential $W$ and their intersection properties.
To start with, observe that $W$ is equivariant with respect to the
diagonal action of the cyclic group $\Z/(a+b).$
Therefore, the vanishing cycles $L_j\subset\Sigma_0$ (which are Lagrangian
0-spheres, i.e.\ pairs of points) form a single $\Z/(a+b)$-orbit,
and $L_j=\zeta^{-j}\cdot L_0.$

\begin{figure}[t]
\centering
\setlength{\unitlength}{0.8cm}
\begin{picture}(5,3.2)(-2.5,-1.5)
\put(-2.5,-1.5){\small $\lambda=0$}
\put(-1.5,0){\vector(1,0){3}}
\put(0,-1.5){\vector(0,1){3}}
\put(-1.31,0){\circle*{0.08}}
\put(-0.81,-1.02){\circle*{0.08}}
\put(-0.81,1.02){\circle*{0.08}}
\put(0.29,-1.28){\circle*{0.08}}
\put(0.29,1.28){\circle*{0.08}}
\put(1.18,-0.57){\circle*{0.08}}
\put(1.18,0.57){\circle*{0.08}}
\qbezier[50](0.29,-1.28)(0.7,0)(0.29,1.28)
\put(1.55,0.1){\tiny Re\,$x$}
\put(0.1,1.55){\tiny Im\,$x$}
\end{picture}
\begin{picture}(5,3)(-1.5,-1.5)
\put(0.5,-1.5){\small $\lambda=\lambda_0$} 
\put(-1,0){\vector(1,0){4.2}}
\put(0,-1.5){\vector(0,1){3}}
\put(-0.66,0){\circle*{0.08}}
\put(-0.14,-0.72){\circle*{0.08}}
\put(-0.14,0.72){\circle*{0.08}}
\put(1.48,0){\circle*{0.08}}
\put(1.48,0){\circle*{0.08}}
\put(2.78,-0.35){\circle*{0.08}}
\put(2.78,0.35){\circle*{0.08}}
\end{picture}
\begin{picture}(6,3)(-1.5,-1.5)
\put(0.5,-1.5){\small $\lambda\to+\infty$} 
\put(-1,0){\vector(1,0){5}}
\put(0,-1.5){\vector(0,1){3}}
\put(-0.59,0){\circle*{0.08}}
\put(-0.09,-0.64){\circle*{0.08}}
\put(-0.09,0.64){\circle*{0.08}}
\put(0.74,0){\circle*{0.08}}
\put(2.69,0){\circle*{0.08}}
\put(3.41,-0.29){\circle*{0.08}}
\put(3.41,0.29){\circle*{0.08}}
\end{picture}
\caption{The fiber of $W$ for $\lambda\in\R_+$ ($(a,b)=(4,3)$)}
\label{fig:vcline}
\end{figure}
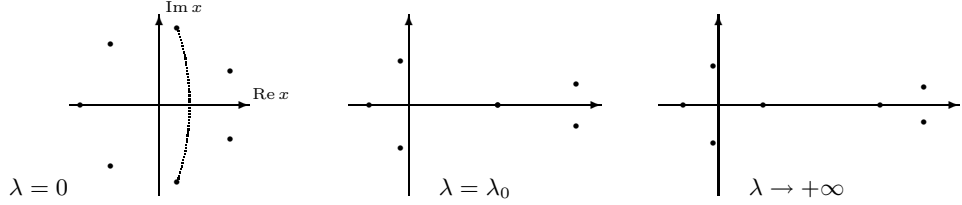

In order to determine $L_0$, we study how the fiber $W^{-1}(\lambda)$
varies as $\lambda$ increases along the positive real axis (see Figure
\ref{fig:vcline}). For $\lambda=0$,
the fiber $\Sigma_0$ consists of $a+b$ points whose first coordinates are
the roots of the equation $x^{a+b}=(-1)^b$ (these form a $\Z/(a+b)$-orbit,
hence the points of $\Sigma_0$ can naturally be identified with the elements
of $\Z/(a+b)$ up to a translation). As $\lambda$ increases towards
$\lambda_0$, two complex conjugate points of the fiber converge towards
each other, and become real points for $\lambda>\lambda_0$. By
considering the situation for $\lambda\to +\infty$, where the solutions of
$x^a(\lambda-x)^b=1$ split into two groups, one consisting of $a$ roots
near the origin, and the other consisting of $b$ roots near $\lambda$, one
easily checks that the vanishing cycle $L_0$ consists of the two points of
$\Sigma_0$ with first coordinate $x=\exp(\pm \frac{i\pi\,b}{a+b})$.

Hence, for a suitable identification of the fiber $\Sigma_0$ with
$\Z/(a+b)$, the vanishing cycle associated to the arc
$\gamma_0=[0,\lambda_0]$ is $L_0=\{0,b\}$. It follows immediately that
$L_j=\zeta^{-j}\cdot L_0 = \{-j,b-j\}$ for all $j=0,1,\dots,a+b-1.$

Given $0\le i<j<a+b$, the vanishing cycles $L_i$ and $L_j$ intersect
if and only if the subsets $\{-i,b-i\}$ and $\{-j,b-j\}$ of $\Z/(a+b)$
have non-empty intersection, i.e.\ if $j=i+a$ or $j=i+b$. Therefore,
we have:

\begin{lemma}\label{l:isectsp1}
The direct sum $\bigoplus_{i<j} CF^*(L_i,L_j)$ is a free module of total
rank $(a+b)$ over the coefficient ring, generated by the
intersection points
$$x_i\in CF^*(L_i,L_{i+a})\quad (0\le i<b)\qquad\mathrm{and}\qquad
y_i\in CF^*(L_i,L_{i+b})\quad (0\le i<a).$$
\end{lemma}

Because $\Sigma_0$ is a discrete set, all pseudo-holomorphic curves in
$\Sigma_0$ must be constant maps. However, each point of $\Sigma_0$
occurs exactly once as an intersection between two vanishing cycles
(there are no triple intersections), which implies that the Floer
differentials and products are trivial. Therefore, we have:

\begin{lemma}
The differentials and products $m_k,\ k\ge 1$ in the
$A_\infty$-category $\FS(W,\{\gamma_i\})$
are all identically zero, with the exception of the obvious
ones $m_2(\cdot,id)$ and $m_2(id,\cdot).$
\end{lemma}

This of course greatly simplifies the argument, eliminating the need for
many of the arguments required in the case of higher-dimensional weighted
projective spaces. At this point, our only remaining task is to determine
the Maslov indices of the various intersection points, by choosing graded
Lagrangian lifts of the vanishing cycles. A word of warning is in order
here: because we are actually dealing with graded Lagrangian submanifolds
in a Calabi-Yau 0-fold, the argument is very specific (see \S 2 of
\cite{Se2} for a discussion of graded Lagrangian submanifolds of
0-dimensional symplectic manifolds) and does not give
a good intuition of the higher-dimensional case.

\begin{lemma}\label{l:gradingp1}
There exists a natural choice of gradings for which $\deg(x_i)=\deg(y_i)=1$.
\end{lemma}

\proof
Equip the curve $X=\{x^ay^b=1\}\subset(\C^*)^2$ with the complex
structure induced by the standard one. The holomorphic volume form
$d\log x\wedge d\log y$ on $(\C^*)^2$ induces a $(1,0)$-form $\Omega$ on $X,$
characterized by the property that it is the restriction to $X$ of a
$1$-form (which we also call $\Omega$) such that $\Omega\wedge d(x^ay^b)=d\log x\wedge d\log y$,
i.e., using the fact that $x^ay^b=1$ along $X$,
$$\Omega\wedge \bigl(\frac{a}{x}\,dx+\frac{b}{y}\,dy\bigr)=\frac{dx\wedge
dy}{xy}.$$
Outside of the branch points of $W$, the 1-form $\Omega$ can be expressed
as $\Theta\,dw$, for some meromorphic function $\Theta$ with simple poles
at the branch points. The above equation becomes
$\Theta(\frac{b}{y}-\frac{a}{x})=\frac{1}{xy}$, i.e.\ $\Theta=(bx-ay)^{-1}=
((a+b)x-aw)^{-1}$.
In particular, near $\Sigma_0=W^{-1}(0)$, we have $\arg \Theta= -\arg x$.

The complex-valued function $\Theta$ is (up to scaling by a positive real
factor) the natural holomorphic
volume form induced by $\Omega$ on the 0-dimensional manifold
$\Sigma_0=W^{-1}(0)$. Let $L_0=\{p_-,p_+\}$, where the $x$-coordinate
of $p_{\pm}$ is $x_\pm=\exp(\pm \frac{i\pi b}{a+b})$. The {\it phase} of $L_0$
is the function $\phi_{L_0}:L_0\to\R/\pi\Z$ defined by
$$\phi_{L_0}(p_\pm)=\arg \Theta(p_\pm)=\mp \frac{\pi b}{a+b}.$$ Note that
an orientation on $L_0$ determines a lift of $\phi_{L_0}$ to
a $\R/2\pi\Z$-valued function; in order
to define the Maslov index, we need to view $L_0$ as a
{\it graded Lagrangian} submanifold, i.e.\ to choose a real lift
$\tilde{\phi}_{L_0}:L_0\to\R$ of the phase function. Although there is
{\it a priori} a $\Z^2$-space of such choices, one has to restrict
oneself to only those lifts which are compatible
with a graded Lagrangian lift of the Lefschetz thimble
$D_0$ (which reduces the space of choices to $\Z$, as expected since
vanishing cycles are only defined up to shifts). If we orient $D_0$ from $p_-$ towards $p_+$,
then the phase of $D_0$ (the function $\phi_{D_0}:D_0\to\R/2\pi\Z$ defined
by $\phi_{D_0}(p)=\arg \Omega(v)$ for any $p\in D_0$ and $v\in T_pD_0-\{0\}$
compatible with the orientation) has the property
that $$\phi_{D_0}(p_-)=\frac{\pi b}{a+b}\ \ \mathrm{and}\ \ \phi_{D_0}(p_+)=
\frac{\pi a}{a+b}.$$ Moreover, it is easy to check that
$\phi_{D_0}(p)\in (0,\pi)$ for all $p\in D_0$ (because
$\Omega=\frac{1}{b}d\log x$, and $\arg x$ is monotonically increasing
along $D_0$).
Hence, there exists a graded Lagrangian lift of $D_0$ for which
the phase function takes values in $(0,\pi),$ which means that we can
choose a graded lift of $L_0$ by setting
$$\tilde{\phi}_{L_0}(p_-)=\frac{\pi b}{a+b}\ \ \mathrm{and}
\ \ \tilde{\phi}_{L_0}(p_+)=\frac{\pi a}{a+b}.$$
Arguing similarly for the other vanishing cycles (or using the
$\Z/(a+b)$-equivariance), we can choose graded lifts of
$L_j=\{p_{j,-},p_{j,+}\}$ (where $\arg x_{j,\pm}=\frac{1}{a+b}(\pm \pi
b-2\pi j)$) by setting
$$\tilde{\phi}_{L_j}(p_{j,-})=\frac{\pi (b+2j)}{a+b}\ \ \mathrm{and}
\ \ \tilde{\phi}_{L_j}(p_{j,+})=\frac{\pi (a+2j)}{a+b}.$$
Now, the degree of the morphism $x_j$, corresponding to
$p_{j,+}=p_{j+a,-}\in L_j\cap L_{j+a}$, is given by the difference of
phases:
$$\deg x_j=\frac{1}{\pi}(\tilde{\phi}_{L_{j+a}}(p_{j+a,-})-
\tilde{\phi}_{L_j}(p_{j,+}))=\frac{b+2(j+a)}{a+b}-\frac{a+2j}{a+b}=1.$$
Similarly for $y_j$:
$$\deg y_j=\frac{1}{\pi}(\tilde{\phi}_{L_{j+b}}(p_{j+b,+})-
\tilde{\phi}_{L_j}(p_{j,-}))=\frac{a+2(j+b)}{a+b}-\frac{b+2j}{a+b}=1.$$
\endproof

Theorem \ref{thm:fswp1} now follows immediately from Lemmas
\ref{l:isectsp1}--\ref{l:gradingp1}; as in the case of weighted projective
planes, the only difference between the DG-categories $\FS(W,\{\gamma_i\})$
and $\Cc$ is that the objects of $\Cc$ are numbered ``backwards'', so the
equivalence of categories takes the objects $L_0,\dots,
L_{a+b-1}$ of $\FS(W,\{\gamma_i\})$ to the
objects $w_{a+b-1},\dots,w_0$ of $\Cc.$

\section{Mirrors of weighted projective planes}\label{sec:fsp2}

\subsection{The mirror Landau-Ginzburg model and its fiber $\Sigma_0$}
\label{ss:sigma0}
The mirror to the weighted projective
plane $\CP^2(a,b,c)$ is the affine hypersurface $X=\{x^ay^bz^c=1\}\subset
(\C^*)^3$, equipped with the superpotential $W=x+y+z$, and a symplectic form
$\omega$ that we leave unspecified for the moment. During most of the
argument, we will assume $\omega$ to be anti-invariant under complex conjugation
(which implies exactness) and under the diagonal action of the cyclic group
$\Z/(a+b+c)$, but these assumptions will be weakened at the end. Of course,
since $X$ is non-compact, we also need to choose $\omega$ in such a way as
to ensure that the Lefschetz thimbles and vanishing cycles considered below
are well-defined. It is easy to check that,
among many other possibilities, a symplectic form such as
$$\omega=i\sum_{i,j=1}^3
a_{ij} \frac{dz_i}{z_i}\wedge\frac{d\bar{z}_j}{\bar{z}_j}$$ (where $(a_{ij})$
is a positive definite Hermitian matrix, with real coefficients if we
require complex conjugation anti-invariance) generates a horizontal
distribution for which parallel transport is well-defined,
because, with respect to the induced K\"ahler metric, $X$ is complete and
the gradient vector of $W$ has norm bounded from below outside of a compact
set.

Topologically, $X$ is just a complex torus $(\C^*)^2$, at least if
$\delta=gcd(a,b,c)=1$; otherwise $X$ is disconnected, and each of its
$\delta$ components is a complex torus.

For each $\lambda\in\C$,
the fiber $\Sigma_\lambda=W^{-1}(\lambda)\subset X$ is an affine curve
given by the equation $x^ay^b(\lambda-x-y)^c=1$; this curve is smooth
unless $\lambda$ is one of the $a+b+c$ critical values of $W$.
We will view $\Sigma_\lambda$ as a branched
covering of $\C^*$, by projecting to the $x$ axis (this choice is arbitrary,
and we will occasionally use the symmetry between the variables $x,y,z$ in
the argument). For a generic value of $x\in\C^*$, the polynomial
$x^a y^b(\lambda-x-y)^c-1$ of degree $b+c$ in the variable $y$ admits
$b+c$ distinct simple roots; therefore, the projection
$\pi_x:\Sigma_\lambda\to \C^*$ is a $(b+c)$-fold covering. The branch
points of $\pi_x$ are those values of $x$ for which there is a double
root, i.e.\ a value of $y$ such that $P(y)=x^a y^b(\lambda-x-y)^c=1$ and
$P'(y)=0$. Since $$\frac{P'(y)}{P(y)}=\frac{b}{y}-\frac{c}{\lambda-x-y},$$
the condition $P'(y)=0$ implies that $cy=b(\lambda-x-y)$, i.e.\
$y=\frac{b}{b+c}(\lambda-x)$. Substituting into the equation of
$\Sigma_\lambda$, we obtain the equation
\begin{equation} \label{eq:branchx}
x^a (\lambda-x)^{b+c}=\frac{(b+c)^{b+c}}{b^b\,c^c}
\end{equation}
for the branch points of $\pi_x$. Since this is a polynomial equation of
degree $a+b+c$, for a generic value of $\lambda$ there are $a+b+c$ distinct
branch points, all of which are simple (i.e.\ isolated non-degenerate
critical points of $\pi_x$).

In the remainder of this section, we set $\lambda=0$, and describe the
curve $\Sigma_0$ in detail, by computing the
monodromy of the $(b+c)$-fold branched covering $\pi_x:\Sigma_0\to\C^*$
around the origin and around its $a+b+c$ branch points.

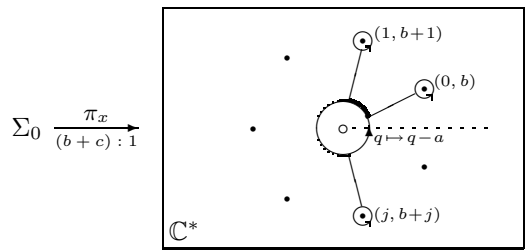
\begin{figure}[b]
\centering
\setlength{\unitlength}{0.8cm}
\begin{picture}(8,4.2)(-5,-2)
\put(-5,0){\makebox(0,0)[rc]{$\Sigma_0$}}
\put(-4.8,0){\vector(1,0){1.4}}
\put(-4.1,-0.1){\makebox(0,0)[ct]{\tiny $(b+c):1$}}
\put(-4.1,0.1){\makebox(0,0)[cb]{\small $\pi_x$}}
\put(-3,-2){\line(1,0){6}}
\put(-3,2){\line(1,0){6}}
\put(-3,-2){\line(0,1){4}}
\put(3,-2){\line(0,1){4}}
\put(0,0){\circle{0.15}}
\put(-2.9,-1.85){$\C^*$}
\put(1.35,0.65){\circle*{0.1}}
\put(1.35,-0.65){\circle*{0.1}}
\put(0.33,-1.46){\circle*{0.1}}
\put(0.33,1.46){\circle*{0.1}}
\put(-0.93,1.17){\circle*{0.1}}
\put(-0.93,-1.17){\circle*{0.1}}
\put(-1.5,0){\circle*{0.1}}
\multiput(0.15,0)(0.2,0){12}{\line(1,0){0.05}}
\put(0,0){\circle{0.8}}
\put(0.43,-0.05){\vector(0,1){0.1}}
\put(0.5,-0.23){\tiny $q\!\mapsto\!q\!-\!a$}
\put(0.41,0.2){\circle*{0.08}}
\put(0.39,0.18){\line(2,1){0.8}}
\put(1.35,0.65){\circle{0.3}}
\put(1.47,0.55){\line(-1,0){0.08}}
\put(1.47,0.55){\line(0,-1){0.08}}
\put(1.5,0.7){\tiny $(0,b)$}
\qbezier[30](0.43,0.21)(0.34,0.43)(0.10,0.48)
\put(0.10,0.48){\line(1,4){0.208}}
\put(0.33,1.46){\circle{0.3}}
\put(0.45,1.36){\line(-1,0){0.08}}
\put(0.45,1.36){\line(0,-1){0.08}}
\put(0.5,1.5){\tiny $(1,b\!+\!1)$}
\qbezier[30](0.405,0.195)(0.32,0.41)(0.095,0.455)
\qbezier[20](-0.095,0.455)(0,0.475)(0.095,0.455)
\qbezier[5](-0.425,0.20)(-0.32,0.41)(-0.095,0.455)
\qbezier[20](-0.095,-0.455)(0,-0.475)(0.095,-0.455)
\qbezier[5](-0.425,-0.20)(-0.32,-0.41)(-0.095,-0.455)
\put(0.095,-0.455){\line(1,-4){0.208}}
\put(0.33,-1.46){\circle{0.3}}
\put(0.45,-1.56){\line(-1,0){0.08}}
\put(0.45,-1.56){\line(0,-1){0.08}}
\put(0.5,-1.5){\tiny $(j,b\!+\!j)$}
\end{picture}
\caption{The projection $\pi_x:\Sigma_0\to \C^*$ (of degree $b+c$, with
$a+b+c$ branch points)}
\label{fig:sigma0}
\end{figure}

\begin{lemma}\label{l:sigma0}
The fiber of $\pi_x:\Sigma_0\to\C^*$ can be identified with $\Z/(b+c)$ in
such a way that the monodromy of $\pi_x$ around the origin in $\C^*$ is
given by $q\mapsto q-a$, and the monodromies around the $a+b+c$ branch
points are given by the transpositions $(j,j+b)$, $0\le j<a+b+c$
(see Figure~\ref{fig:sigma0}).
\end{lemma}

To understand this statement, first observe
that, when $x=\epsilon e^{i\theta}$ is close to $0$, the $b+c$ roots
of the equation
\begin{equation}\label{eq:s0}
x^a y^b (-x-y)^c=1
\end{equation} lie close to
those of the equation $$(-1)^c y^{b+c}=\epsilon^{-a} e^{-ia\theta}.$$ Hence,
we can choose an identification of the fiber of $\pi_x$ above a small real
positive value $x=\epsilon$ (or any other $\epsilon e^{i\theta}$ fixed in
advance) with the cyclic group $\Z/(b+c)$ in a manner
compatible with the cyclic ordering of the points. Moreover, varying $\theta$ from
$0$ to $2\pi$, we obtain that the monodromy of $\pi_x$ around
the origin is given by the translation $q\mapsto q-a$ in $\Z/(b+c)$ (i.e.,
the permutation sending the root $y_q$ of
$x^ay^b (-x-y)^c=1$ to $y_{q-a}$).

Next, consider a critical value of $\pi_x$, i.e.\ a
root $x_0$ of (\ref{eq:branchx}) for $\lambda=0$, and the
radial half-line $\ell$ through $x_0$, i.e.\ the set of all $x\in \C^*$ with
argument equal to $\theta_0=\arg x_0$. Moving $x$ along $\ell$ starting
from a point $x_*=\epsilon e^{i\theta_0}$ close to the origin, two of the
$b+c$ roots of (\ref{eq:s0}) become equal to each other
as $x$ approaches $x_0$; this determines the monodromy of $\pi_x$ around
$x_0$, namely a transposition in the symmetric group $S_{b+c}$ acting on
a fiber of $\pi_x$. We claim that, identifying the fiber $\pi_x^{-1}(x_*)$
with $\Z/(b+c)$ as above, this transposition
exchanges two elements $q_0$ and $q_0+b$. This can be seen as follows.

\begin{figure}[t]
\centering
\setlength{\unitlength}{0.8cm}
\begin{picture}(5,3.5)(-2.5,-1.5)
\put(-2.5,-1.5){\small $x\to 0$} 
\put(-1.5,0){\vector(1,0){3}}
\put(0,-1.5){\vector(0,1){3}}
\put(-1.29,-0.51){\circle*{0.08}}
\put(-1.29,0.51){\circle*{0.08}}
\put(-0.57,-1.23){\circle*{0.08}}
\put(-0.57,1.23){\circle*{0.08}}
\put(0.45,-1.23){\circle*{0.08}}
\put(0.45,1.23){\circle*{0.08}}
\put(1.17,-0.51){\circle*{0.08}}
\put(1.17,0.51){\circle*{0.08}}
\multiput(0.45,-1.23)(0,0.2){12}{\line(0,1){0.1}}
\put(1.55,0.1){\tiny Re\,$y$}
\put(0.1,1.55){\tiny Im\,$y$}
\end{picture}
\begin{picture}(5,3)(-3.5,-1.5)
\put(-3.5,-1.5){\small $x=x_0$} 
\put(-2.5,0){\vector(1,0){3}}
\put(0,-1.5){\vector(0,1){3}}
\put(-2.26,-0.28){\circle*{0.08}}
\put(-2.26,0.28){\circle*{0.08}}
\put(-1.73,-0.61){\circle*{0.08}}
\put(-1.73,0.61){\circle*{0.08}}
\put(-0.67,0){\circle*{0.08}}
\put(-0.67,0){\circle*{0.08}}
\put(0.19,-0.22){\circle*{0.08}}
\put(0.19,0.22){\circle*{0.08}}
\end{picture}
\begin{picture}(5,3)(-4,-1.5)
\put(-3.5,-1.5){\small $x\to+\infty$} 
\put(-3.5,0){\vector(1,0){4}}
\put(0,-1.5){\vector(0,1){3}}
\put(-3.14,-0.22){\circle*{0.08}}
\put(-3.14,0.22){\circle*{0.08}}
\put(-2.70,-0.43){\circle*{0.08}}
\put(-2.70,0.43){\circle*{0.08}}
\put(-2.31,0){\circle*{0.08}}
\put(-0.15,0){\circle*{0.08}}
\put(0.08,-0.12){\circle*{0.08}}
\put(0.08,0.12){\circle*{0.08}}
\multiput(-2.31,0.02)(0.2,0){11}{\line(1,0){0.1}}
\end{picture}
\caption{The roots of $x^a y^b (-x-y)^c=1$ for $x\in\R_+$ ($(a,b,c)=(1,3,5)$)}
\label{fig:rootsab}
\end{figure}
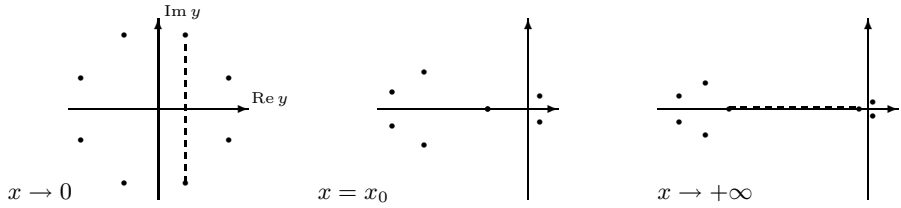

Assume for simplicity that $b+c$ is even and that
$x_0$ is the positive real root of
(\ref{eq:branchx}) for $\lambda=0$; the general case is handled similarly,
inserting factors $e^{i\theta_0}$ where needed. For $x\to 0$, as explained
above, the $b+c$ roots of (\ref{eq:s0}) are close to those of
$$y^{b+c}=(-1)^c x^{-a},$$ i.e.\ $b+c$ evenly spaced
points on a circle (Figure~\ref{fig:rootsab}, left).
As $x$ increases, two complex conjugate roots
$y,\bar{y}$ approach the real axis and eventually become equal for $x=x_0$
(Figure~\ref{fig:rootsab}, center),
so that there are two additional real roots for $x>x_0$. As $x\to
+\infty$, the roots of (\ref{eq:s0}) are divided into two groups,
$b$ roots close to the origin, approximated by those of $$y^b=(-1)^c
x^{-(a+c)},$$ and $c$ roots close to $-x$, corresponding to values of
$z=-x-y$ close to the origin and approximated by the roots of $$z^c=(-1)^b
x^{-(a+b)}$$ (Figure~\ref{fig:rootsab}, right). Hence, identifying the fiber
of $\pi_x$ for $x$ small with $\Z/(b+c)$ in a manner compatible with the
cyclic ordering, the two points which merge for
$x=x_0$ (the vanishing cycle of $\pi_x$ at $x_0$) differ from each other
by exactly $b$ (this can also be checked by numerical experimentation).

The above argument gives us that the monodromy around one of the
branch points $x_0$ of $\pi_x$, e.g.\ the branch point located on
the positive real axis or immediately above it, is a transposition
$(q_0,q_0+b)$; changing the identification between the reference fiber
of $\pi_x$ above $x_*$ and the cyclic group
$\Z/(b+c)$ if necessary, we can assume that $q_0=0$.

We now find the monodromy around the other branch points of $\pi_x$. For this
purpose, observe that the group $G=\Z/(a+b+c)$ acts on $X$ by
$(x,y,z)\mapsto (x\zeta^j,y\zeta^j,z\zeta^j)$, where $\zeta=\exp(\frac{2\pi
i}{a+b+c})$, and that this action preserves $\Sigma_0$, mapping the fiber
of $\pi_x$ above $x$ to the fiber above $x\zeta^j$. Hence, denoting by
$y',y''$ the two points of the fiber above $x_*=\epsilon e^{i\theta_0}$ which
converge to each other as $x$ moves radially outwards to $x_0$ (those
labelled $0$ and $b$), we know that the two points of the fiber above $x_*
\zeta^j$ which converge to each other as $x$ moves radially
outwards to $x_0\zeta^j$ are $y'\zeta^j$ and $y''\zeta^j$. We now transport
these two values of $y$ from the fiber $\pi_x^{-1}(x_*\zeta^j)$ to
$\pi_x^{-1}(x_*)$ along the arc $x(t)=x_*e^{2\pi it}$
for $t\in [0,\frac{j}{a+b+c}]$. Approximating the $b+c$ points of
$\pi_x^{-1}(\epsilon e^{i\theta})$ by the roots of
$(-1)^c y^{b+c}=\epsilon^{-a}e^{-ia\theta}$, the parallel transport along
the considered arc induces a multiplication by
$\exp(2\pi i \frac{a}{b+c}\frac{j}{a+b+c})$.
Observing that $$\zeta^j \exp(2\pi i\tfrac{j\,a}{(b+c)(a+b+c)})=
\exp(2\pi i\tfrac{j}{b+c}),$$ we obtain that the two points of
$\pi_x^{-1}(x_*)$ which become equal as $x$ is moved first counterclockwise
around the origin and then radially outwards to $x_0\zeta^j$ are those which
correspond to the elements $j$ and $b+j$ of $\Z/(b+c)$. Hence, the monodromy
of $\pi_x$ around $x_0\zeta^j$ (joining $x_*$ to $x_0\zeta^j$ in the
prescribed way) is the transposition $(j,b+j)$, which completes the proof of
Lemma~\ref{l:sigma0}. By the way, remark that the comparison
between the values $j=0$ and $j=a+b+c$ is consistent with our
determination of the monodromy around $x=0$.

\subsection{The vanishing cycles}\label{ss:vcs}
Now that the fiber $\Sigma_0$ is
well-understood, we compute the vanishing cycles of the Lefschetz
fibration $W:X\to\C$ by studying the degeneration of $\Sigma_\lambda$
as $\lambda$ approaches a critical value of $W$.

The curve $\Sigma_\lambda$ becomes singular when two branch points
of the projection $\pi_x:\Sigma_\lambda\to\C^*$
merge with each other, giving rise to a nodal point. This occurs
whenever (\ref{eq:branchx}) admits a double root. Considering the
logarithmic derivative of the left-hand side,
we obtain the relation $\frac{a}{x}-\frac{b+c}{\lambda-x}=0$, which leads to
$x=\frac{a}{a+b+c}\lambda$ for a double root of (\ref{eq:branchx}), and
substituting we obtain the equation
\begin{equation} \label{eq:critw}
\lambda^{a+b+c}=\frac{(a+b+c)^{a+b+c}}{a^a\,b^b\,c^c}
\end{equation}
for the $a+b+c$ critical values of $W$ (this equation can also be obtained
directly).

For symmetry and for simplicity, we will choose the smooth curve
$\Sigma_0=W^{-1}(0)$ as our reference fiber of the Lefschetz fibration
$W:X\to\C$, and we will choose straight lines for the arcs $\gamma_j$
joining the origin to the various critical values $\lambda_j=\lambda_0
\zeta^{-j}$ of $W$ ($0\le j<a+b+c$), where $\lambda_0$ is
the real positive root of (\ref{eq:critw}) and
$\zeta=\exp(\frac{2\pi i}{a+b+c})$.
Hence, in order to construct the category of Lagrangian vanishing cycles of $W$, we need
to understand how the smooth fiber $\Sigma_0$ above the reference point
$0$ degenerates to the nodal curve $\Sigma_{\lambda_j}$ when $\lambda$ moves
radially from $0$ to $\lambda_j$.

\begin{figure}[b]
\centering
\setlength{\unitlength}{0.8cm}
\begin{picture}(5,3.2)(-2.5,-1.5)
\put(-2.5,-1.5){\small $\lambda=0$}
\put(-1.5,0){\vector(1,0){3}}
\put(0,-1.5){\vector(0,1){3}}
\put(-1.31,0){\circle*{0.08}}
\put(-0.81,-1.02){\circle*{0.08}}
\put(-0.81,1.02){\circle*{0.08}}
\put(0.29,-1.28){\circle*{0.08}}
\put(0.29,1.28){\circle*{0.08}}
\put(1.18,-0.57){\circle*{0.08}}
\put(1.18,0.57){\circle*{0.08}}
\qbezier[50](0.29,-1.28)(0.7,0)(0.29,1.28)
\put(1.55,0.1){\tiny Re\,$x$}
\put(0.1,1.55){\tiny Im\,$x$}
\put(0.5,-1){\small $\delta_0$}
\put(0,0){\circle{0.1}}
\end{picture}
\begin{picture}(5,3)(-1.5,-1.5)
\put(0.5,-1.5){\small $\lambda=\lambda_0$} 
\put(-1,0){\vector(1,0){4.2}}
\put(0,-1.5){\vector(0,1){3}}
\put(-0.66,0){\circle*{0.08}}
\put(-0.14,-0.72){\circle*{0.08}}
\put(-0.14,0.72){\circle*{0.08}}
\put(1.48,0){\circle*{0.08}}
\put(1.48,0){\circle*{0.08}}
\put(2.78,-0.35){\circle*{0.08}}
\put(2.78,0.35){\circle*{0.08}}
\end{picture}
\begin{picture}(6,3)(-1.5,-1.5)
\put(0.5,-1.5){\small $\lambda\to+\infty$} 
\put(-1,0){\vector(1,0){5}}
\put(0,-1.5){\vector(0,1){3}}
\put(-0.59,0){\circle*{0.08}}
\put(-0.09,-0.64){\circle*{0.08}}
\put(-0.09,0.64){\circle*{0.08}}
\put(0.74,0){\circle*{0.08}}
\put(2.69,0){\circle*{0.08}}
\put(3.41,-0.29){\circle*{0.08}}
\put(3.41,0.29){\circle*{0.08}}
\end{picture}
\caption{The branch points of $\pi_x$ for $\lambda\in\R_+$ ($(a,b,c)=(4,2,1)$)}
\label{fig:vc0}
\end{figure}
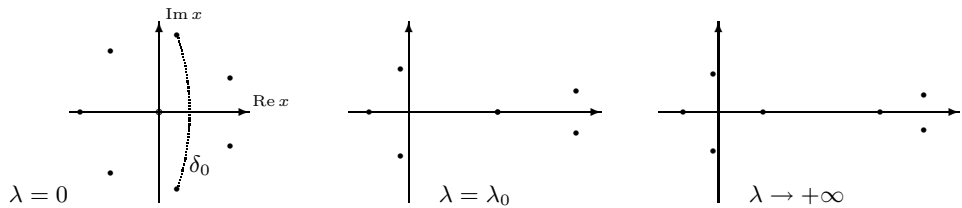

We first consider the motion of the branch points of $\pi_x$ as $\lambda$
increases along the positive real axis from $0$ to the critical value
$\lambda_0$. For each value of $\lambda$, the $a+b+c$ branch points are
given by the roots of (\ref{eq:branchx}). When $\lambda=0$, they all
lie on a circle centered at the origin, as represented in
Figure~\ref{fig:sigma0}. As $\lambda\to\lambda_0$, two complex conjugate
branch points converge to each other, so that for $\lambda=\lambda_0$ the
equation (\ref{eq:branchx}) has a double root $x=\frac{a}{a+b+c}\lambda_0$
on the positive real axis (Figure~\ref{fig:vc0}, center). Finally, for
$\lambda\to +\infty$, the roots of
(\ref{eq:branchx}) split into two groups, one of $a$ points close to the
origin that can be approximated by the roots of
$x^a=K_{b,c}\lambda^{-(b+c)}$ (where $K_{b,c}=b^{-b}c^{-c}(b+c)^{b+c}$),
and one of $b+c$ points close to $\lambda$ for which $\xi=\lambda-x$ can
be approximated by the roots of $\xi^{b+c}=K_{b,c}\lambda^{-a}$
(Figure~\ref{fig:vc0}, right). Hence, it
can be checked that the two branch points of $\pi_x:\Sigma_0\to\C^*$ which
merge for $\lambda\to\lambda_0$ are those with argument $\arg x=
\pm \frac{b+c}{a+b+c}\pi$, and that the projection to $\C^*$ of the
corresponding vanishing cycle is an arc $\delta_0$ which is symmetric
with respect to the real axis, intersects it only once in its
positive part, and remains everywhere inside the circle containing the
critical values of $\pi_x$ (Figure~\ref{fig:vc0}, left).

More precisely, the above discussion gives us a {\it topological}
description of the vanishing cycle $L_0\subset \Sigma_0$, up to homotopy.
Namely, two of the $b+c$ lifts to $\Sigma_0$ of the arc
$\delta_0\subset\C^*$ have common end points (the ramification points of
$\pi_x$ lying above the end points of $\delta_0$), and their union forms
a closed loop $L'_0$ in $\Sigma_0$. This loop is a topological vanishing
cycle, i.e.\ it shrinks to a point in $\Sigma_\lambda$ when $\lambda\to
\lambda_0$, but a priori it is only homotopic to the symplectic vanishing
cycle $L_0$ (obtained by parallel transport using the symplectic connection).

The actual position of the vanishing cycle $L_0$ inside $\Sigma_0$ depends on
the choice of the symplectic form $\omega$ on $X$; for a given $\omega$ it
can be calculated numerically (and it can be checked that for ``reasonable''
choices of $\omega$, $L_0$ and $L'_0$ intersect all other vanishing cycles
in the same manner). However, this calculation is
unnecessary for our purposes.
Indeed, if we endow $X$ with a symplectic form that is anti-invariant
by complex conjugation, then the vanishing cycle
$L_0$ is invariant by complex conjugation, i.e.\ complex conjugation maps
$L_0$ to itself in an orientation-preserving manner, and the same is true
of $L'_0$. Since $L_0$ and $L'_0$ are homotopic to each other in $\Sigma_0$,
their (oriented) invariance under complex conjugation is
sufficient to imply that they are {\it Hamiltonian isotopic}, which
means that for the purpose of determining categories of vanishing cycles, $L_0$ and
$L'_0$ are interchangeable.

If we deform $\omega$ to a non-exact form, complex conjugation invariance is
lost. The intersection patterns between vanishing cycles remain the
same for small deformations (and can be forced to remain the same even for
large deformations by performing suitable Hamiltonian isotopies), but
the calculation of the coefficient assigned to a given pseudo-holomorphic
curve involves its symplectic area and hence requires one to work with the
actual vanishing cycles rather than their topological approximations.
Hence, we may obtain non-trivial deformations of the category of vanishing
cycles; however, these deformations only amount to modifications of the structure
constants of the products $m_k$, rather than changes in the Floer complexes
themselves or in the types of pseudo-holomorphic curves that may arise.

In any case, except at the very end of the argument, we will always be
considering symplectic forms that are anti-invariant under complex conjugation,
in which case the approximation of $L_0$ by $L'_0$ is legitimate.

We now consider the other vanishing cycles $L_j$ of the Lefschetz fibration
$W$. Recall that the group $G=\Z/(a+b+c)$ acts on $X$, in a manner that
preserves $\Sigma_0$; moreover, $W:X\to\C$ is $G$-equivariant. If we assume
the symplectic form $\omega$ to be $G$-invariant, the symplectic connection
and the associated parallel transport will also be $G$-equivariant.
Therefore, since the arc $\gamma_j\subset \C$ joining the origin to
$\lambda_j=\lambda_0\zeta^{-j}$ is the image of $\gamma_0$ by the action of
$\zeta^{-j}$ (where $\zeta=\exp(\frac{2\pi i}{a+b+c})$), the same is true
of the corresponding Lefschetz thimbles, and hence of the vanishing cycles
in $\Sigma_0$. This gives us a description of $L_j$ for all values of $j$.
As in the case of $L_0$, we will consider, rather than $L_j$ itself, a loop
$L'_j\subset\Sigma_0$ which is homotopic to $L_j$ and can be obtained as
a double lift via $\pi_x:\Sigma_0\to\C^*$ of an embedded arc
$\delta_j\subset \C^*$. The loop $L'_j$ is defined to be the image of $L'_0$
by the action of $\zeta'_j$, which means that $\delta_j$ is the image of
$\delta_0$ by a rotation of angle $-\frac{2\pi\,j}{a+b+c}$. If, in addition
to its $G$-invariance, $\omega$ is assumed to be anti-invariant under complex
conjugation, then $L'_j$ is Hamiltonian isotopic to $L_j$, so we can work
with $L'_j$ instead of $L_j$.

\begin{figure}[ht]
\centering
\setlength{\unitlength}{1cm}
\begin{picture}(4,3)(-2,-1.6)
\put(-1.31,0){\circle*{0.08}}
\put(-0.81,-1.02){\circle*{0.08}}
\put(-0.81,1.02){\circle*{0.08}}
\put(0.29,-1.28){\circle*{0.08}}
\put(0.29,1.28){\circle*{0.08}}
\put(1.18,-0.57){\circle*{0.08}}
\put(1.18,0.57){\circle*{0.08}}
\qbezier[150](0.29,-1.28)(0.7,0)(0.29,1.28)
\qbezier[150](1.18,-0.57)(0.436,0.547)(-0.81,1.02)
\qbezier[150](1.18,0.57)(0.436,-0.547)(-0.81,-1.02)
\qbezier[150](1.18,0.57)(-0.155,0.68)(-1.31,0)
\qbezier[150](1.18,-0.57)(-0.155,-0.68)(-1.31,0)
\qbezier[150](0.29,1.28)(-0.63,0.30)(-0.81,-1.02)
\qbezier[150](0.29,-1.28)(-0.63,-0.30)(-0.81,1.02)
\put(0.55,-1.05){\makebox(0,0)[cc]{\small $\delta_0$}}
\put(1.18,-0.22){\makebox(0,0)[cc]{\small $\delta_6$}}
\put(0.90,0.77){\makebox(0,0)[cc]{\small $\delta_5$}}
\put(-0.07,1.18){\makebox(0,0)[cc]{\small $\delta_4$}}
\put(-0.95,0.71){\makebox(0,0)[cc]{\small $\delta_3$}}
\put(-1.14,-0.30){\makebox(0,0)[cc]{\small $\delta_2$}}
\put(-0.47,-1.10){\makebox(0,0)[cc]{\small $\delta_1$}}
\put(-2,-1.7){\small $(a,b,c)=(4,2,1)$}
\put(0,0){\circle{0.1}}
\end{picture}
\begin{picture}(4,2.6)(-2,-1.6)
\put(1,0){\circle*{0.08}}
\put(-0.5,0.866){\circle*{0.08}}
\put(-0.5,-0.866){\circle*{0.08}}
\qbezier[180](-0.5,0.866)(1,0)(-0.5,-0.866)
\qbezier[180](-0.5,0.866)(-0.5,-0.866)(1,0)
\qbezier[180](1,0)(-0.5,0.866)(-0.5,-0.866)
\put(0.75,0.4){\makebox(0,0)[cc]{\small $\delta_2$}}
\put(-0.05,-0.9){\makebox(0,0)[cc]{\small $\delta_0$}}
\put(-0.7,0.5){\makebox(0,0)[cc]{\small $\delta_1$}}
\put(-1.6,-1.7){\small $(a,b,c)=(1,1,1)$}
\put(0,0){\circle{0.1}}
\end{picture}
\caption{The vanishing cycles $L_j\subset\Sigma_0$}
\label{fig:vcs}
\end{figure}
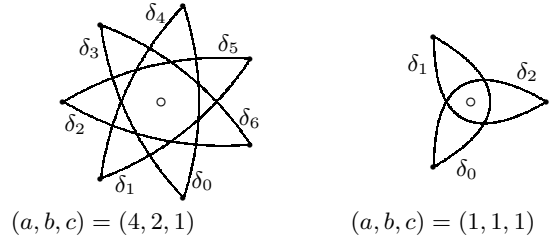

Hence, to summarize the above discussion, we have the following lemma:

\begin{lemma}\label{l:vcs}
The vanishing cycles $L_j\subset\Sigma_0$ $(0\le j<a+b+c)$ are homotopic
(and, if $\omega$ is invariant under the action of $\Z/(a+b+c)$ and
anti-invariant under complex conjugation, Hamiltonian isotopic) to closed loops
$\smash{L'_j}\subset\Sigma_0$ which
project by $\pi_x$ to arcs $\delta_j\subset\C^*$ as represented in Figure
\ref{fig:vcs} (the end points of $\delta_j$ are the branch
points of $\pi_x$ for which $\arg x=-2\pi\frac{j}{a+b+c}\pm \pi
\frac{b+c}{a+b+c})$.
\end{lemma}

In the following sections, we assume that $\omega$ is $\Z/(a+b+c)$-invariant
and anti-invariant under complex conjugation, and we implicitly identify $L_j$
with $L'_j$.

\subsection{The Floer complexes}\label{ss:isects}
The objects of the category $\FS(W,\{\gamma_j\})$ are described
by Lemma \ref{l:vcs}; we now determine its morphisms by studying the intersections
between the closed loops $L_j\subset \Sigma_0$. This simply involves
looking carefully at Figures \ref{fig:sigma0} and \ref{fig:vcs} in order
to determine, among the intersections between $\delta_i$ and $\delta_j$,
which ones lift to intersections between $L_i$ and $L_j$.

\begin{lemma}\label{l:isects}
The direct sum $\bigoplus_{i<j} CF^*(L_i,L_j)$ is a free module
of total rank $3(a+b+c)$ over the coefficient ring, generated by
the following intersection points:
\medskip

\begin{tabular}{lllll}
$x_i\in CF^*(L_i,L_{i+a})$ & $(0\le i<b+c)$, & \qquad &
$\bar{x}_i\in CF^*(L_i,L_{i+b+c})$ & $(0\le i<a)$,\\
$y_i\in CF^*(L_i,L_{i+b})$ & $(0\le i<a+c)$, & &
$\bar{y}_i\in CF^*(L_i,L_{i+a+c})$ & $(0\le i<b)$,\\
$z_i\in CF^*(L_i,L_{i+c})$ & $(0\le i<a+b)$, & &
$\bar{z}_i\in CF^*(L_i,L_{i+a+b})$ & $(0\le i<c)$.
\end{tabular}
\medskip

\noindent Moreover, the Floer differential is trivial, i.e.\ $m_1=0$.
\end{lemma}

To determine $CF^*(L_i,L_j)$ for given $0\le i<j<a+b+c$, one must look for
intersection points between the projected arcs $\delta_i$ and
$\delta_j$. The arcs $\delta_i$ and $\delta_j$ intersect only if $j-i\le
b+c$ or $j-i\ge a$; in all other cases, $\delta_i\cap \delta_j=\emptyset$
and hence $CF^*(L_i,L_j)=0$. More precisely, $\delta_i\cap \delta_j$
contains one point if $j-i\le b+c$, and one point if $j-i\ge a$; if both
conditions hold simultaneously, then $|\delta_i\cap\delta_j|=2$ (see Lemma
\ref{l:vcs} and Figure \ref{fig:vcs}). Moreover,
if equality holds ($j-i=b+c$ or $j-i=a$), then the corresponding intersection
occurs at an end point of $\delta_i$ and $\delta_j$, i.e.\ a branch point of
$\pi_x$. In this case, the intersection of $\delta_i$ and $\delta_j$ always
lifts to a transverse intersection of $L_i$ and $L_j$, at the corresponding
critical point of $\pi_x$; this accounts for the generators $x_i$ and
$\bar{x}_i$ mentioned in the statement of Lemma \ref{l:isects}.

When $j-i<b+c$ or $j-i>a$, we need to consider the structure of the branched
covering $\pi_x$ in order to determine whether intersections between
$\delta_i$ and $\delta_j$ lift to intersections between $L_i$ and $L_j$.
Call $p_i$ the branch point of $\pi_x$ with argument
$\arg x=-2\pi\frac{j}{a+b+c}-\pi\frac{b+c}{a+b+c}$, which is an end point
of $\delta_i$, and define similarly $p_j$. When $j-i<b+c$, consider the
corresponding intersection point $q\in \delta_i\cap \delta_j$, and use the
arcs joining $p_j$ to $q$ in $\delta_j$ and $q$ to $p_i$ in $\delta_i$ to
define an arc $\eta\subset\C^*$ joining $p_j$ to $p_i$, with a rotation
angle of $2\pi\frac{j-i}{a+b+c}$ around the origin. It follows from Lemma
\ref{l:sigma0} (cf.\ also Figure \ref{fig:sigma0}) that, over a neighborhood
of $\eta$, we can consistently label the sheets of the covering $\pi_x$ by
elements of $\Z/(b+c)$, in such a way that the monodromies around the branch
points $p_i$ and $p_j$ are transpositions of the form $(k_i,k_i+b)$ and
$(k_j,k_j+b)$, with $k_i-k_j=j-i$. Hence, near the point $q$, the vanishing
cycle $L_i$ lies in the two sheets of $\pi_x$ labelled $k_i$ and $k_i+b$,
and similarly for $L_j$; the intersections of $L_i$ with $L_j$ above $q$
correspond to the elements of $\{k_i,k_i+b\}\cap \{k_j,k_j+b\}$. Since
$0<k_i-k_j=j-i<b+c$, this intersection is empty unless $k_i=k_j+b\mod b+c$,
i.e.\ $j-i=b$, which corresponds to the generator $y_i$ of the Floer
complex, or $k_j=k_i+b\mod b+c$, i.e.\ $j-i=c$, which corresponds to the
generator $z_i$. When $j-i>a$, one proceeds similarly,
introducing an arc in $\C^*$ joining $p_j$ to $p_i$ through the
relevant intersection point $q'$ of $\delta_i$ with $\delta_j$, with a rotation
angle of $2\pi(\frac{j-i}{a+b+c}-1)$ around the origin.
The sheets of $\pi_x$ containing $L_i$ and $L_j$ above the intersection point
$q'$ are now labelled $k'_i,k'_i+b$ and $k'_j,k'_j+b$, with $k'_i$ and
$k'_j$ two constants in $\Z/(b+c)$ such that
$k'_i-k'_j=j-i-(a+b+c)=j-i-a\mod b+c$. Therefore, the two cases where $L_i$
and $L_j$ intersect above $q'$ are when $i+j=a+b$, which corresponds to the
generator $z'_i$ of the Floer complex, and when $i+j=a+c$, which corresponds
to $y'_i$.

At this point it is worth observing that, for generic values of $(a,b,c)$,
each Floer complex $CF^*(L_i,L_j)$ has total rank at most one, so that the
Floer differential is necessarily zero. However, for specific values of
$(a,b,c)$ we may have numerical coincidences leading to more than one
intersection between two vanishing cycles; the most striking example is
that of the usual projective plane, $(a,b,c)=(1,1,1)$, for which $|L_i\cap
L_j|=3$ $\forall i<j$ (cf.\ Figure~\ref{fig:vcs}). Nonetheless, even in
these cases, the Floer differential vanishes, because $L_i$ and $L_j$
always realize the minimal geometric intersection number between closed
loops in their homotopy classes, as can be checked by
enumerating the various posible cases. This minimality of
intersection implies that $\Sigma_0$ contains no non-constant immersed disc
with boundary in $L_i\cup L_j$, and hence that the Floer differential
vanishes.

Another way to prove the vanishing of the Floer differential is to endow
$\Sigma_0$ and $\C^*$ with almost-complex structures which make the
projection $\pi_x$ holomorphic, and to observe that the projection to $\C^*$
of a pseudo-holomorphic disc in $\Sigma_0$ with boundary in $L_i\cup L_j$
is a pseudo-holomorphic disc in $\C^*$ with boundary in $\delta_i\cup
\delta_j$. If $|\delta_i\cap \delta_j|=1$, the maximum principle implies
that the projected pseudo-holomorphic disc is a constant map, and hence that
the disc in $\Sigma_0$ is contained in a fiber of $\pi_x$, which implies
that it is also constant. If $|\delta_i\cap \delta_j|=2$, one reaches the
same conclusion by observing the respective positions of the two intersection
points in $\C^*$ (a non-constant disc would have to pass through the origin).
As before, one concludes that the absence of non-trivial
pseudo-holomorphic discs makes the Floer differential identically zero,
which completes the proof of Lemma \ref{l:isects}.

\subsection{The product structures}\label{ss:products}
The aim of this section is to prove the following results concerning the
category $\FS(W,\{\gamma_j\})$:

\begin{lemma}\label{l:m3}
The higher products $m_k$ $(k\ge 3)$ are all identically zero.
\end{lemma}

\begin{lemma}\label{l:m2}
There exist non-zero constants $\alpha_{uv,i}$ such that\medskip

\begin{tabular}{clc}
$m_2(x_i,y_{i+a})=\alpha_{xy,i}\,\bar{z}_i$,& &
$m_2(x_i,z_{i+a})=\alpha_{xz,i}\,\bar{y}_i$, \\
$m_2(y_i,z_{i+b})=\alpha_{yz,i}\,\bar{x}_i$,& &
$m_2(y_i,x_{i+b})=\alpha_{yx,i}\,\bar{z}_i$, \\
$m_2(z_i,x_{i+c})=\alpha_{zx,i}\,\bar{y}_i$,& &
$m_2(z_i,y_{i+c})=\alpha_{zy,i}\,\bar{x}_i$.
\end{tabular}\medskip

\noindent
All other compositions (except those involving identity morphisms) vanish.
\end{lemma}

These results follow from a careful observation of the boundary structure
of a pseudo-holomorphic disc in $\Sigma_0$ with boundary in $\bigcup L_j$.
Endow $\Sigma_0$ with any almost-complex structure,
and let $u:D^2\to \Sigma_0$ be a pseudo-holomorphic map from the disc with
$k+1\ge 3$ marked points on its boundary
to $\Sigma_0$, mapping each segment on the boundary to an arc in one of the
Lagrangian submanifolds $L_j$. Each ``corner'' of the image of $u$
corresponds to an intersection point between two of the vanishing cycles,
and as such it corresponds to a generator of the Floer complex.

In accordance with Lemma \ref{l:isects}, we can classify the generators of the Floer
complex into three families, those of {\it type} $x$ (corresponding to
generators $x_i,\bar{x}_i$), those of {\it type} $y$ (generators $y_i,
\bar{y}_i$), and those of {\it type} $z$ (generators $z_i,\bar{z}_i$).
Moreover, observe that the total intersection of each $L_i$ with
all other vanishing cycles consists of $6$ points, two of each type:
depending on the value of $i$, $L_i$ is either the source of the morphism
$x_i$ or the target of $\bar{x}_{i-b-c}$, and it is either the source of
$\bar{x}_i$ or the target of $x_{i-a}$; similarly for types $y$ and $z$.

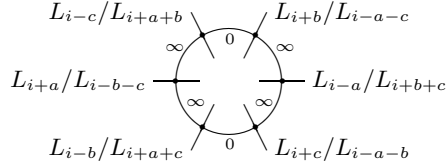
\begin{figure}[t]
\centering
\setlength{\unitlength}{1cm}
\begin{picture}(4,2.1)(-2,-1)
\put(0,0){\circle{1.42}}
\put(-0.71,0){\circle*{0.08}}
\put(0.71,0){\circle*{0.08}}
\put(0.35,0.61){\circle*{0.08}}
\put(-0.35,0.61){\circle*{0.08}}
\put(-0.35,-0.61){\circle*{0.08}}
\put(0.35,-0.61){\circle*{0.08}}
\put(-1,0){\line(1,0){0.6}}
\put(1,0){\line(-1,0){0.6}}
\put(0.2,0.31){\line(1,2){0.3}}
\put(-0.2,0.31){\line(-1,2){0.3}}
\put(0.2,-0.31){\line(1,-2){0.3}}
\put(-0.2,-0.31){\line(-1,-2){0.3}}
\put(1.1,0){\makebox(0,0)[lc]{\small $L_{i-a}/L_{i+b+c}$}}
\put(-1.1,0){\makebox(0,0)[rc]{\small $L_{i+a}/L_{i-b-c}$}}
\put(0.6,0.91){\makebox(0,0)[lc]{\small $L_{i+b}/L_{i-a-c}$}}
\put(-0.6,0.91){\makebox(0,0)[rc]{\small $L_{i-c}/L_{i+a+b}$}}
\put(0.6,-0.91){\makebox(0,0)[lc]{\small $L_{i+c}/L_{i-a-b}$}}
\put(-0.6,-0.91){\makebox(0,0)[rc]{\small $L_{i-b}/L_{i+a+c}$}}
\put(-0.85,0.4){\tiny $\infty$}
\put(-0.05,0.5){\tiny $0$}
\put(0.62,0.4){\tiny $\infty$}
\put(0.34,-0.33){\tiny $\infty$}
\put(-0.57,-0.33){\tiny $\infty$}
\put(-0.05,-0.9){\tiny $0$}
\end{picture}
\caption{The intersections of $L_i$ with the other vanishing cycles}
\label{fig:vcinters}
\end{figure}

The manner in which these points are arranged along the loop $L_i$ can be
seen easily by looking at Figure \ref{fig:vcs} and recalling the discussion
in the previous section.
Recall that $L_i$ passes through two branch points of $\pi_x$, which split
it into two halves (lifts of $\delta_i$ lying in different sheets of
$\pi_x$). One of these branch points corresponds to $x_i$ or
$\bar{x}_{i-b-c}$, while the other corresponds to $\bar{x}_i$ or $x_{i-a}$.
In between them, we have, on one half of $L_i$, one intersection of type $y$
(either $y_i$ or $\bar{y}_{i-a-c}$) and one of type $z$ (either $\bar{z}_i$
or $z_{i-c}$); on the other half of $L_i$, we have similarly one
intersection of type $y$ (either $\bar{y}_i$ or $y_{i-b}$) and one of
type $z$ (either $z_i$ or $\bar{z}_{i-a-b}$).
This structure is summarized in Figure \ref{fig:vcinters}.

An important property is that, for every one of the six portions of
$L_i$ delimited by these intersection points, one of the two immediately
adjacent components of $\Sigma_0-\bigcup L_j$ (on either side of $L_i$)
is unbounded (it is denoted by $0$ or $\infty$ on Figure
\ref{fig:vcinters} depending on whether its image under $\pi_x$ contains
the origin or the point at infinity in $\C^*$). These unbounded components
form an alternating pattern around $L_i$, changing side (left or right)
every time one of the intersection points is crossed.

On the other hand, the image of the pseudo-holomorphic map
$u$ may not intersect any of the unbounded components of $\Sigma_0-\bigcup
L_j$, because otherwise the maximum principle would imply that the image of
$u$ is unbounded. This imposes very strong constraints on the behavior of $u$
along the boundary of the disc. Namely, consider two consecutive marked
points (``corners''), such that the portion of boundary (``edge'')
in between them is mapped to an arc $\eta$ (oriented according to the
boundary orientation of the unit disc) contained in the vanishing cycle
$L_i$. Then, $\eta$ is exactly one of the six portions of $L_i$
delimited by its intersections with the other vanishing cycles, and
its orientation is determined by the requirement that the component
of $\Sigma_0-\bigcup L_j$ immediately to the left of $\eta$ be bounded
(see Figure \ref{fig:vcinters}). Moreover, the local behavior of $u$ at
an end point $p$ of $\eta$ is ``convex'', i.e.\ $u$ locally maps into only
one of the four regions delimited locally by the two vanishing cycles meeting
at $p$. In other words, the boundary of $\mathrm{Im}(u)$ is an oriented
piecewise smooth curve $\theta\subset\bigcup L_j$ which always turns left at every
intersection point it encounters. This boundary behavior has several important
consequences.

\begin{lemma}\label{l:3corners}
Among three consecutive corners of the image of $u$, there is always exactly
one of each type $x,y,z$.
\end{lemma}

\proof
Observe that two consecutive corners of the image of $u$ are necessarily of
different types (because two adjacent intersections of $L_i$ with other
vanishing cycles are always of different types).
Let $p,q,r$ be three consecutive
corners of the image of $u$, such that the edge from $p$ to $q$ lies in a
vanishing cycle $L_i$ and the edge from $q$ to $r$ lies in a vanishing cycle
$L_j$. The knowledge of the types of the points $p$ and $q$ completely
determines them, which in turn determines the type of $r$. For example, if
$p$ is of type $y$ and $q$ is of type $z$, then on the diagram of Figure
\ref{fig:vcinters} the edge joining them is the lowermost portion of $L_i$;
in particular the edge from $p$ to $q$ is adjacent to an unbounded component
of $\Sigma_0$ whose image under $\pi_x$ contains the origin. Considering the
intersection diagram for $L_j$ (similar to Figure~\ref{fig:vcinters}), the
point $q$ can be located by comparison with the diagram for $L_i$ (in our
example, $q$ is the point to the upper left of the diagram). Moreover, the
direction from which $\theta$ reaches $q$ can be determined by identifying
the unbounded component to which it is adjacent (in our example, the component
whose image under $\pi_x$ contains the origin, so $\theta$ reaches $q$ from
the innermost side of the diagram); since $\theta$ turns left at $q$, this
determines the edge from $q$ to $r$ and hence the type of $r$ (in our
example, $r$ is the left-most point on the intersection diagram, and hence
of type $x$). It can be checked easily that in all six cases, the type of
$r$ is different from those of $p$ and $q$.
\endproof

Next, recall that
by definition the successive edges of the image of $u$ lie inside
vanishing cycles $L_{i_0},L_{i_1},\dots,L_{i_k}$ with $i_0<i_1<\dots<i_k$
(see Definition \ref{def:fs}), and observe that following $\theta$ at a corner
of $u$ leads from a vanishing cycle $L_i$ to another vanishing cycle $L_j$,
with $i<j$ if and only if the intersection point is $x_i$, $y_i$ or $z_i$,
and $i>j$ if and only if the intersection point is $\bar{x}_j$, $\bar{y}_j$
or $\bar{z}_j$ (see Figure \ref{fig:vcinters}). Therefore, all corners of
$u$ but one correspond to generators of the Floer complexes among
$\{x_i,y_i,z_i\}$, while the last corner (between the edge on $L_{i_k}$ and
the edge on $L_{i_0}$) correspond to a generator among $\{\bar{x}_i,
\bar{y}_i,\bar{z}_i\}$.

With this observation, Lemma \ref{l:m3} follows immediately from Lemma
\ref{l:3corners}. Indeed, assume that there exists a pseudo-holomorphic
map $u$ from a disc with $k+1$ marked points to $\Sigma_0$, with edges
lying in vanishing cycles
$L_{i_0},L_{i_1},\dots,L_{i_k}$ ($0\le i_0<i_1<\dots<i_k<a+b+c$),
contributing to the product $m_k$, for some $k\ge 3$. Among the first three
corners of $u$, one is among the generators $x_i$, one is among the $y_i$,
and one is among the $z_i$. Therefore, $i_3=i_0+a+b+c$, which contradicts
the inequality $i_3<a+b+c$. Hence the moduli spaces of pseudo-holomorphic
curves involved in the definition of $m_k$ are all empty for $k\ge 3$, which
implies that $m_k=0$.

Lemma \ref{l:m2} also follows immediately at this point: in the case of a
pseudo-holomorphic map $u$ from a disc with $3$ marked points, the three
corners $p,q,r$ are all of different types (by Lemma \ref{l:3corners}),
and the first two corners $p,q$ correspond to generators among
$\{x_i,y_i,z_i\}$ while the last one $r$ corresponds to a generator
among $\{\bar{x}_i,\bar{y}_i,\bar{z}_i\}$. Therefore, $p$ and $q$ completely
determine $r$, and moreover it is easy to check from the above discussion
and from Figures~\ref{fig:vcs} and~\ref{fig:vcinters} that the image of the pseudo-holomorphic
map $u$ is also uniquely determined by the pair $(p,q)$. For example, if $p$
is of type $x$ and $q$ is of type $y$, then necessarily there exists $i<c$
such that $p=x_i$, $q=y_{i+a}$, and $r=\bar{z}_i$; moreover, it is easy to
check (see Lemma \ref{l:vcs} and Figure \ref{fig:vcs}) that the moduli
space determining the coefficient of $\bar{z}_i$ in $m_2(x_i,y_{i+a})$
consists of a single curve, regular, whose image $T_{xy,i}$ is the
triangular region
of $\Sigma_0$ delimited by arcs joining $p,q,r$ in the vanishing cycles
$L_i$, $L_{i+a}$, $L_{i+a+b}$. Therefore, we have
$m_2(x_i,y_{i+a})=\alpha_{xy,i}\,\bar{z}_i$, where $\alpha_{xy,i}=\pm
\exp(-\mathrm{Area}(T_{xy,i}))$. The situation is the same in all other cases.
\medskip

{\bf Remark.} The $a+b+c$
triangles $T_{xy,i}$ ($i<c$), $T_{yz,i}$ ($i<a$), $T_{zx,i}$ ($i<b$) are
all related to each other via the action of the cyclic group $\Z/(a+b+c)$.
Indeed, the diagonal multiplication by a power of $\zeta=\exp(\frac{2\pi
i}{a+b+c})$ induces a permutation of the vanishing cycles and of the
intersection points, preserving the cyclic ordering of the $L_i$ and the
types of their intersection points, and hence mapping every triangle in
$\Sigma_0$ with boundary in $\bigcup L_i$ to another such triangle.
A similar description holds for the triangles $T_{yx,i}$, $T_{zy,i}$, $T_{xz,i}$.

\subsection{Maslov index and grading}\label{ss:grading}
The aim of this section is to define a $\Z$-grading on the Floer complexes
$CF^*(L_i,L_j)$, and to compute the degree of the various generators.
Using the triviality of the canonical bundles of $\Sigma_0$ and $X$, it is
easy to prove (by considering the Lefschetz thimbles) that the Maslov class
of $L_i$ is trivial, and hence that it is possible to lift each vanishing cycle
to a {\it graded Lagrangian} submanifold of $\Sigma_0$, that we denote again
by $L_i$. This lets us associate a degree to each generator of the Floer
complex.

\begin{lemma}\label{l:grading}
There exists a natural choice of gradings, for which
$\deg(x_i)=\deg(y_i)=\deg(z_i)=1$ and
$\deg(\bar{x}_i)=\deg(\bar{y}_i)=\deg(\bar{z}_i)=2$.
\end{lemma}

Assume for simplicity that the symplectic form $\omega$ is compatible with
the standard complex structure of $\Sigma_0$ inherited from that of
$(\C^*)^3$, which allows us to define explicitly a holomorphic volume form
$\Omega$ on $\Sigma_0$ (i.e., a non-vanishing holomorphic $1$-form). Then,
given an oriented Lagrangian submanifold $L\subset \Sigma_0$, the {\it phase}
of $L$ is the function $\phi_L:L\to \R/2\pi\Z$ whose value at every point is
the argument of the (non-zero) complex number obtained by evaluating $\Omega$
on an oriented volume element in $L$ (in the 1-dimensional case,
$\phi_L(x)=\arg \Omega(v)$ for $v$ a tangent vector to $L$ at $x$ defining
the orientation of $L$). The Maslov class is the 1-cocycle representing
the obstruction to lift $\phi_L$
to a real-valued function; if it vanishes, then $L$ can be lifted to a {\it
graded Lagrangian} submanifold, i.e.\ we can choose a real-valued lift
of the phase, $\tilde{\phi}_L:L\to \R$. In the 1-dimensional case, the
relationship between Maslov index and phase is very simple: given a
transverse intersection point $p$ between two graded Lagrangians
$L,L'\subset\Sigma_0$, the Maslov index of $p\in CF^*(L,L')$ is equal to
the smallest integer greater than $\frac{1}{\pi}(\phi_{L'}(p)-\phi_L(p))$.

The holomorphic volume form $\Omega$ on $\Sigma_0$ can be defined from the
standard holomorphic volume form $\Omega_0=d\log x\wedge d\log y\wedge d\log
z$ on $(\C^*)^3$ by taking residues first along the hypersurface $X$ of
equation $x^ay^bz^c=1$ and then along the level set $W=0$. We can
characterize $\Omega$ as follows: $\Omega$ is the restriction to $\Sigma_0$
of a 1-form (that we denote again by $\Omega$) such that $\Omega\wedge
dw\wedge d(x^ay^bz^c)=\Omega_0$, i.e.\ (using the fact that $x^ay^bz^c=1$
along $X$)
$$\Omega\wedge (dx+dy+dz)\wedge
(\frac{a}{x}\,dx+\frac{b}{y}\,dy+\frac{c}{z}\,dz)=\frac{dx\wedge dy\wedge
dz}{xyz}.$$
(In fact the 1-form $\Omega$ determined in this way may differ from the
``usual'' one by a real positive factor, irrelevant for our purposes).
At this point it is easy to see why the Maslov class of $L_i$ is trivial:
indeed, $\Omega\wedge dw$ extends to a non-vanishing $(2,0)$-form on $X$,
whose phase over the Lefschetz thimble $D_i$ admits a real lift; because
$W$ maps $D_i$ to an embedded arc, the phase of $\Omega\wedge dw$ over the
boundary of $D_i$ and the phase of $\Omega$ over $L_i$ differ by a constant
term, so that the latter also admits a real lift.

At every point of $\Sigma_0$ except for the branch points of $\pi_x$, the
1-form $\Omega$ can be expressed as $\Theta\,dx$, for some meromorphic
function $\Theta$ over $\Sigma_0$ (with simple poles at the branch points
of $\pi_x$). The above equation becomes: $\Theta(\frac{c}{z}-\frac{b}{y})=
\frac{1}{xyz}$, which determines $\Theta$. At this point, the most direct
method of determination of the phases of the vanishing cycles $L_i$ at their
intersection points (and hence of the corresponding Maslov indices)
involves computer calculations; however we will attempt to give a sketch
of a geometric argument.

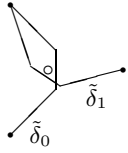
\begin{figure}[t]
\centering
\setlength{\unitlength}{1cm}
\begin{picture}(4,2)(-2,-1)
\put(1,0){\circle*{0.08}}
\put(-0.5,0.866){\circle*{0.08}}
\put(-0.5,-0.866){\circle*{0.08}}
\put(-0.5,0.866){\line(1,-1){0.6}}
\put(-0.5,-0.866){\line(1,1){0.6}}
\put(0.1,-0.266){\line(0,1){0.532}}
\put(1,0){\line(-4,-1){0.84}}
\put(0.16,-0.21){\line(-3,2){0.39}}
\put(-0.5,0.866){\line(1,-3){0.27}}
\put(-0.1,-0.85){\makebox(0,0)[cc]{\small $\tilde{\delta}_0$}}
\put(0.65,-0.35){\makebox(0,0)[cc]{\small $\tilde{\delta}_1$}}
\put(0,0){\circle{0.1}}
\end{picture}
\caption{The deformed cycles $\tilde{L}_j$ ($(a,b,c)=(1,1,1)$)}
\label{fig:vcsgrade}
\end{figure}

If we restrict ourselves to the
domain where $x$ is very small, then we have $y\simeq -z$, so that
$\Theta\simeq \frac{1}{(b+c)xy}$. Therefore, $\arg\Theta\simeq -\arg x-\arg
y$ in this region of $\Sigma_0$. Hence, the calculations are simplified if
we can deform the vanishing cycles $L_i$ in such a way that the intersection
points of a given type ($y$ or $z$) occur close to the origin in $\C^*$. Of
course this process preserves gradings and Maslov indices only
if the intersection pattern between the relevant vanishing cycles is not
affected by the deformation. We consider a deformation where $L_i$ is
replaced by a loop $\tilde{L}_i\subset\Sigma_0$, obtained as a double lift
of a piecewise smooth arc $\tilde{\delta}_i\subset\C^*$ joining two branch
points of $\pi_x$ (a deformation of $\delta_i$ with fixed end points).
The arc $\tilde{\delta}_0$ consists of three line
segments, two joining the end points $p,\bar{p}\in\mathrm{crit}(\pi_x)$ to
two complex conjugate points $q,\bar{q}$ very close to the origin, and such
that $0<\mathrm{Re}\,q\ll \mathrm{Im}\,q\ll 1$. The other arcs
$\tilde{\delta}_i$ are obtained from $\tilde{\delta}_0$ by the action of
$\Z/(a+b+c)$ (see Figure \ref{fig:vcsgrade}).

Assuming that $b<a+c$, this deformation can be carried out for intersections
of type $y$ without affecting the intersection pattern between $L_i$ and
$L_{i+b}$ or $L_{i+a+c}$, and in such a way that the
intersection occurs in the central portion of $\tilde{\delta}_i$
(see Figure \ref{fig:vcsgrade}).
The same is true for intersections of type $z$
when $c<a+b$. If we choose $a\ge b\ge c$ then these two assumptions hold,
so we can use this method to determine the degrees of $y_i,z_i,\bar{y}_i,
\bar{z}_i$.

We start by considering the portion of $\tilde{L}_0$ lying above the central
segment in $\tilde{\delta}_0$ (joining $q$ to $\bar{q}$). Recall that, for
$x$ small, the $b+c$ sheets of the covering $\pi_x$ (i.e.\ the $b+c$ roots of
$x^a y^b (-x-y)^c=1$) can be approximated by the roots of
$y^{b+c}=(-1)^c x^{-a}$. Hence, the possible values for the argument of $y$
are $\arg y\simeq -\frac{a}{b+c}\,\arg x+\pi \frac{c}{b+c} \mod \frac{2\pi}{b+c}$.
It follows from Lemma \ref{l:sigma0} that the two sheets of $\pi_x$
containing $\tilde{L}_0$ are those where $\arg y\simeq -\frac{a}{b+c}\,\arg
x+\epsilon \pi\frac{c}{b+c}$, for $\epsilon=\pm 1$. Hence, we have $\arg\Theta\simeq
\frac{a-b-c}{b+c}\,\arg x-\epsilon\pi\frac{c}{b+c}$. We choose to orient
$\tilde{L}_0$ in such a way that its projection goes counterclockwise around
the origin in the sheet corresponding to $\epsilon=1$, and clockwise in the
sheet corresponding to $\epsilon=-1$. With this understood, since the
projection of oriented tangent vector to $\tilde{L}_0$ is positively
proportional to $\epsilon i$, we obtain the following formula for the phase
of the central portion of $\tilde{L}_0$, modulo $2\pi$:
\begin{equation}\label{eq:phase0}
\phi(\tilde{L}_0)\simeq\frac{a-b-c}{b+c}\,\arg x+\epsilon\Bigl(\frac{\pi}{2}-
\frac{\pi\,c}{b+c}\Bigr).
\end{equation}
We choose a lift of $\tilde{L}_0$ (and hence also $L_0$ via the isotopy
between them) as a graded Lagrangian by setting the (real-valued) phase
of $\tilde{L}_0$ to be given by (\ref{eq:phase0}), choosing the
determination of $\arg x$ with the smallest absolute value; checking that the
choices made in the two portions of $\tilde{L}_0$ corresponding to
$\epsilon=\pm 1$
are consistent with each other is a tedious task, best left to a computer
program.

The phase of $\tilde{L}_j=\zeta^{-j}\cdot\tilde{L}_0$ is easily deduced from
the above calculations for $\tilde{L}_0$. Indeed, the above formula for
$\Theta$ implies that the value of $\arg\Theta$ at
the point $\zeta^{-j}\cdot p$ differs from that at the point $p$ by
$4\pi\frac{j}{a+b+c}$. On the other hand, the argument of the $x$ component
of the tangent vector to $\tilde{L}_j$ at $\zeta^{-j}\cdot p$ differs from
that of the tangent vector to $\tilde{L}_0$ at $p$ by
$-2\pi\frac{j}{a+b+c}$. Therefore, (\ref{eq:phase0}) implies that
$$\phi(\tilde{L}_j)\simeq \frac{a-b-c}{b+c}\,\Bigl(\arg x+\frac{2\pi j}{a+b+c}
\Bigr)+\epsilon\Bigl(\frac{\pi}{2}-\frac{\pi\,c}{b+c}\Bigr)+\frac{2\pi
j}{a+b+c},$$
or equivalently
\begin{equation}\label{eq:phases}
\phi(\tilde{L}_j)\simeq \frac{a-b-c}{b+c}\,\arg x
+\epsilon\Bigl(\frac{\pi}{2}-\frac{\pi\,c}{b+c}\Bigr)
+\frac{2\pi j\,a}{(a+b+c)(b+c)}.
\end{equation}
This formula can also be obtained directly by observing that the two
sheets of $\pi_x$ containing $\tilde{L}_j$ are those where $\arg y\simeq
-\frac{a}{b+c}\arg x-2\pi\frac{j}{b+c}+\epsilon\pi\frac{c}{b+c}$, for
$\epsilon=\pm 1$, by Lemmas~\ref{l:sigma0} and~\ref{l:vcs}.
As in the case of $\tilde{L}_0$, we choose a lift of $\tilde{L}_j$ whose
(real-valued) phase is given by (\ref{eq:phases}), using the determination
of $\arg x$ closest to $-2\pi\frac{j}{a+b+c}$.

We are now in a position to compare the phases of two vanishing cycles at
one of their intersection points. Consider an intersection point between
$\tilde{L}_i$ and $\tilde{L}_{i+b}$, corresponding to the intersection $y_i$
between $L_i$ and $L_{i+b}$. Comparing the values of $\arg y$ on both
vanishing cycles, it is easy to see that the intersection occurs in the
$\epsilon=1$ part of $L_i$ and in the $\epsilon=-1$ part of $L_{i+b}$.
Therefore, (\ref{eq:phases}) yields that, at the intersection point,
$$\phi(\tilde{L}_{i+b})-\phi(\tilde{L}_i)\simeq -2\Bigl(\frac{\pi}{2}-
\frac{\pi\,c}{b+c}\Bigr)+\frac{2\pi\,b\,a}{(a+b+c)(b+c)}
=\pi-\frac{2\pi\,b}{a+b+c},$$
which is between $0$ and $\pi$ since we have assumed that $b<a+c$.
Therefore, we have $\deg y_i=1$. Similarly, the intersection between
$\tilde{L}_i$ and $\tilde{L}_{i+c}$ corresponding to $z_i$ occurs in
the $\epsilon=-1$ part of $\tilde{L}_i$ and the $\epsilon=1$ part of
$\tilde{L}_{i+c}$, so that (\ref{eq:phases}) yields
$$\phi(\tilde{L}_{i+c})-\phi(\tilde{L}_i)\simeq 2\Bigl(\frac{\pi}{2}-
\frac{\pi\,c}{b+c}\Bigr)+\frac{2\pi\,c\,a}{(a+b+c)(b+c)}=
\pi-\frac{2\pi\,c}{a+b+c},$$
which is also between $0$ and $\pi$ since $c<a+b$. Therefore, $\deg z_i=1$.
In the case of $\bar{y}_i$, things are similar, but with one new subtlety:
in accordance with the above prescriptions, the determinations of $\arg x$
at the intersection point to be used for $\tilde{L}_i$ and
$\tilde{L}_{i+a+c}$ differ by $2\pi$. Therefore, from (\ref{eq:phases})
we now get (taking $\epsilon=-1$ for $\tilde{L}_i$ and $+1$ for
$\tilde{L}_{i+a+c}$)
$$\phi(\tilde{L}_{i+a+c})-\phi(\tilde{L}_i)\simeq -2\pi\,\frac{a-b-c}{b+c}
+2\Bigl(\frac{\pi}{2}-
\frac{\pi\,c}{b+c}\Bigr)+\frac{2\pi\,(a+c)\,a}{(a+b+c)(b+c)}=
\pi+\frac{2\pi\,b}{a+b+c},$$
which is between $\pi$ and $2\pi$; therefore, $\deg \bar{y}_i=2$.
Similarly, for $\bar{z}_i$ one finds that
$$\phi(\tilde{L}_{i+a+b})-\phi(\tilde{L}_i)\simeq -2\pi\,\frac{a-b-c}{b+c}
-2\Bigl(\frac{\pi}{2}-
\frac{\pi\,c}{b+c}\Bigr)+\frac{2\pi\,(a+b)\,a}{(a+b+c)(b+c)}=
\pi+\frac{2\pi\,c}{a+b+c},$$
which is also between $\pi$ and $2\pi$, so that $\deg \bar{z}_i=2$.

Finally, the degrees of $x_i$ and $\bar{x}_i$ can be deduced from those of
the intersections of types $y$ and $z$ by considering e.g.\ the triangles
$T_{xy,i}$, which gives that $\deg x_i+\deg y_{i+a}=\deg \bar{z}_i$, and
hence $\deg x_i=1$, and $T_{yz,i}$, which gives that $\deg y_i+\deg z_{i+b}=
\deg\bar{x}_i$, and hence $\deg\bar{x}_i=2$. This completes the proof of
Lemma \ref{l:grading}.

\subsection{The exterior algebra structure}\label{ss:antisym}

The aim of this section is to determine the coefficients appearing in
Lemma~\ref{l:m2}, by studying the orientations of the moduli spaces of
pseudo-holomorphic curves and the symplectic areas of their images
($T_{xy,i},\dots$).

\begin{lemma} \label{l:antisym}
If the symplectic form $\omega$ is anti-invariant under complex conjugation
and invariant under the action of \,$\Z/(a+b+c)$, then there exists
a constant
$\alpha\in\C^*$ such that $\alpha_{xy,i}=\alpha_{yz,i}=\alpha_{zx,i}=\alpha$ and
$\alpha_{yx,i}=\alpha_{zy,i}=\alpha_{xz,i}=-\alpha$ for all $i$.
Therefore, $m_2(x_i,y_{i+a})=-m_2(y_i,x_{i+b})$,
$m_2(y_i,z_{i+b})=-m_2(z_i,y_{i+c})$, and
$m_2(z_i,x_{i+c})=-m_2(x_i,z_{i+a})$.
\end{lemma}

The coefficients $\alpha_{xy,i},\dots$ are determined up to sign by the
symplectic areas of the triangular regions $T_{xy,i},\dots$ inside $\Sigma_0$.
To simplify notations, define
$$T_i=\begin{cases}
T_{xy,i} & \mathrm{if}\ 0\le i<c,\\
T_{zx,i-c} & \mathrm{if}\ c\le i<b+c,\\
T_{yz,i-b-c} & \mathrm{if}\ b+c\le i<a+b+c,
\end{cases}
\ \mathrm{and}\
T'_i=\begin{cases}
T_{xz,i} & \mathrm{if}\ 0\le i<b,\\
T_{yx,i-b} & \mathrm{if}\ b\le i<b+c,\\
T_{zy,i-b-c} & \mathrm{if}\ b+c\le i<a+b+c,
\end{cases}
$$
so that $T_i$ and $T'_i$ are the two triangles having either $x_i$ or
$\bar{x}_{i-b-c}$ as one of their vertices. We similarly define $\alpha_i$
and $\alpha'_i$ to be the coefficients associated to $T_i$ and $T'_i$ in
the formula giving $m_2$, namely $\alpha_i=\pm \exp(-\mathrm{Area}(T_i))$ and
$\alpha'_i=\pm \exp(-\mathrm{Area}(T'_i))$.
Then, as observed at the end of \S \ref{ss:products}, the invariance
properties of $\omega$ imply that the $a+b+c$ triangles $T_i$ form a
single orbit under the action of $\Z/(a+b+c)$, with $\zeta^{-q}\cdot T_i=
T_{i+q}$, and similarly for the other triangles $T'_i$, with $\zeta^{-q}\cdot
T'_i=T'_{i+q}$. Moreover, complex conjugation exchanges these two families
of triangular regions, by mapping $T_i$ to $T'_{b+c-i}$ (see Figure
\ref{fig:vcs}). It follows that all
of these triangles have the same symplectic area, and therefore that
the various constants $\alpha_i$ and $\alpha'_i$
are all equal up to sign.

In order to identify the signs, one needs to orient the relevant moduli
spaces of pseudo-holomorphic discs in some consistent way, which requires
the choice of a spin structure over each Lagrangian $L_i$. As explained
at the end of \S \ref{ss:fs}, we need to endow each $L_i$ with the spin
structure which extends to the corresponding thimble, i.e.\ the non-trivial
one.

We now describe a convenient rule for determining the correct
signs in the one-dimensional case, due to
Seidel \cite{Se3}. We start with the case of trivial spin structures.
Then to each intersection
point $p\in L_i\cap L_j$ ($i<j$) one can associate an orientation line
$\mathcal{O}_p$. This orientation line is canonically trivial when
$\deg p$ is even, whereas in the odd degree case, a choice of
trivialization of $\mathcal{O}_p$ is equivalent to a choice of orientation
of the line $T_p L_j$. If one considers a pseudo-holomorphic map
$u:D^2\to\Sigma_0$ contributing to $m_k$, whose image is a polygonal
region with $k+1$ vertices $p_0,\dots,p_k$, then the corresponding sign
factor is
actually an element of the tensor product $\Lambda=\mathcal{O}_{p_0}
\otimes\dots\otimes \mathcal{O}_{p_k}$. We can define a preferred
trivialization of $\Lambda$ by choosing, at each vertex of odd degree,
the orientation of the vanishing cycle which agrees with the positive
orientation on the boundary of the image of $u$.
The sign factor associated to $u$ is then equal to $+1$ with respect to this
trivialization of $\Lambda$ (or $-1$ with respect to the
other trivialization).
In the presence of non-trivial spin structures, this rule needs to be
modified as follows: fix a marked point on each $L_i$ carrying a non-trivial
spin structure (distinct from its
intersection points with the other vanishing cycles); then the sign
associated to $u$ is affected by a factor of $-1$ for each marked point
that the boundary of $u$ passes through \cite{Se3}.

It is worth mentioning that, while it is clear from the above construction
that the individual sign factors fail to be canonical and depend on some
choices, the various possibilities yield equivalent categories, since
the coefficients of Floer homology and Floer products simply differ by the
conjugation action of some diagonal matrix with $\pm 1$ coefficients.

In our case, we choose trivializations of the orientation lines as follows:
for every intersection point $p\in L_i\cap L_j$ of degree $1$ (i.e.,
one of $x_i,y_i,z_i$), we orient $T_p L_j$ consistently with the
boundary orientation of the single triangular region among
$T_0,\dots,T_{a+b+c-1}$ having $p$ among its vertices.
If we consider trivial spin
structures, then with this convention the sign factor associated to each triangle
$T_i$ is by definition equal to $+1$. In the case of $T'_i$, at
each of the two vertices of degree $1$ the chosen trivialization of $T_pL_j$
disagrees with the boundary orientation of the triangular
region, so that for trivial spin structures we get a sign factor of
$(-1)^2=+1$ again. Since we need to consider non-trivial spin structures,
we must introduce a marked point on each $L_i$; for example, we choose
this marked point in the portion of $L_i$ that corresponds to the top-most
edge on Figure \ref{fig:vcinters}. With this choice, the boundary of each
$T'_i$ passes through exactly one marked point (between the vertex of type
$z$ and that of type $y$), while the boundary of $T_i$ does not meet any
marked point. Therefore, with these conventions, the sign factors are
$+1$ for all $T_i$ and $-1$ for all $T'_i$; this completes the proof of
Lemma \ref{l:antisym}.

\subsection{Non-exact symplectic forms and non-commutative deformations}
\label{ss:deform}

The purpose of this section is to describe the effect on the category of
Lagrangian vanishing cycles
of $W$ of relaxing the assumptions made above on the symplectic form, losing
in particular its exactness. In order to make the vanishing cycle
construction well-defined, we will keep assuming that $\omega$ induces a
complete K\"ahler metric on $X$ and that the gradient of $W$ with respect
to this metric is bounded from below outside of a compact set. For example,
choosing a $3\times 3$ positive definite Hermitian matrix $(a_{ij})$, we
can endow $X$ with the symplectic form $$\omega=i\sum_{i,j=1}^3 a_{ij}
\frac{dz_i}{z_i}\wedge \frac{d\bar{z}_j}{\bar{z}_j}.$$ Observe that
$H_2(X,\Z)\simeq \Z$ is generated by the torus $T=\{(x,y,z)\in X,
\ |x|=|y|=|z|=1\}$ (for simplicity we assume $gcd(a,b,c)=1$). An easy
calculation shows that  \begin{equation}\label{eq:cohomclass}
[\omega]\cdot[T]=4\pi^2i\,(a\,(a_{23}-a_{32})+b\,
(a_{31}-a_{13})+c\,(a_{12}-a_{21})).\end{equation}
Many other choices of symplectic form are equally acceptable, and it
is important to mention that the most sensible course of action in presence
of a non-explicit symplectic form is to search for a {\it topological}
interpretation of the category of Lagrangian vanishing cycles, involving only topological
quantities such as the cohomology class of $\omega$.

In comparison to the restrictive situation considered above, the vanishing
cycles $L_j$ remain in the same smooth isotopy classes, because one can
continuously deform from one symplectic structure to the other. Hence,
the vanishing cycles are smoothly isotopic to the loops $L'_j\subset\Sigma_0$
introduced in \S \ref{ss:vcs}, but not necessarily Hamiltonian isotopic to
them. Nonetheless, because the ends of the non-compact Riemann surface
$\Sigma_0$ all have infinite volume, we can easily deform $L'_j$ into
loops $L''_j\subset \Sigma_0$ that are Hamiltonian isotopic to the vanishing
cycles, without modifying the pattern of the intersections between them.
More precisely, recall from \S \ref{ss:vcs} that each $L'_j$ is the double
lift via $\pi_x:\Sigma_0\to\C^*$ of an arc joining two branch points of
$\pi_x$. Then, by ``pulling'' a suitable portion of one of the two lifts
towards an end of $\Sigma_0$ (either towards infinity or towards zero in the
$x$-axis projection), we can make $L'_j$ sweep through an arbitrarily large
amount of symplectic area to obtain the desired $L''_j$, without affecting
the intersection points with the other vanishing cycles.

Since the vanishing cycles are Hamiltonian isotopic to the loops $L''_j$,
we may use $L''_j$ instead of the actual vanishing cycles in order to
determine the category $D(\FS(W))$. Hence, the symplectic
deformation does not affect in any way the generators of the Floer complexes
and the types of pseudo-holomorphic maps to be considered. The only
significant change has to do with the coefficients assigned to the various
pseudo-holomorphic discs appearing in the definition of $m_2$, as the
symplectic areas of the various triangular regions $T_i$ and $T'_i$
($i=0,\dots,a+b+c-1$) inside $\Sigma_0$ may now take more or less
arbitrary values instead of all being equal to each other. Because the
description of $\omega$ and of the vanishing cycles is not explicit, it is
hopeless (and useless) to calculate the individual coefficients $\alpha_i$
and $\alpha'_i$. However, we can state the following result:

\begin{lemma}\label{l:deform}
Lemmas \ref{l:isects}--\ref{l:grading} remain valid in the
more general case of an arbitrary symplectic form inducing a complete
K\"ahler metric on $X$ for which $|\nabla W|$ is bounded from below at
infinity. Moreover, the structure constants
for the composition $m_2$ are related by the identity
$$\frac{\prod_{i=0}^{a+b+c-1}\alpha_i}{\prod_{i=0}^{a+b+c-1}\alpha'_i}=
\frac{\prod\alpha_{xy,i}\,\prod\alpha_{yz,i}\,\prod\alpha_{zx,i}}
{\prod\alpha_{yx,i}\,\prod\alpha_{zy,i}\,\prod\alpha_{xz,i}}=
(-1)^{a+b+c}\exp(-[\omega]\cdot[T]).$$
\end{lemma}
\smallskip

The assumption of completeness of the induced K\"ahler metric can be
dropped if we have some other way of ensuring that the vanishing cycles are
well-defined and that the deformation from $L'_j$ to $L''_j$ can be carried
out without introducing new intersection points. In fact, the invariance of
Floer homology under Hamiltonian isotopies essentially implies that the
introduction of new intersection points in the deformation does not have
any particular impact on the derived category, so the only thing that
matters is actually the well-definedness of the vanishing cycles.

Although Lemma \ref{l:deform} seems to give only very partial
information about the constants $\alpha_i$ and $\alpha'_i$, it actually
completely determines the category $D(\FS(W))$. Indeed, simply by
rescaling the generators of the Floer complexes we can modify the
coefficients $\alpha_i$ and $\alpha'_i$ almost at will: for example,
replacing $x_i$ with $\lambda\,x_i$ has the effect of simultaneously
multiplying $\alpha_i$ and $\alpha'_i$ by $\lambda^{-1}$; similarly,
rescaling the generator $y_i$ simultaneously affects $\alpha_{i-a}$ (or
$\alpha_{i+b+c}$) and $\alpha'_{i+b}$. Still assuming $gcd(a,b,c)=1$,
it is not hard to check that the
only quantity left invariant by all rescalings of the generators is the
ratio $\prod \alpha_i/\prod \alpha'_i$, which is therefore sufficient to
characterize the derived category. This observation that the symplectic
deformations of $D(\FS(W))$ are governed by a single parameter is
naturally related to the fact that
the second Betti number of $X$ is equal to 1.

\proof[Proof of Lemma \ref{l:deform}]
The key observation to be made here is that the boundary of the 2-chain
$C=\sum T_i-\sum T'_i\subset\Sigma_0$ is exactly $\partial C=-\sum L_i$
(for a suitable choice of orientation of the $L_i$).
Indeed, looking at Figure \ref{fig:vcinters}, each of the six portions
of $L_i$ arises exactly once as an edge of one of the triangular
regions, and the boundary orientation of the triangular region is the
``clockwise'' orientation of $L_i$ in the case of $T_0,\dots,T_{a+b+c-1}$,
and the ``counterclockwise'' orientation
in the case of $T'_0,\dots,T'_{a+b+c-1}$. Recalling
that each vanishing cycle $L_i$ bounds a Lefschetz thimble $D_i$ in $X$,
we can build a 2-cycle $\tilde{C}\subset X$ by capping $C$ with
these $a+b+c$ Lagrangian discs. Next, observe that the sign factors
arising from the orientations of the moduli spaces remain the same as in
\S \ref{ss:antisym}, and that $\int_{D_i}\omega=0$, so that
$$\frac{\prod \alpha_i}{\prod\alpha'_i}=(-1)^{a+b+c}\,\frac{\prod \exp(
-\int_{T_i}\omega)}{\prod
\exp(-\int_{T'_i}\omega)}=(-1)^{a+b+c}\exp(-\int_C\omega)=
(-1)^{a+b+c}\exp(-[\omega]\cdot[\tilde{C}]).$$

Hence, the last step in the proof is to show that $[\tilde{C}]$ and $[T]$
are the same elements of $H_2(X,\Z)\simeq \Z$. A simple way to achieve this
is to compute the intersection pairing of $\tilde{C}$ with the relative cycle
$R=\{(x,y,z)\in X,\ x,y,z\in\R^+\}$, which intersects $T$ transversely
once at the point $(1,1,1)$.

To understand how $R$ intersects $\tilde{C}$, we compare the values of $W$
over $R$ and over $\tilde{C}$. By construction, $\tilde{C}$ is the union of
the 2-chain $C\subset\Sigma_0$, over which $W$ vanishes identically, and the
various Lefschetz thimbles $D_j$, which $W$ maps to straight line segments
joining the origin to the critical values $\lambda_j$. On the other hand,
the restriction to $R$ of $W=x+y+z$ is a proper function which takes real
positive values. With respect to the standard complex structure, $R$
is totally real and $W$ is holomorphic, so any critical point of $W_{|R}$ is
also a critical point of $W$, and in particular the minimum of $W$ over $R$
is a critical value of $W$. Indeed, a simple computation shows that the
minimum of $W$ over $R$ is exactly $(a+b+c)(a^a b^b c^c)^{-1/(a+b+c)}=
\lambda_0$, achieved at the critical point $p_0$ of $W$ corresponding to the
critical value $\lambda_0$.

It follows that the only point where $\tilde{C}$ and $R$ intersect is $p_0$.
Moreover, by considering the local model near $p_0$, it is easy to check
that this intersection is transverse, since the Hessian of $W$ at $p_0$
restricts to the tangent space $T_{p_0}D_0$ as a negative definite real
quadratic form, and to $T_{p_0}R$ as a positive definite real quadratic
form. Therefore the intersection number between $\tilde{C}$ and $R$ is equal
to $1$ (for a suitable choice of orientation that we will not discuss here),
and it follows that $[\tilde{C}]=[T]$ in $H_2(X,\Z)$.
\endproof

\subsection{B-fields and complexified deformations}\label{ss:deform2}

So far we have identified a real one-parameter family of deformations of
the category of Lagrangian vanishing cycles of $W$. To extend this to a complex family of
deformations, we need to introduce a non-trivial B-field, i.e.\ a closed
2-form $B\in\Omega^2(X,\R)$. The presence of a B-field affects Fukaya
categories by modifying the nature of the objects to be considered: namely,
one should consider pairs consisting of a Lagrangian submanifold and a
vector bundle over it equipped with a {\it projectively flat} (rather than
flat) connection with curvature equal to $-iB\otimes\mathrm{Id}$ (depending
on conventions, a factor of $2\pi$ is sometimes added).

In our case, we are considering Lagrangian vanishing cycles $L_j\simeq S^1$
arising as boundaries of the Lefschetz thimbles $D_j$. Since $\dim L_j=1$,
over $L_j$ every bundle is trivial and every connection is flat; moreover,
we can safely restrict ourselves to the case of line bundles. However,
the presence of the B-field results in a nontrivial holonomy. By Stokes'
theorem, if a U(1)-connection $\nabla_j=d+i\alpha_j$ is the restriction
to $L_j$ of a U(1)-connection with curvature $-iB$ over $D_j$, then the
holonomy of $\nabla_j$ around $L_j$ is given by
$\mathrm{hol}_{\nabla_j}(L_j)=\exp(\int_{L_j}i\alpha_j)=\exp(\int_{D_j}
i\,d\alpha_j)=\exp(-i\int_{D_j} B).$
Since this property characterizes
the connection $\nabla_j$ uniquely up to gauge, we can drop the line
bundle and the connection from the notation when considering the objects
$(L_j,E_j,\nabla_j)$ of $\FS(W,\{\gamma_j\})$.

However, we do need to take the holonomy of $\nabla_j$
into account when computing the twisted Floer differential and
compositions $m_k$, since the weight attributed to a given
pseudo-holomorphic disc $u:(D^2,\partial D^2)\to (\Sigma_0,\bigcup L_j)$
is modified by a factor corresponding to the holonomy along its boundary,
and becomes
$\pm\,\mathrm{hol}(u(\partial D^2))\,\exp(i\int_{D^2} u^*(B+i\omega))$.
More precisely, for each intersection point $p\in L_i\cap L_j$ we need to
fix an isomorphism between the fibers $(E_i)_{|p}$ and $(E_j)_{|p}$;
then it becomes possible to define the holonomy along the closed loop
$u(\partial D^2)$ using the parallel transport induced by $\nabla_j$ from
one ``corner'' of $u$ to the next one, and the chosen isomorphism at each
corner.

In this context, we now have the following result characterizing
$D(\FS(W))$:

\begin{lemma}\label{l:deform2}
Lemmas \ref{l:isects}--\ref{l:grading} remain valid for
an arbitrary symplectic form inducing a complete K\"ahler metric on $X$
for which $|\nabla W|$ is bounded from below at
infinity, and an arbitrary B-field. Moreover, the structure
constants for the composition $m_2$ are related by the identity
$$\frac{\prod_{i=0}^{a+b+c-1}\alpha_i}{\prod_{i=0}^{a+b+c-1}\alpha'_i}=
\frac{\prod\alpha_{xy,i}\,\prod\alpha_{yz,i}\,\prod\alpha_{zx,i}}
{\prod\alpha_{yx,i}\,\prod\alpha_{zy,i}\,\prod\alpha_{xz,i}}=
(-1)^{a+b+c}\exp(i[B+i\omega]\cdot[T]).$$
\end{lemma}
\smallskip

\proof
We again consider the 2-chain
$C=\sum T_i-\sum T'_i\subset\Sigma_0$, with boundary $\partial C=-\sum
L_j$, and the 2-cycle $\tilde{C}\subset X$ obtained by capping $C$ with
the Lagrangian discs $D_j$. We now have:
\begin{multline*}
\frac{\prod \alpha_i}{\prod\alpha'_i}=\frac{(-1)^{a+b+c}}{
\prod \mathrm{hol}_{\nabla_j}(L_j)}\,\frac{\prod \exp(
i\int_{T_i}B+i\omega)}{\prod
\exp(i\int_{T'_i}B+i\omega)}
=\frac{(-1)^{a+b+c}}{\prod \exp(\int_{D_j}\, -iB)} \exp(i\int_C B+i\omega)\\
=(-1)^{a+b+c}\exp(i[B+i\omega]\cdot[\tilde{C}]).
\end{multline*}

This completes the proof since $[\tilde{C}]=[T]$.
\endproof

It is interesting to observe that this statement reinterpretes the quantity
$\prod \alpha_i/\prod \alpha'_i$ in purely topological terms, thus avoiding
the pitfall of having to compute the individual coefficients attached to
the various pseudo-holomorphic discs in $\Sigma_0$. This outcome is rather
unsurprising since, whereas the individual coefficients
$\alpha_i$ and $\alpha'_i$ are heavily dependent on a number of
arbitrary choices, the underlying derived category of Lagrangian vanishing
cycles is expected
to depend only on the meaningful parameters -- in our case, the
cohomology class $[B+i\omega]$.

We would like to suggest that this feature reflects a general {\it
principle}.
Namely, the various structure coefficients of the Floer differentials and
products involved in the definition of the category $\FS(W)$ depend on many
choices and have no precise meaning in general. However, different sets of
values of the structure coefficients may become equivalent after a suitable
rescaling of the generators of the Floer complexes or other similarly benign
operations. Hence, we can reduce to a much smaller set of parameters (certain
combinations of the individual Floer coefficients) that actually govern the
structure of the category. Then, we expect the following statement to hold
in much greater generality than the examples studied here:

\begin{property} The structure of the derived category of Lagrangian vanishing
cycles is governed by deformation parameters which are all of the form
$\exp(i[B+i\omega]\cdot [C_j])$ for suitable 2-cycles $C_j\subset X$.
\end{property}

This is of course ultimately related to the fact that Floer homology and
Floer products can be defined over {\it Novikov rings}, counting
pseudo-holomorphic discs with coefficients that reflect relative
homology classes rather than actual symplectic areas; the version with
complex coefficients that we used here is then recovered from the version
with Novikov ring coefficients by evaluation at the point $[B+i\omega]$.

\section{Hirzebruch surfaces}

We now consider the case of Hirzebruch surfaces $\mathbb{F}_n$, for which
the mirror Landau-Ginzburg model consists of $X=(\C^*)^2$ equipped with
a superpotential of the form $$W=x+y+\frac{a}{x}+\frac{b}{x^n y}$$ for
some non-zero constants $a,b$. For simplicity we will only consider the
case of an exact symplectic form. Since different values of the
constants $a,b$ lead to mutually isotopic exact symplectic Lefschetz
fibrations, the actual choices do not matter (we can e.g.\ assume
$a=b=1$ or any other convenient choice).

\subsection{The case of $\mathbb{F}_0$ and $\mathbb{F}_1$}\label{ss:f0f1}

The first two Hirzebruch surfaces
 $\mathbb{F}_0=\CP^1\times\CP^1$ and $\mathbb{F}_1$
(i.e., $\CP^2$ blown up at one point) need to be considered separately.

\begin{prop}\label{prop:f0}
When $n=0$, there exists a system of arcs $\{\gamma_i\}$ such that
$\FS(W,\{\gamma_i\})$ is equivalent to the full subcategory of
$\db{\coh(\mathbb{F}_0)}$ whose objects are $\mathcal{O},$ $\mathcal{O}(1,0),$
$\mathcal{O}(0,1),$ $\mathcal{O}(1,1).$
Therefore, $D(\FS(W))\simeq \db{\coh(\mathbb{F}_0)}$.
\end{prop}

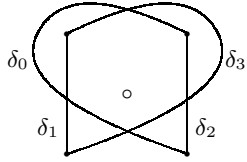
\begin{figure}[ht]
\centering
\setlength{\unitlength}{1cm}
\begin{picture}(3,2)(-1.5,-0.8)
\put(0.8,0.8){\circle*{0.08}}
\put(-0.8,0.8){\circle*{0.08}}
\put(0.8,-0.8){\circle*{0.08}}
\put(-0.8,-0.8){\circle*{0.08}}
\put(0,0){\circle{0.1}}
\put(-0.8,-0.8){\line(0,1){1.6}}
\put(0.8,-0.8){\line(0,1){1.6}}
\qbezier[250](0.8,0.8)(-0.8,1.6)(-1.2,0.8)
\qbezier[250](-1.2,0.8)(-1.6,0)(0.8,-0.8)
\qbezier[250](-0.8,0.8)(0.8,1.6)(1.2,0.8)
\qbezier[250](1.2,0.8)(1.6,0)(-0.8,-0.8)
\put(1.05,-0.45){\makebox(0,0)[cc]{\small $\delta_2$}}
\put(-1.05,-0.45){\makebox(0,0)[cc]{\small $\delta_1$}}
\put(-1.45,0.45){\makebox(0,0)[cc]{\small $\delta_0$}}
\put(1.45,0.45){\makebox(0,0)[cc]{\small $\delta_3$}}
\end{picture}
\caption{The vanishing cycles for $\mathbb{F}_0$}
\label{fig:vcf0}
\end{figure}

\proof
The four critical values of $W=x+y+\frac{a}{x}+\frac{b}{y}$ are
$\pm 2a^{1/2}\pm 2b^{1/2}$. Up to an exact deformation which does not
affect the category of Lagrangian vanishing cycles, we can choose $a>b>0$, and assume the
symplectic form to be anti-invariant under reflection about the
imaginary axis $(x,y)\mapsto (-\bar{x},-\bar{y})$. We choose
$\Sigma_0=W^{-1}(0)$ as our reference fiber, and join it to the singular
fibers by considering arcs $\gamma_i$ that pass below the real axis in
$\C$, so that the clockwise ordering of the critical values agrees with
their natural ordering
$-2a^{1/2}-2b^{1/2}<-2a^{1/2}+2b^{1/2}<2a^{1/2}-2b^{1/2}<2a^{1/2}+2b^{1/2}$.
The projection $\pi_x$ to the $x$ variable realizes $\Sigma_0$ as a double
cover of $\C^*$ branched at four points, and the vanishing cycles $L_i$
can be represented as double lifts of the arcs $\delta_i\subset\C^*$ shown
in Figure \ref{fig:vcf0}.

It follows that $\mathrm{Hom}(L_1,L_2)=0$, while
$\mathrm{Hom}(L_0,L_1)$, $\mathrm{Hom}(L_2,L_3)$, $\mathrm{Hom}(L_0,L_2)$,
and $\mathrm{Hom}(L_1,L_3)$ are two-dimensional; label the corresponding
intersection points $L_0\cap L_1=\{s,t\}$, $L_2\cap L_3=\{s',t'\}$,
$L_0\cap L_2=\{u,v\}$, $L_1\cap L_3=\{u',v'\}$. Finally,
$\mathrm{Hom}(L_0,L_3)$ has rank 4.
By considering the triangular regions delimited by the vanishing cycles in
$\Sigma_0$, and using the symmetry of the configuration with respect to
$(x,y)\mapsto (-\bar{x},-\bar{y})$, we can easily show that
$m_2(s,u')=m_2(s',u)$, $m_2(t,u')=m_2(t',u)$, $m_2(s,v')=m_2(s',v)$, and
$m_2(t,v')=m_2(t',v)$; these four elements of
$\mathrm{Hom}(L_0,L_3)$ are proportional to the generators. All other
products vanish ($m_k=0$ for $k\neq 2$). Finally, gradings can be chosen
so that all morphisms have degree $0$ (the verification is left to the
reader).

Therefore, the category $\FS(W,\{\gamma_i\})$ is indeed equivalent
to the full subcategory of $\db{\coh(\mathbb{F}_0)}$ whose objects are
$\mathcal{O},$ $\mathcal{O}(1,0),$
$\mathcal{O}(0,1),$ $\mathcal{O}(1,1),$
as can be seen by thinking of $(s,t)$ and
$(u,v)$ as homogeneous coordinates on the two factors of
$\mathbb{F}_0=\CP^1\times\CP^1$. Since these four line bundles form a
full strong exceptional collection generating $\db{\coh(\mathbb{F}_0)}$, the result follows.
\endproof

Alternatively, Proposition \ref{prop:f0} can also be obtained as a direct corollary of
a general product formula for categories of Lagrangian vanishing cycles of Lefschetz fibrations
of the form $(X_1\times X_2,W_1+W_2)$ (\cite{AKOS}, cf.\ also \S
\ref{ss:product}).

\begin{prop}\label{prop:f1}
When $n=1$, there exists a system of arcs $\{\gamma_i\}$ such that
$\FS(W,\{\gamma_i\})$ is equivalent to the full subcategory of
$\db{\coh(\mathbb{F}_1)}$ whose objects are $\mathcal{O},$
$\pi^*(T_{\mathbb{P}^2}(-1)),$ $\pi^*(\mathcal{O}_{\mathbb{P}^2}(1)),$
$\mathcal{O}_E$
(where $E$ is the exceptional curve and $\pi:\mathbb{F}_1\to\CP^2$ is
the blow-up map).
Therefore, $D(\FS(W))\simeq \db{\coh(\mathbb{F}_1)}$.
\end{prop}

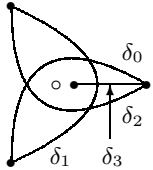
\begin{figure}[ht]
\centering
\setlength{\unitlength}{1.2cm}
\begin{picture}(4,2)(-2,-1)
\put(1,0){\circle*{0.08}}
\put(-0.5,0.866){\circle*{0.08}}
\put(-0.5,-0.866){\circle*{0.08}}
\put(0.2,0){\circle*{0.08}}
\qbezier[240](-0.5,0.866)(1.4,0)(-0.5,-0.866)
\qbezier[220](-0.5,0.866)(-0.5,-0.866)(1,0)
\qbezier[220](1,0)(-0.5,0.866)(-0.5,-0.866)
\put(0.2,0){\line(1,0){0.8}}
\put(0.85,0.35){\makebox(0,0)[cc]{\small $\delta_0$}}
\put(0.05,-0.8){\makebox(0,0)[cc]{\small $\delta_1$}}
\put(0.85,-0.35){\makebox(0,0)[cc]{\small $\delta_2$}}
\put(0.63,-0.8){\makebox(0,0)[cc]{\small $\delta_3$}}
\put(0.6,-0.6){\vector(0,1){0.6}}
\put(0,0){\circle{0.1}}
\end{picture}
\caption{The vanishing cycles for $\mathbb{F}_1$}
\label{fig:vcf1}
\end{figure}

\proof
We choose $a=b=1$, and equip $X$ with a symplectic form that is
anti-invariant under complex conjugation.
Let $(\lambda_i)_{0\le i\le 3}$ be the four critical values of
$W=x+y+\frac{1}{x}+\frac{1}{xy}$, ordered clockwise around the origin
so that $\mathrm{Im}(\lambda_0)>0$, $\lambda_1\in\R_+$,
$\mathrm{Im}(\lambda_2)<0$, and $\lambda_3\in\R_-$.
We choose $\Sigma_0=W^{-1}(0)$ as reference
fiber, and choose the arcs $\gamma_i$ joining $0$ to $\lambda_i$ to be
straight lines.
The projection $\pi_x$ to the $x$ variable realizes $\Sigma_0$ as a double
cover of $\C^*$ branched at four points, and the vanishing cycles $L_i$
can be represented as double lifts of the arcs $\delta_i\subset\C^*$ shown
in Figure \ref{fig:vcf1}.

The corresponding category of vanishing cycles can then be studied explicitly. In fact,
much of the work has already been carried out in \S \ref{sec:fsp2}, since
the situation for $L_0,L_1,L_2$ is rigorously identical (including grading
and orientation issues) to that previously considered for the three
vanishing cycles of the Lefschetz fibration mirror to $\CP^2$. While the
choice of grading used in \S \ref{sec:fsp2} yields morphisms in degrees 1
and 2, a different choice of gradings (shifting $L_1$ by 1 and $L_2$ by 2)
ensures that all morphisms between $L_0,L_1,L_2$ have degree 0.
This readily
implies that a category equivalent to the derived category of $\CP^2$ can
be realized inside $D(\FS(W))$ as a full subcategory, with the
exceptional collection $L_0,L_1,L_2$ corresponding to the exceptional
collection $\mathcal{O},T_{\mathbb{P}^2}(-1),\mathcal{O}(1)$ dual to the
standard one. (This claim can of course also be verified ``by hand''
following the same outline of argument as in \S \ref{sec:fsp2}).

 From Figure \ref{fig:vcf1} it is clear that $\mathrm{Hom}(L_0,L_3)$ and
$\mathrm{Hom}(L_2,L_3)$ are one-dimensional, (call their generators $p_0$
and $p_2$), while $\mathrm{Hom}(L_1,L_3)$ has rank 2 (call its generators
$q$ and $q'$). To be consistent with the notation of \S \ref{sec:fsp2}, call
$x_0,y_0,z_0$ (resp.\ $x_1,y_1,z_1$; resp.\ $\bar{x},\bar{y},\bar{z}$)
the generators of $\mathrm{Hom}(L_0,L_1)$ (resp.\ $\mathrm{Hom}(L_1,L_2)$;
resp.\ $\mathrm{Hom}(L_0,L_2)$). Then, looking at the various
pseudo-holomorphic discs in $\Sigma_0$ (including a constant one at the
triple intersection of $L_0,L_2,L_3$), we have:
$m_2(x_0,q)=m_2(x_0,q')=0,$ $m_2(y_0,q)=\alpha\,p_0,$ $m_2(y_0,q')=0,$
$m_2(z_0,q)=0,$ $m_2(z_0,q')=\alpha'\,p_0,$ $m_2(x_1,p_2)=0,$
$m_2(y_1,p_2)=-\alpha\,q',$ $m_2(z_1,p_2)=\alpha'\,q,$
$m_2(\bar{x},p_2)=p_0,$ $m_2(\bar{y},p_2)=m_2(\bar{z},p_2)=0$
(for some non-zero constants $\alpha,\alpha'$).
Moreover, for a suitable choice of grading of $L_3$ it can be checked that
all morphisms have degree 0.

It is then easy to check that these formulas correspond exactly to the
composition formulas in the full subcategory of $\db{\coh(\mathbb{F}_1)}$ whose
objects are the pull-backs $\mathcal{O}$,
$\pi^*(T_{\mathbb{P}^2}(-1))$, $\pi^*(\mathcal{O}_{\mathbb{P}^2}(1))$, and
the structure sheaf $\mathcal{O}_E$ of the exceptional curve
(If one follows the analogy suggested by the notation between the
morphisms from $L_0$ to $L_2$ and the homogeneous
coordinates on $\CP^2$, then the blow-up point is located at $(1:0:0)$). The result follows.
\endproof

\subsection{Other Hirzebruch surfaces}

For larger values of $n$, the situation becomes different:

\begin{lemma}
If $n\ge 2$, then the Lefschetz fibrations over $(\C^*)^2$ defined by
$W=x+y+\frac{1}{x}+\frac{1}{x^ny}$ and $\tilde{W}=x+y+\frac{1}{x^ny}$ are
isotopic. Therefore, $D(\FS(W))\simeq D(\FS(\tilde{W}))\simeq
\db{\coh(\CP^2(n,1,1))}$.
\end{lemma}

\proof
Consider the maps $W_a=x+y+\frac{a}{x}+\frac{1}{x^ny}$ for $a\in
[0,1]$. The key observation is that the $n+2$ critical points of $W_a$
remain distinct and stay in a compact subset of $(\C^*)^2$. Indeed, the
critical points of $W_a$ are the solutions of
$$\begin{cases} 1-\frac{a}{x^2}-\frac{n}{x^{n+1}y}=0\\
1-\frac{1}{x^ny^2}=0,\end{cases}$$
i.e.\ $$y=nx^{1-n}(x^2-a)^{-1},\ \mathrm{and}\ x^{n-2}(x^2-a)^2-n^2=0.$$
It is easy to check that for $|a|\le 1$ the roots of this
equation satisfy $1\le |x|\le \sqrt{n+1}$. It follows that $|x^2-a|=
n|x|^{1-\frac{n}{2}}$ is bounded between two positive constants, and hence
that $y=nx^{1-n}(x^2-a)^{-1}=(x^2-a)/nx$ is also bounded between two
positive constants independently of $a$. Hence the critical points of $W_a$
remain inside a compact subset of $(\C^*)^2$. Moreover, the polynomial
$P(x)=x^{n-2}(x^2-a)^2-n^2$ always has simple roots when $|a| \le 1$,
since the roots of $P'(x)=x^{n-3}(x^2-a)((n+2)x^2-(n-2)a)$ are $0$, $\pm
\sqrt{a}$, and $\pm \sqrt{\frac{n-2}{n+2}a}$, where $P$ never vanishes.
In fact, even though this is not necessary for the argument, the critical
values of $W_a$ also remain distinct throughout the deformation, since
at a critical point we have $W_a=\frac{n+2}{n}x+\frac{n-2}{n}\frac{a}{x}$,
which as a function of $x$ is injective over $\{|x|\ge 1\}$.

Therefore, $W_a$ defines an exact symplectic Lefschetz fibration on
$(\C^*)^2$ for all $a\in [0,1]$, which allows us to match the vanishing
cycles of $W_1=W$ with those of $W_0=\tilde{W}$. The resulting
categories of vanishing cycles differ at most by a deformation of the structure coefficients
of the compositions $m_2$, but since the isotopy is through {\it exact}
Lagrangian vanishing cycles, we need not worry about those (see also the
argument for Lemma \ref{l:deform}).

We can therefore conclude that $D(\FS(W))\simeq D(\FS(\tilde{W}))$.
Since $((\C^*)^2,\tilde{W})$ is exactly the mirror to $\CP^2(n,1,1)$
studied at length in \S \ref{sec:fsp2}, our result for weighted projective
planes implies that this category
is also equivalent to $\db{\coh(\CP^2(n,1,1))}$.
\endproof

For $n=2$, it is well-known that $\db{\coh(\mathbb{F}_2)}\simeq
\db{\coh(\CP^2(2,1,1))}$, so we get the expected result.
However, for $n\ge 3$ this is no longer true. Namely, the fully faithful
functor $MK_n$ constructed in \S \ref{ss:dbfn} allows us to view the
category $\db{\coh(\mathbb{F}_n)}$
as a full subcategory of $\db{\coh(\CP^2(n,1,1))},$ generated by the exceptional
collection $(\mathcal{O},\mathcal{O}(1),\mathcal{O}(n),\mathcal{O}(n+1))$.
It is
therefore a natural question to ask whether this subcategory can be singled
out on the mirror side, by selecting $4$ of the $n+2$ critical points of
$W$. It turns out that this is indeed the case. Our first result in this
direction is the following:

\begin{lemma}
For $n\ge 3$, in the limit $b\to 0$, $n-2$ of the critical values
of the superpotentials $W_b=x+y+\frac{1}{x}+\frac{b}{x^ny}$ go to infinity,
while the remaining four critical points stay in a bounded region.
\end{lemma}

\proof
The $x$ coordinates of the critical points of $W_b$ are the solutions of
$$x^{n-2}(x^2-1)^2-n^2b=0.$$
As $b\to 0$, four roots of this equation converge to $\pm 1$, while
the remaining $n-2$ converge to $0$. Since at a critical point we also have
$y=nbx^{1-n}(x^2-1)^{-1}=\frac{1}{n}(x-\frac{1}{x})$ and
$W_b=\frac{n+2}{n}x+\frac{n-2}{n}\frac{1}{x}$, we conclude that four
critical points of $W_b$ converge to $(\pm 1,0)$, with the corresponding
critical values converging to $\pm 2$, while the others escape to infinity.
\endproof

This suggests that the deformation $b\to 0$ singles out a subcategory of
$D(\FS(W_b))$, obtained by restricting oneself to the preimage of a
disc containing only four critical values of $W_b$. We start by describing
the case $n=3$.

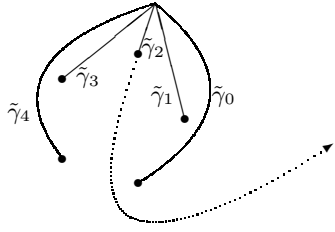
\begin{figure}[t]
\centering
\setlength{\unitlength}{9mm}
\begin{picture}(4,3.2)(-2,-1.5)
\put(1,0){\circle*{0.1}}
\put(0.31,0.95){\circle*{0.1}}
\put(0.31,-0.95){\circle*{0.1}}
\put(-0.81,0.59){\circle*{0.1}}
\put(-0.81,-0.59){\circle*{0.1}}
\qbezier[40](0.31,0.95)(-0.5,-1.3)(0.3,-1.5)
\qbezier[40](0.3,-1.5)(1.1,-1.7)(3,-0.5)
\put(3.2,-0.36){\vector(3,2){0}}
\qbezier[200](-0.81,-0.59)(-2,0.9)(0.55,1.7)
\put(-0.81,0.59){\line(5,4){1.36}}
\put(0.31,0.95){\line(1,3){0.25}}
\put(1,0){\line(-1,4){0.43}}
\qbezier[200](0.31,-0.95)(2.3,0.3)(0.55,1.7)
\put(-1.62,0){\small $\tilde{\gamma}_4$}
\put(-0.6,0.52){\small $\tilde{\gamma}_3$}
\put(0.36,0.95){\small $\tilde{\gamma}_2$}
\put(0.5,0.25){\small $\tilde{\gamma}_1$}
\put(1.4,0.25){\small $\tilde{\gamma}_0$}
\end{picture}
\caption{The deformation $b\to 0$ for $n=3$}
\label{fig:deformf3}
\end{figure}

For $n=3$, we can study explicitly the deformation process as
$b$ changes from $1$ to a value close to $0$. For $b=1$ the five critical
values of $W_b$ form a pentagon roughly centered at the origin (and can
for all practical purposes be identified with the critical values of the
superpotential mirror to $\CP^2(3,1,1)$). As $b$
decreases along the real axis, two things happen: first, the two complex
conjugate critical points with $\mathrm{Re}(W_b)>0$ merge and turn into
two real critical points; then, one of these two real critical points
escapes to infinity as $b\to 0$. The process is easier to visualize if one
avoids the two values of
$b$ in the interval $(0,1)$ for which two critical values of $W_b$ coincide,
by considering e.g.\ a deformation from
$b=1$ to $b=0$ where the imaginary part of $b$ is kept positive.
It is then easy to check that, as $b\to 0$, two critical values converge to $2$
and two others converge to $-2$, while the fifth one escapes to infinity
in the manner represented on Figure \ref{fig:deformf3}.

\begin{figure}[ht]
\centering
\setlength{\unitlength}{8mm}
\begin{picture}(6,2)(-3,0)
\put(-2.5,1.5){\circle*{0.1}}
\put(-1.5,1.5){\circle*{0.1}}
\put(-2,0){\line(-1,3){0.5}}
\put(-2,0){\line(1,3){0.5}}
\put(-2.8,0.8){\small $\gamma_i$}
\put(-1.65,0.8){\small $\gamma_{i+1}$}
\put(0,0.9){\vector(1,0){0.6}}
\put(0,0.9){\vector(-1,0){0.6}}
\put(3,1.5){\circle*{0.1}}
\put(2,1.5){\circle*{0.1}}
\put(2.5,0){\line(-1,3){0.5}}
\put(2.5,0){\line(-2,3){0.6}}
\qbezier[150](1.9,0.9)(1,2.5)(3,1.5)
\put(2.3,0.8){\small $\gamma_i$}
\put(0.75,0.8){\small $L\gamma_{i+1}$}
\end{picture}

\caption{The (left) sliding operation $(\gamma_i,\gamma_{i+1})\longleftrightarrow
(L\gamma_{i+1},\gamma_i)$}
\label{fig:sliding}
\end{figure}
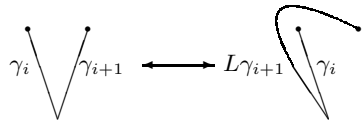

Therefore, if we consider the category of Lagrangian vanishing cycles
associated to the system of arcs $\tilde{\gamma}_0,\dots,\tilde{\gamma}_4$
represented on Figure \ref{fig:deformf3}, the deformation $b\to 0$ singles
out the full subcategory generated by the four vanishing cycles $\tilde{L}_0,
\tilde{L}_1,\tilde{L}_3,\tilde{L}_4$ (where $\tilde{L}_i$ is the vanishing
cycle associated to $\tilde{\gamma}_i$). The collection of arcs
$\{\tilde{\gamma}_i\}$ looks very different from the collection $\{\gamma_i\}$
considered in \S \ref{sec:fsp2}, but they are related to each other by
a sequence of elementary sliding transformations performed on consecutive
arcs (see Figure \ref{fig:sliding}).

It follows immediately from Definition \ref{def:fs} that every
ordered collection of arcs yields a full exceptional collection
generating $D(\FS(W))$;
it was shown by Seidel that (left or right) sliding operations on
collections of arcs correspond to (left or right) mutations of the
corresponding exceptional collections \cite{Se1}. With this is mind,
and identifying implicitly the critical points of $W_1$ with those of
the superpotential mirror to $\CP^2(3,1,1),$
it is easy to check that the left
dual to the exceptional collection $(\tilde{L}_0,\dots,\tilde{L}_4)$
associated to the arcs $\{\tilde{\gamma}_i\}$ is equivalent (up to some
shifts) to the exceptional collection associated to the arcs
$(\gamma_2,\gamma_3,\gamma_4,\gamma_0,\gamma_1)$. Moreover, using
$\Z/5$-equivariance for $\CP^2(3,1,1),$ there exists an auto-equivalence
of $D(\FS(W_1))$
which maps this exceptional collection to the one associated to the
collection of arcs $(\gamma_0,\dots,\gamma_4)$ considered in \S
\ref{sec:fsp2}.

Recall that the two exceptional collections for
$\db{\coh(\CP^2(3,1,1))}$ presented in \S \ref{sec:algebra} are mutually
dual (cf.\ Example \ref{ex:dual}), and that Theorem \ref{thm:fswp2}
identifies the exceptional collection associated to the arcs
$(\gamma_0,\dots,\gamma_4)$ with that given by Corollary \ref{cor:excoll}.
Therefore, there is an equivalence of categories which maps the exceptional
collection $(\tilde{L}_0,\dots,\tilde{L}_4)$ for $D(\FS(W_1))$ to the
exceptional collection $(\mathcal{O},\dots,\mathcal{O}(4))$ for
$\db{\coh(\CP^2(3,1,1))}$. The full subcategory of $D(\FS(W_1))$ singled out
by the deformation $b\to 0$ is that generated by the exceptional
collection $(\tilde{L}_0,\tilde{L}_1,\tilde{L}_3,\tilde{L}_4)$, which
corresponds under the above identification to
the full subcategory of $\db{\coh(\CP^2(3,1,1))}$ generated by the exceptional
collection
$(\mathcal{O},\mathcal{O}(1),\mathcal{O}(3),\mathcal{O}(4))$, which is in turn
known to be equivalent to the derived category of the Hirzebruch surface
$\mathbb{F}_3$ (see \S \ref{ss:dbfn}).
\medskip

A similar analysis of the deformation $b\to 0$ can be carried out for all
values of $n,$ and leads to the following result:

\begin{prop}\label{prop:fnsubcat}
Given any $n\ge 3$ and $R\gg 2$, and assuming that $b$ is sufficiently
close to $0$, the full subcategory of $D(\FS(W_b))$ arising from
restriction to the open domain $\{|W_b|<R\}$ is equivalent to
$\db{\coh(\mathbb{F}_n)}.$
\end{prop}

In order to prove this proposition we need a lemma about mutations in
the standard full exceptional collection
$\left( \O,\O(1),\dots,\O(n+1)\right)$ on the weighted
projective plane $\CP^2(n,1,1).$ Let us fix a pair
$\left(\O(k),\O(k+1)\right)$ with $2<k<n.$
Denote by $F_{k+2}$ the  mutation of the object $\O(k+2)$ to the left through
$\O(k),\O(k+1),$ i.e.\ $F_{k+2}\cong L^{(2)}\O(k+2).$
Performing the same mutations on $\O(k+3),\dots,\O(n+1)$ we obtain
exceptional objects $F_{i}=L^{(2)}\O(i)$ for $k+2\le i\le n+1$ and a new
exceptional
collection
$$
\left( \O,\dots,\O(k-1), F_{k+2},\dots,F_{n+1},\O(k),\O(k+1)\right).
$$
Denote by $G_k, G_{k+1}$ the left mutations of $\O(k), \O(k+1)$
respectively through all $F_i.$
We get an exceptional collection
$$
\left( \O,\dots,\O(k-1), G_k, G_{k+1}, F_{k+2},\dots,F_{n+1}\right).
$$
Denote by $\D$ the triangulated subcategory of the category
$\bD^b(\coh(\CP^2(n,1,1)))$ generated by the collection
$\left(\O,\O(1),G_k,G_{k+1}\right)$
\begin{lemma}\label{l:fncatmutation}
The triangulated subcategory $\D$ coincides with the subcategory
$$
\langle\O,\O(1),\O(n),\O(n+1)\rangle.
$$
\end{lemma}
\begin{proof} This Lemma is equivalent to the statement that
the subcategory $\langle G_k, G_{k+1}\rangle$ coincides with the
subcategory
$\langle \O(n), \O(n+1)\rangle.$
First, let us show that $\O(n)$ and $\O(n+1)$ belong to $\langle G_k,
G_{k+1}\rangle.$
Since $\Hom(\O(l), \O(s))=0$ for $l=n,n+1$ and $0\le s<k,$ we can
immediately conclude that $\O(n)$ and $\O(n+1)$ belong to $\langle
G_k,G_{k+1},F_{k+2},\dots,F_{n+1}\rangle.$
Therefore, it is sufficient to check that
$$
\Hom^{\bdot}(F_i, \O(n))=0, \qquad \Hom^{\bdot}(F_i, \O(n+1))=0
$$
for all $k+2\le i\le n+1.$

By definition of $F_i$ there are   distinguished triangles
\begin{align}
&T_i\lto V_i\otimes\O(k+1)\lto\O(i),\label{tr1}\\
&F_i\lto W_i\otimes \O(k)\lto T_i,\label{tr2}
\end{align}
with $V_i=\Hom(\O(k+1),\O(i))$ and
$W_i=\Hom(\O(k), T_i).$
It is clear that $V_i\cong S^{i-k-1}U,$
where $U$ is the two dimensional vector space $H^0(\CP^2(n,1,1),\O(1)).$
Considering the sequence of Hom's from $\O(k)$ to the
triangle (\ref{tr1}), it is easy to check that $W_i\cong S^{i-k-2}U$
(we use an isomorphism $\Lambda^2 U\cong \kk$).

We have isomorphisms
$$
\Hom(V_i\otimes \O(k+1), \O(n+1))= S^{i-k-1} U^*\otimes S^{n-k}U\cong
\bigoplus_{j=0}^{i-k-1} S^{n-i+1+2j}U,
$$
which implies that
$$
\Hom(T_i, \O(n+1))\cong\bigoplus_{j=1}^{i-k-1} S^{n-i+1+2j}U.
$$
On the other hand, there are isomorphisms
$$
\Hom(W_i\otimes \O(k), \O(n+1))= S^{i-k-2} U^*\otimes S^{n-k+1}U\cong
\bigoplus_{j=1}^{i-k-1} S^{n-i+1+2j}U,
$$
and, moreover, it can be checked that the natural
morphism $\Hom(T_i, \O(n+1))\to\Hom(W_i\otimes\O(k), \O(n+1))$
is an isomorphism.
Hence, $\Hom^{\bdot}(F_i, \O(n+1))=0$ for all $k+2\le i\le n+1.$
By the same reasons $\Hom^{\bdot}(F_i, \O(n))=0$ for all $k+2\le i\le
n+1.$
Thus the subcategory $\langle \O(n), \O(n+1)\rangle$ is contained in
$\langle G_k, G_{k+1}\rangle.$

Since $\Hom(G_k,G_{k+1})\cong U\cong \Hom(\O(n),\O(n+1))$,
these two categories are both equivalent to the derived  category
of representations of the quiver with two vertices and two arrows
$\bullet\rightrightarrows\bullet,$ and, as consequence, it can be easily
shown that
they are equivalent.
\end{proof}

\proof[Proof of Proposition \ref{prop:fnsubcat}]
The argument is similar to the case $n=3$: in the initial configuration, for
$b=1,$ the $n+2$ critical values of $W_b$ approximate a regular polygon,
and can essentially be identified with the critical values of the
superpotential mirror to $\CP^2(n,1,1).$ We label these critical values
by integers from $0$ to $n+1,$ with $0$ corresponding to the positive
real critical value, and continuing counterclockwise. As the value of $b$
is decreased towards $0$, pairs of complex conjugate critical values of
$W_b$ (those labelled $k$ and $n+2-k$, for $1\le k\le \frac{n}{2}$),
successively converge towards each other. For $2\le k<\frac{n}{2}$, the
corresponding vanishing cycles are disjoint, and the two complex conjugate
critical values essentially exchange their positions before escaping to
infinity (with complex arguments close to $\mp \frac{k-1}{n-2}\,2\pi$)
for $b\to 0$.
On the other hand, for $k=1$ the two complex conjugate critical points
labelled $1$ and $n+1$ merge and turn into two real critical points, one
of which escapes to infinity as $b\to 0$; similarly for $k=\frac{n}{2}$
if $n$ is even.

\begin{figure}[ht]
\centering
\setlength{\unitlength}{14mm}
\begin{picture}(3,3.2)(-1.5,-1.5)
\put(1,0){\circle*{0.07}}
\put(0.31,0.95){\circle*{0.07}}
\put(0.31,-0.95){\circle*{0.07}}
\put(-0.81,0.59){\circle*{0.07}}
\put(-0.81,-0.59){\circle*{0.07}}
\put(-1,0){\circle*{0.07}}
\put(-0.31,0.95){\circle*{0.07}}
\put(-0.31,-0.95){\circle*{0.07}}
\put(0.81,0.59){\circle*{0.07}}
\put(0.81,-0.59){\circle*{0.07}}
\qbezier[60](0.81,0.59)(0.2,-1.3)(1.5,-0.5)
\put(1.5,-0.5){\vector(3,2){0}}
\qbezier[60](-0.81,-0.59)(-0.2,1.3)(-1.5,0.5)
\put(-1.5,0.5){\vector(-3,-2){0}}
\qbezier[30](0.31,0.95)(0.23,0.08)(0.15,-0.8)
\qbezier[15](0.15,-0.8)(0.1,-1.2)(0.7,-1.5)
\put(0.7,-1.5){\vector(2,-1){0}}
\qbezier[30](0.31,-0.95)(0.39,-0.08)(0.47,0.8)
\qbezier[15](0.47,0.8)(0.52,1.2)(1.1,1.5)
\put(1.1,1.5){\vector(2,1){0}}
\qbezier[30](-0.31,-0.95)(-0.23,-0.08)(-0.15,0.8)
\qbezier[15](-0.15,0.8)(-0.1,1.2)(-0.7,1.5)
\put(-0.7,1.5){\vector(-2,1){0}}
\qbezier[30](-0.31,0.95)(-0.39,0.08)(-0.47,-0.8)
\qbezier[15](-0.47,-0.8)(-0.52,-1.2)(-1.1,-1.5)
\put(-1.1,-1.5){\vector(-2,-1){0}}
\end{picture}
\qquad \qquad
\begin{picture}(3,3)(-1.5,-1.5)
\put(1,0){\circle*{0.07}}
\put(0.31,0.95){\circle*{0.07}}
\put(0.23,-0.95){\circle*{0.07}}
\put(-0.81,0.59){\circle*{0.07}}
\put(-0.81,-0.59){\circle*{0.07}}
\put(-1,0){\circle*{0.07}}
\put(-0.23,0.95){\circle*{0.07}}
\put(-0.31,-0.95){\circle*{0.07}}
\put(0.81,0.59){\circle*{0.07}}
\put(0.81,-0.59){\circle*{0.07}}
\qbezier[60](0.81,0.59)(0.2,-1.3)(1.5,-0.5)
\put(1.5,-0.5){\vector(3,2){0}}
\qbezier[60](-0.81,-0.59)(-0.2,1.3)(-1.5,0.5)
\put(-1.5,0.5){\vector(-3,-2){0}}
\qbezier[30](0.31,0.95)(0.43,0.18)(0.55,-0.6)
\qbezier[15](0.55,-0.6)(0.62,-1)(1.02,-1.2)
\put(1.02,-1.2){\vector(2,-1){0}}
\qbezier[50](-0.23,0.95)(0.21,-0.12)(0.71,-1.3)
\put(0.71,-1.3){\vector(1,-2){0}}
\qbezier[30](-0.31,-0.95)(-0.43,-0.18)(-0.55,0.6)
\qbezier[15](-0.55,0.6)(-0.62,1)(-1.02,1.2)
\put(-1.02,1.2){\vector(-2,1){0}}
\qbezier[50](0.23,-0.95)(-0.21,0.12)(-0.71,1.3)
\put(-0.71,1.3){\vector(-1,2){0}}
\qbezier[120](0.7,1.5)(1.1,0.8)(1,0)
\qbezier[200](0.7,1.5)(1.7,0.4)(0.81,-0.59)
\put(0.7,1.5){\line(-3,-1){0.9}}
\put(0.7,1.5){\line(-5,-1){1}}
\qbezier[35](-0.2,1.2)(-0.4,1.13)(-0.31,0.9)
\qbezier[45](-0.3,1.3)(-0.55,1.25)(-0.41,0.9)
\put(-0.31,0.9){\line(2,-5){0.72}}
\put(-0.41,0.9){\line(2,-5){0.72}}
\qbezier[55](0.31,-0.9)(0.45,-1.3)(0,-1.4)
\qbezier[70](0.41,-0.9)(0.6,-1.4)(0,-1.5)
\qbezier[180](0,-1.4)(-1.2,-1.45)(-1,0)
\qbezier[220](0,-1.5)(-1.4,-1.6)(-1.2,0)
\qbezier[70](-1.2,0)(-1.18,0.35)(-0.81,0.59)
\put(0.65,1){\small $\tilde{\gamma}_1$}
\put(1.1,1){\small $\tilde{\gamma}_0$}
\put(0,1.1){\small $\gamma'$}
\put(-0.1,1.45){\small $\gamma''$}
\end{picture}

\caption{The deformation $b\to 0$ $(n=8)$}
\label{fig:deformfn}
\end{figure}
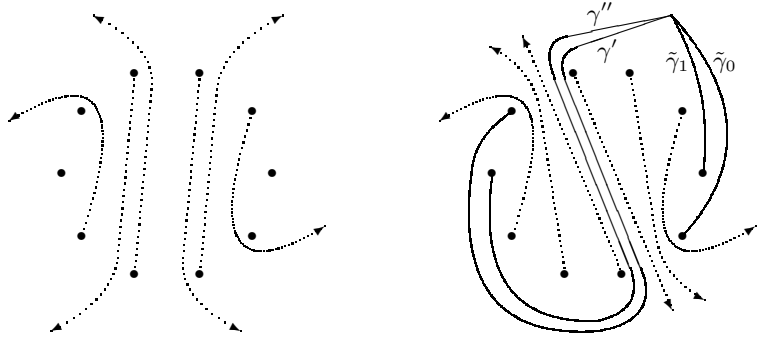

If instead of following the real axis we carry out the deformation $b\to 0$
with $\mathrm{Im}(b)$ small positive, then we can avoid all the values of $b$
for which two critical values of $W_b$ coincide, which allows us to keep
track of the manner in which $n-2$ of the critical values escape to
infinity. This is represented on Figure \ref{fig:deformfn} (left).

Observe that the vanishing cycles at the critical points corresponding
to labels in the range $1\le k<\frac{n}{2}$ are disjoint from those at the
critical points with labels in the range $\frac{n}{2}+2\le k\le n$.
Therefore, for the purposes of determining the remaining vanishing
cycles as $b\to 0$, the family of Lefschetz fibrations $W_b$ is equivalent
to one where the various critical values escape to infinity in a slightly
different manner, with the critical values coming from the
$\mathrm{Im}\,W<0$ half-plane staying ``to the left'' (towards the negative
real axis) of those coming from the $\mathrm{Im}\,W>0$ half-plane,
as pictured on Figure \ref{fig:deformfn} (right).

Therefore, if we consider the category of Lagrangian vanishing cycles
associated to a system of arcs containing the four arcs $\tilde{\gamma}_0,
\tilde{\gamma}_1,\gamma',\gamma''$ represented on Figure \ref{fig:deformfn}
right, then the full subcategory singled out by the deformation $b\to 0$ is
that generated by the four vanishing cycles $\tilde{L}_0,
\tilde{L}_1,L',L''$ associated to these arcs. A suitable collection of arcs
can be built by a sequence of sliding operations, starting from a
collection $\{\tilde{\gamma}_i,\ 0\le i\le n+1\}$ where
$\tilde{\gamma}_0$ and $\tilde{\gamma}_1$ are as pictured, and all the
$\tilde{\gamma}_i$ remain outside of the unit disc.
Identify implicitly the critical points of $W_1$ with those of
the superpotential mirror to $\CP^2(n,1,1)$, and recall that sliding
operations correspond to mutations. Then the left dual to the exceptional
collection $(\tilde{L}_0,\dots,\tilde{L}_{n+1})$
associated to the arcs $\{\tilde{\gamma}_i\}$ is equivalent (up to some
shifts) to the exceptional collection associated to the arcs
$(\gamma_2,\gamma_3,\dots,\gamma_{n+1},\gamma_0,\gamma_1)$ (using the
notation of \S \ref{sec:fsp2}). Using $\Z/(n+2)$-equivariance, the latter
is equivalent to the exceptional collection associated to the
system of arcs $(\gamma_0,\dots,\gamma_{n+1})$ considered in \S
\ref{sec:fsp2}.

Recall that the two exceptional collections for
$\db{\coh(\CP^2(n,1,1))}$ presented in \S \ref{sec:algebra} are mutually
dual (cf.\ Example \ref{ex:dual}), and that Theorem \ref{thm:fswp2}
identifies the exceptional collection associated to the arcs
$(\gamma_0,\dots,\gamma_{n+1})$ with that given by Corollary \ref{cor:excoll}.
Therefore, there is an equivalence of categories which maps the exceptional
collection $(\tilde{L}_0,\dots,\tilde{L}_{n+1})$ for $D(\FS(W_1))$ to the
exceptional collection $(\mathcal{O},\dots,\mathcal{O}(n+1))$ for
$\db{\coh(\CP^2(n,1,1))}$.

Next, let $k=\lfloor\frac{n+3}{2}\rfloor$, so that $\gamma'$ and $\gamma''$
have the same endpoints as $\tilde{\gamma}_k$ and $\tilde{\gamma}_{k+1}$
respectively. First slide $\tilde{\gamma}_{k+2},\dots,\tilde{\gamma}_{n+1}$
to the left of $\tilde{\gamma}_k$ and $\tilde{\gamma}_{k+1}$ to obtain
another system of arcs $(\tilde{\gamma}_0,\dots,\tilde{\gamma}_{k-1},
\eta_{k+2},\dots,\eta_{n+1},\tilde{\gamma}_k,\tilde{\gamma}_{k+1})$.
Then the arcs obtained by sliding $\tilde{\gamma}_k$ and
$\tilde{\gamma}_{k+1}$ to the left of $\eta_{k+2},\dots,\eta_{n+1}$ are
homotopic to $\gamma'$ and $\gamma''$. This gives us a new system of arcs
$(\tilde{\gamma}_0,\tilde{\gamma}_1,\dots,\tilde{\gamma}_{k-1},
\gamma',\gamma'',\eta_{k+2},\dots,\eta_{n+1})$, which determines a
full exceptional collection
$(\tilde{L}_0,\tilde{L}_1,\dots,\tilde{L}_{k-1},L',L'',
\Lambda_{k+2},\dots,\Lambda_{n+1})$ in $D(\FS(W_1))$.

By construction, the full subcategory $\langle \tilde{L}_0,
\tilde{L}_1,L',L''\rangle$ of the category $D(\FS(W_1))$ is
equivalent to the triangulated subcategory $\langle
\mathcal{O},\mathcal{O}(1),G_k,G_{k+1} \rangle$ of
$\db{\coh(\CP^2(n,1,1))}$, which by Lemma \ref{l:fncatmutation}
coincides with $\langle
\mathcal{O},\mathcal{O}(1),\mathcal{O}(n),\mathcal{O}(n+1)\rangle.$
As seen in \S \ref{ss:dbfn} this category is equivalent to the
derived category of the Hirzebruch surface $\mathbb{F}_n,$ which
completes the proof.
\endproof

It is also possible to prove Proposition \ref{prop:fnsubcat} by a direct
calculation involving the monodromy of $W_1,$ instead of using
Lemma \ref{l:fncatmutation}. Starting from the description of the vanishing
cycles associated to the arcs $\gamma_i$ in \S \ref{sec:fsp2}, one can
determine first the vanishing cycles
$\tilde{L}_i$ associated to $\tilde{\gamma}_i$ for all $i$, and then
those associated to $\gamma'$ and $\gamma''.$ It is
then possible to check that, although the vanishing cycles associated to
$\gamma'$ and $\gamma''$ do not quite
correspond to $\tilde{L}_n$ and $\tilde{L}_{n+1}$, after
sliding $\gamma'$ and $\gamma''$ around each other a certain number of
times one obtains two vanishing cycles that are Hamiltonian isotopic
to $\tilde{L}_n$ and $\tilde{L}_{n+1}$.

\section{Further remarks}
\nobreak
\subsection{Higher-dimensional weighted projective spaces}
\label{ss:higherdim}
\nobreak
Many of the arguments in \S \ref{sec:fsp2} extend to
higher-dimensional weighted projective spaces, working by induction
on dimension in a manner similar to the ideas in \S 5 of \cite{AuG}.
Indeed, the mirror to the weighted projective space $\CP^n(a_0,\dots,a_n)$
is the affine hypersurface $X=\{x_0^{a_0}\dots x_n^{a_n}=1\}\subset
(\C^*)^{n+1}$, equipped with the superpotential $W=x_0+\dots+x_n$ and
an exact symplectic form $\omega$ that we can choose to be invariant
under the diagonal action of $\Z/(a_0+\dots+a_n)$ and anti-invariant
under complex conjugation for simplicity. It is easy to check that $W$ has
$a_0+\dots+a_n$ critical points over $X$, all isolated and non-degenerate;
the corresponding critical values are the roots $\lambda_j$ of
$$\lambda^{a_0+\dots+a_n}=\frac{(a_0+\dots+a_n)^{a_0+\dots+a_n}}{a_0^{a_0}
\dots a_n^{a_n}}.$$
As in the two-dimensional case we use $\Sigma_0=W^{-1}(0)$ as our
reference fiber, and join it to the singular fibers of $W$ via straight
line segments $\gamma_j\subset\C$ joining the origin to $\lambda_j$.

In order to understand the vanishing cycles $L_j\subset\Sigma_0$, we consider
as before the projection to one of the coordinate axes, for example
$\pi_0:(x_0,\dots,x_n)\mapsto x_0$. For generic values of $\lambda$,
the map $\pi_0:\Sigma_\lambda\to\C^*$ defines an affine Lefschetz fibration
on $\Sigma_\lambda=W^{-1}(\lambda)$, with $a_0+\dots+a_n$ singular fibers.
These singular fibers are the preimages of the critical values of $\pi_0$ over
$\Sigma_\lambda$, which are the roots of
\begin{equation}\label{eq:branchx0}
x_0^{a_0}\,(\lambda-x_0)^{a_1+\dots+a_n}=
\frac{(a_1+\dots+a_n)^{a_1+\dots+a_n}}{a_1^{a_1}\dots a_n^{a_n}}
\end{equation}
(compare with (\ref{eq:branchx})). This equation acquires a double root
whenever $\lambda$ is one of the $\lambda_j$; the manner in which
two of the roots approach each other as one moves from $\lambda=0$ to
$\lambda=\lambda_j$ along the arc $\gamma_j$ defines an arc
$\delta_j\subset\C^*$, which is a {\it matching path} for the
Lefschetz fibration $\pi_0:\Sigma_0\to\C^*$. As in the two-dimensional
case, the Lagrangian vanishing cycle $L_j\subset\Sigma_0$ is isotopic to a
sphere $L'_j$ which lies above the arc $\delta_j$;
the generic fiber of $\pi_{0|L'_j}:L'_j\to\delta_j\subset\C^*$ is now a
Lagrangian $(n-2)$-sphere inside the fiber of $\pi_0$.

Because of the similarity between equations (\ref{eq:branchx0}) and
(\ref{eq:branchx}), it is easy to check that Lemma \ref{l:vcs}
extends almost verbatim to the higher-dimensional
case, substituting $a_0$ for $a$ and $a_1+\dots+a_n$ for $b+c$.

In order to determine the Floer complexes $CF^*(L_i,L_j)$, or equivalently
$CF^*(L'_i,L'_j)$, we need to understand, for each point of $\delta_i\cap
\delta_j$, how $L'_i$ and $L'_j$ intersect each other inside the corresponding
fiber of $\pi_0$. Because $L'_i$ and $L'_j$ each arise from matching pairs
of vanishing cycles of the Lefschetz fibration $\pi_0$, this can be done
by studying in more detail the topology of the fiber of
$\pi_0:\Sigma_0\to\C^*$ and the manner in which it degenerates as one
moves from a generic value of $x_0$ to one of the critical values.
In fact, we can use the same approach to study the vanishing cycles of
$\pi_0:\Sigma_0\to\C^*$ as in the case of $W:X\to\C$, namely project the
fiber $F_\mu=\pi_0^{-1}(\mu)$ to one of the coordinates, e.g.\ $x_1$. This
gives rise to a map $\pi_1:F_\mu\to\C^*$, which is again a Lefschetz
fibration (whose fibers are now $(n-3)$-dimensional), with $a_1+\dots+a_n$
singular fibers corresponding to values of $x_1$ that solve the equation
$$\mu^{a_0}x_1^{a_1}(-\mu-x_1)^{a_2+\dots+a_n}=
\frac{(a_2+\dots+a_n)^{a_2+\dots+a_n}}{a_2^{a_2}\dots a_n^{a_n}},$$
which presents a double root precisely when $\mu$ is a solution
of (\ref{eq:branchx0}) (for $\lambda=0$). The process can go on
similarly, considering successive restrictions to fibers and
coordinate projections until we reach the easily understood case of
$0$-dimensional fibers; once this process is completed, it
becomes possible to describe explicitly $CF^*(L'_i,L'_j)$ in terms of the
available combinatorial data. The final result is the following:

\begin{prop}\label{prop:isects3}
For $i<j$, the vanishing cycles $L'_i$ and $L'_j$ intersect
transversely, and
$$
|L'_i\cap L'_j|=\#\{I\subset\{0,\dots,n\},\ \sum_{k\in I}
a_k=j-i\}. $$
 Hence the Floer complex $CF^*(L'_i,L'_j)$ is
naturally isomorphic to the degree $j-i$ part of the exterior
algebra on $n+1$ generators of respective degrees $a_0,\dots,a_n$.
Moreover, the Floer differential is trivial, i.e.\ $m_1=0$.
\end{prop}

\begin{figure}[t]
\centering
\setlength{\unitlength}{1cm}
\begin{picture}(4,2)(-2,-1)
\put(0.8,0.8){\circle*{0.08}}
\put(-0.8,0.8){\circle*{0.08}}
\put(0.8,-0.8){\circle*{0.08}}
\put(-0.8,-0.8){\circle*{0.08}}
\put(0.9,0.8){\tiny $3$}
\put(0.9,-0.85){\tiny $0$}
\put(-1,-0.85){\tiny $1$}
\put(-1,0.8){\tiny $2$}
\put(0.38,-0.47){\tiny $p$}
\put(0.1,0.05){\tiny $\mu\!_0$}
\put(0.2,0){\circle*{0.01}}
\qbezier[120](0.8,0.8)(-0.4,1)(-0.4,0)
\qbezier[120](0.8,-0.8)(-0.4,-1)(-0.4,0)
\qbezier[120](-0.8,0.8)(0.4,1)(0.4,0)
\qbezier[120](-0.8,-0.8)(0.4,-1)(0.4,0)
\qbezier[120](0.8,-0.8)(1,0.4)(0,0.4)
\qbezier[120](-0.8,-0.8)(-1,0.4)(0,0.4)
\qbezier[120](0.8,0.8)(1,-0.4)(0,-0.4)
\qbezier[120](-0.8,0.8)(-1,-0.4)(0,-0.4)
\put(0.55,-1.10){\makebox(0,0)[cc]{\small $\delta_2$}}
\put(-1.05,-0.45){\makebox(0,0)[cc]{\small $\delta_3$}}
\put(-1.05,0.45){\makebox(0,0)[cc]{\small $\delta_1$}}
\put(-0.47,-1.10){\makebox(0,0)[cc]{\small $\delta_0$}}
\put(0,0){\circle{0.1}}
\end{picture}
\begin{picture}(4,2)(-2,-1)
\put(1,0){\circle*{0.08}}
\put(-0.5,0.866){\circle*{0.08}}
\put(-0.5,-0.866){\circle*{0.08}}
\qbezier[180](-0.5,0.866)(1,0)(-0.5,-0.866)
\qbezier[180](-0.5,0.866)(-0.5,-0.866)(1,0)
\qbezier[180](1,0)(-0.5,0.866)(-0.5,-0.866)
\put(1.6,0.4){\makebox(0,0)[cc]{\small $\pi_1(\beta_0)=\pi_1(\beta_3)$}}
\put(0.25,-0.85){\makebox(0,0)[cc]{\small $\pi_1(\beta_1)$}}
\put(-0.98,0.5){\makebox(0,0)[cc]{\small $\pi_1(\beta_2)$}}
\put(0,0){\circle{0.1}}
\end{picture}
\caption{The case of $\CP^3$: images by $\pi_0$ of the vanishing cycles
$L'_j\subset\Sigma_0$ of $W$ (left), and images by $\pi_1$ of the
vanishing cycles $\beta_j\subset F_{\mu_0}$ of $\pi_0$ (right)}
\label{fig:cp3}
\end{figure}
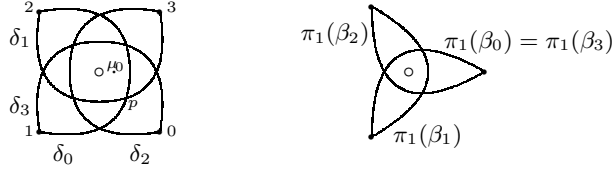

Instead of providing a complete proof, we simply illustrate Proposition
\ref{prop:isects3} by considering the example of the projective space $\CP^3$.
In that case, $\Sigma_0$ is an affine K3 surface, and $\pi_0:\Sigma_0\to\C^*$
is a fibration by affine elliptic curves, with four singular fibers.
The four vanishing cycles $L'_j\subset \Sigma_0$ project to arcs
$\delta_j\subset\C^*$ as shown on Figure \ref{fig:cp3} (left).

Using the projection $\pi_1$ to the second coordinate, we can view
each of the fibers of $\pi_0:\Sigma\to\C^*$ as a double cover
of $\C^*$ branched in 3 points (Figure \ref{fig:cp3}, right). To describe
the monodromy of the elliptic fibration $\pi_0$, we choose a
reference fiber $F_{\mu_0}=\pi_0^{-1}(\mu_0)$ for some $\mu_0\in\C^*$
close to $0$ on the positive real axis. The monodromy of $\pi_0$ around
the origin is the diffeomorphism of $F_{\mu_0}$ obtained by rotating the
three branch points of $\pi_1$ counterclockwise by $2\pi/3$.
To describe the four vanishing cycles of $\pi_0$, we join the regular
value $\mu_0$ of $\pi_0$ to each of the four critical values by considering
arcs which start at $\mu_0$, rotate clockwise around the origin from
$\arg \mu=0$ to $\arg \mu=-\frac{\pi}{4}-j\frac{\pi}{2}$ ($0\le j\le 3$),
and then go radially outwards to the corresponding critical values of
$\pi_0$. The vanishing cycles $\beta_0,\dots,\beta_3$ obtained in this
way are isotopic to the double lifts via $\pi_1:F_{\mu_0}\to\C^*$ of the
arcs shown on Figure \ref{fig:cp3} (right).

Now that the monodromy of $\pi_0$ is well-understood, it is not hard to
visualize the Lagrangian spheres $L'_j\subset\Sigma_0$ lying above the
arcs $\delta_j$, and in particular their intersections. For example,
$L'_0\cap L'_1$ consists of four points, one of which is the critical
point of $\pi_0$ with $\arg x_0=\frac{3\pi}{4}$ (lying above the common
end point of $\delta_0$ and $\delta_1$), while the three others lie in
the fiber above the other point $p$ of $\delta_0\cap\delta_1$ (with
$\arg x_0=-\frac{\pi}{4}$), and correspond (under a suitable parallel
transport operation) to the three intersections between $\beta_1$ and
$\beta_2$ in $F_{\mu_0}$. Similarly, $L'_0\cap L'_2$ consists of 6 points
(three above each point of $\delta_0\cap\delta_2$), and so on.

Finally, we observe that there cannot be any contributions to the Floer
differential $m_1$, for purely topological reasons. Indeed, if we consider
any two intersection points $p,q\in L'_i\cap L'_j$ for some pair $(i,j)$,
and any two arcs $\gamma\subset L'_i$ and $\gamma'\subset L'_j$ joining
$p$ to $q$, then $\gamma$ and $\gamma'$ are never homotopic inside
$\Sigma_0$, as easily seen by considering either $\pi_0(\gamma)$ and
$\pi_0(\gamma')$ (if $\pi_0(p)\neq \pi_0(q)$), or $\pi_1(\gamma)$ and
$\pi_1(\gamma')$ (if $\pi_0(p)=\pi_0(q)$).

The proof of Proposition \ref{prop:isects3} is essentially a careful induction on successive
slices and coordinate projections, where one manages to understand the
structure of the intersections between vanishing cycles by starting with
a 1-dimensional slice of $\Sigma_0$ and then adding one extra dimension
at a time; the main difficulty resides in setting up the induction properly and
in choosing manageable notations for the many objects that appear in the
proof, rather than in the actual calculations which are essentially always
the same.
\medskip

The next step towards understanding the category of vanishing cycles of the
Lefschetz fibration $W:X\to\C$ would be to study the moduli spaces of
pseudo-holomorphic maps from a disc with three or more marked points
to $\Sigma_0$ with boundary on $\bigcup L'_j$, something which falls
beyond the scope of this paper. Nonetheless, a careful observation
suggests that the main features observed in the two-dimensional case,
namely the vanishing of $m_k$ for $k\ge 3$ and the exterior algebra
structure underlying $m_2$, should extend to the higher-dimensional
case.

For example, in the case of $\CP^3$, we can study
$m_2:\mathrm{Hom}(L'_0,L'_1)\otimes \mathrm{Hom}(L'_1,L'_2)\to
\mathrm{Hom}(L'_0,L'_2)$ by looking carefully at Figure \ref{fig:cp3}.
Let $\alpha_0$ (resp.\ $\beta_0$) be the morphism from $L'_0$ to $L'_1$
(resp.\ from $L'_1$ to $L'_2$) which corresponds to their intersection
at a critical point of $\pi_0$, and let $\alpha_1,\alpha_2,\alpha_3$
(resp.\ $\beta_1,\beta_2,\beta_3$) be the
three other morphisms between these two vanishing cycles (labelling them
in a consistent way with respect to the other coordinate projections).
Equipping
$\Sigma_0$ with an almost-complex structure for which the projection $\pi_0$
is holomorphic, pseudo-holomorphic discs project to immersed triangular
regions in $\C^*$ with boundary on $\delta_0\cup\delta_1\cup\delta_2$;
there are three such regions (to the upper-left, to the upper-right, and
to the bottom of Figure \ref{fig:cp3} left). To start with, it is immediate
from an observation of Figure \ref{fig:cp3} that $m_2(\alpha_0,\beta_0)=0$.
Next, by deforming the arcs
$\delta_0$ and $\delta_1$ to make them lie very close to each other near
their common end point,
we can shrink the upper-left region to a very thin triangular sector,
in which case exactly one pseudo-holomorphic map contributes
to the composition of $\alpha_0$ with each of $\beta_1,\beta_2,\beta_3$.
It is then easy to see that composition with $\alpha_0$ induces an isomorphism from
$\mathrm{span}(\beta_1,\beta_2,\beta_3)\subset \mathrm{Hom}(L'_1,L'_2)$
to the subspace of $\mathrm{Hom}(L'_0,L'_2)$ spanned by the three
intersections for which $\arg x_0=\frac{\pi}{2}$. Considering the
upper-right triangular region delimited by $\delta_0,\delta_1,\delta_2$ on
Figure \ref{fig:cp3} left, we can conclude that the same is true for
the compositions of $\alpha_1,\alpha_2,\alpha_3$ with $\beta_0$, and arguing
by symmetry we can check that $m_2(\alpha_0,\beta_i)=\pm
m_2(\alpha_i,\beta_0)$ for $i=1,2,3$ (and, hopefully, a careful study of
orientations should allow one to conclude that the signs are all negative).

By a similar argument, we can study $m_2(\alpha_i,\beta_j)$ for $1\le i,j\le
3$ by shrinking the lower triangular region of Figure \ref{fig:cp3} left
to a single
point, which allows us to localize all the relevant intersection points and
pseudo-holomorphic discs into a single fiber of $\pi_0$. The intersection
pattern inside that fiber of $\pi_0$ is then described by Figure
\ref{fig:cp3} right, so that things become essentially identical to
the discussion carried out in the previous section for the Lefschetz
fibration mirror to $\CP^2$ (observe the similarity between Figures
\ref{fig:cp3} right and \ref{fig:vcs} right). Hence, the same argument as in the
two-dimensional case shows in particular that $m_2(\alpha_i,\beta_i)=0$ for
$1\le i\le 3$ and $m_2(\alpha_i,\beta_j)=\pm m_2(\alpha_j,\beta_i)$ for
$1\le i\neq j\le 3$.


\subsection{Non-commutative deformations of $\CP^2$}\label{ss:ncp2}

As mentioned in the introduction, in the general case one expects the mirror
to be obtained by partial (fiberwise) compactification of the Landau-Ginzburg
model given by the toric mirror ansatz. While not required in the toric Fano
case considered here, this fiberwise compactification allows for more
freedom of deformation, since it enlarges $H^2(X,\C)$; this sometimes makes
it possible to recover more general (non-toric) noncommutative deformations
of the Fano manifold. We now illustrate this by briefly discussing the case
of $\CP^2$. We will show the following:

\begin{prop}
Non-exact symplectic deformations of the fiberwise
compactified Landau-Ginzburg model $(\bar{X},\bar{W})$ correspond to
general noncommutative deformations of the projective plane.
\end{prop}

Moreover, we expect that there is a simple relation between the
cohomology class of the symplectic form on $\bar{X}$ and the
noncommutative deformation parameters for $\CP^2$.
\medskip

Recall that a general noncommutative projective plane is defined by a graded
regular algebra which is presented by 3 generators of degree one
and 3 quadratic relations. All these noncommutative planes were
described in the papers \cite{ATV, AS}, and with another point of
view in \cite{BP}. It was proved in \cite{ATV} that isomorphism
classes of regular graded algebras of dimension 3 generated by 3 elements
of degree 1 are in bijective correspondence with
isomorphism classes of regular triples
${\mathcal T}=(E,\sigma, L),$ where one of the following holds:

\begin{itemize}
\item[1)] $E=\PP^2,$ $\sigma$ is an automorphism of $\PP^2,$ and
$L=\O(1)$;
\item[2)] $E$ is a divisor of degree 3 in $\PP^2,$ $L$ is the restriction of
$\O_{\PP^2}(1),$ and $\sigma$ is an automorphism of $E$ such that
$(\sigma^* L)^2\cong L\otimes \sigma^{2*}L,\quad \sigma^* L\ncong L.$
\end{itemize}

The triples (and the algebras) of the first type are related to the ordinary
commutative $\PP^2$ in the sense that the category $\qgr$ of such an algebra
is equivalent to the category $\coh(\PP^2),$ whereas the triples of the second
type are related to the nontrivial noncommutative projective planes.
For example, the toric noncommutative
deformations of $\PP^2,$ which were discussed above, correspond to the
triples with $E$ isomorphic to a triangle (union of three lines).

Consider now the noncommutative projective planes which correspond to triples
with $E$ isomorphic to a smooth elliptic curve.
We know that sometimes the categories $\qgr$ of two different graded algebras
can be equivalent. In particular, with this point of view any triple with
smooth $E$
is equivalent to a triple with the same $E$ and such that $\sigma$ is a
translation by a point of $E$
(see sect.\ 8 of \cite{BP}). On the other hand, according to \cite{AS}(10.14),
the equations defining a generic
regular graded algebra, which corresponds to a triple $(E,\sigma, L)$ with
$E$ a smooth elliptic curve and $\sigma$ a translation,
can be put into the form
\begin{align*}
f_1&=cx^2 +byz+ azy=0\\
f_2&=axz +cy^2+ bzx=0\\
f_3&=bxy + ayx+ cz^2=0.
\end{align*}

This means that the DG category
${\mathfrak C}$ for these noncommutative projective planes can be described
in the following way. It  has
three objects, say $l_0, l_1, l_2,$ and for $i<j$ the spaces of
morphisms
$\Hom (l_i, l_j)$ are 3-dimensional, with all elements of degree $(j-i)$.
There are bases $x_0,y_0,z_0\in \Hom(l_0,l_1),$
$x_1,y_1,z_1\in \Hom(l_1,l_2),$
$\bar{x},\bar{y},\bar{z}\in \Hom(l_0,l_2)$ in
which the nontrivial compositions are given by the following formulas:
\begin{align*}
m_2(x_0,y_{1})=a\bar{z}, &&
m_2(x_0,z_{1})=b\bar{y}, &&
m_2(x_0,x_{1})=c\bar{x},\\
m_2(y_0,z_{1})=a\bar{x}, &&
m_2(y_0,x_{1})=b\bar{z}, &&
m_2(y_0,y_{1})=c\bar{y},\\
m_2(z_0,x_{1})=a\bar{y}, &&
m_2(z_0,y_{1})=b\bar{x}, &&
m_2(z_0,z_{1})=c\bar{z}.
\end{align*}
All other compositions (except those involving identity morphisms) vanish.

Recall from \S \ref{sec:fsp2} that the mirror of $\CP^2$ is an elliptic
fibration with three singular fibers. In the affine setting, the generic
fibers of $W=x+y+z$ on $X=\{xyz=1\}$ are tori with three punctures, but
it is possible to compactify $X$ partially into an elliptic fibration
$\bar{W}:\bar{X}\to\C$ whose fibers are closed curves; unlike what happens in
more complicated (non-toric) examples, this does not introduce any extra
critical points.

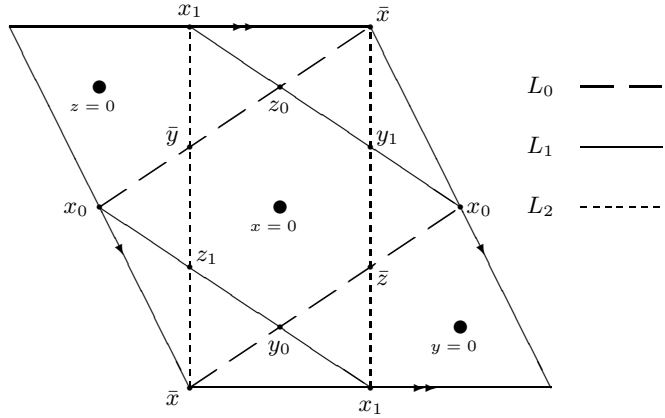
\begin{figure}[ht]
\setlength{\unitlength}{8mm}
\begin{picture}(11,7)(-1,-0.5)
\put(0,6){\line(1,-2){3}}
\put(3,0){\line(1,0){6}}
\put(0,6){\line(1,0){6}}
\put(6,6){\line(1,-2){3}}
\put(3.9,6){\vector(1,0){0}}
\put(4.1,6){\vector(1,0){0}}
\put(6.9,0){\vector(1,0){0}}
\put(7.1,0){\vector(1,0){0}}
\put(1.9,2.2){\vector(1,-2){0}}
\put(7.9,2.2){\vector(1,-2){0}}
\multiput(3,0)(0,0.2){30}{\line(0,1){0.1}}
\multiput(6,0)(0,0.2){30}{\line(0,1){0.1}}
\multiput(1.5,3)(0.66,0.44){7}{\line(3,2){0.45}}
\multiput(3,0)(0.66,0.44){7}{\line(3,2){0.45}}
\put(3,6){\line(3,-2){4.5}}
\put(1.5,3){\line(3,-2){4.5}}
\multiput(9.5,5)(0.8,0){2}{\line(1,0){0.54}}
\put(9.5,4){\line(1,0){1.35}}
\multiput(9.5,3)(0.2,0){7}{\line(1,0){0.1}}
\put(8.6,4.9){\small $L_0$}
\put(8.6,3.9){\small $L_1$}
\put(8.6,2.9){\small $L_2$}
\put(4.5,3){\circle*{0.2}}
\put(1.5,5){\circle*{0.2}}
\put(7.5,1){\circle*{0.2}}
\put(4,2.6){\tiny $x=0$}
\put(1,4.6){\tiny $z=0$}
\put(7,0.6){\tiny $y=0$}
\put(4.5,5){\circle*{0.1}}
\put(4.3,4.6){\small $z_0$}
\put(4.5,1){\circle*{0.1}}
\put(4.3,0.6){\small $y_0$}
\put(1.5,3){\circle*{0.1}}
\put(0.9,2.9){\small $x_0$}
\put(7.5,3){\circle*{0.1}}
\put(7.6,2.9){\small $x_0$}
\put(6,4){\circle*{0.1}}
\put(6.1,4.1){\small $y_1$}
\put(3,6){\circle*{0.1}}
\put(2.8,6.2){\small $x_1$}
\put(6,0){\circle*{0.1}}
\put(5.8,-0.4){\small $x_1$}
\put(3,2){\circle*{0.1}}
\put(3.1,2.1){\small $z_1$}
\put(6,6){\circle*{0.1}}
\put(6.1,6.1){\small $\bar{x}$}
\put(3,0){\circle*{0.1}}
\put(2.6,-0.3){\small $\bar{x}$}
\put(3,4){\circle*{0.1}}
\put(2.6,4.1){\small $\bar{y}$}
\put(6,2){\circle*{0.1}}
\put(6.1,1.7){\small $\bar{z}$}
\end{picture}
\caption{The vanishing cycles of the compactified mirror of $\CP^2$}
\label{fig:compactp2}
\end{figure}

The generic fiber of $\bar{W}$ and the three vanishing cycles are as represented
on Figure \ref{fig:compactp2} (compare with Figure \ref{fig:vcs} right,
which represents the images by $\pi_x$ of the same vanishing cycles; see
also Figure 2 of \cite{Se2}); the bold dots represent the intersections of
the fiber with the compactification divisor.

While it is easy to see that $m_k$ remains trivial for $k\neq 2$, the
compactification modifies the product $m_2$ in the category
$\FS(\bar{W},\{\gamma_i\})$ by introducing an infinite number of immersed
triangular regions with boundary in $L_0\cup L_1\cup L_2$. This induces
a deformation of the product structure, and the uncompactified case considered
in \S \ref{sec:fsp2} now arises as a limiting situation in which the areas
of the hexagonal regions containing the intersections with the
compactification divisor tend to infinity.

For example,
the product $m_2(x_0,y_1)$ remains a multiple of $\bar{z},$ but the relevant
coefficient is now a sum of infinitely many contributions, corresponding
to immersed triangles in which the edge joining $x_0$ to $y_1$ is an
arbitrary immersed arc between these two points inside $L_1.$ The
convergence of the series $\sum_i \pm \exp(-\mathrm{area}(T_i))$
follows directly from the fact that the area grows quadratically with the
number of times that the $x_0y_1$ edge wraps around $L_1.$
Similarly, $m_2(y_0,x_1)$ is a multiple of $\bar{z}$ as in the
uncompactified case, but with a coefficient now given by the sum of an
infinite series of contributions; and similarly for $m_2(y_0,z_1)$
and $m_2(y_1,z_0),$ which remain multiples of $\bar{x},$ and for
$m_2(z_0,x_1)$ and $m_2(x_0,z_1),$ which are proportional to $\bar{y}.$

The important new feature of the compactified Landau-Ginzburg model is
that $m_2(x_0,x_1)$ is now a multiple of $\bar{x}$ (with a coefficient that
may be zero or non-zero depending on the choice of the cohomology class of
the symplectic form); since there are again infinitely many immersed
triangular regions with vertices $x_0,x_1,\bar{x}$ (the smallest two of
which are embedded and easily visible on Figure \ref{fig:compactp2}),
the relevant coefficient is the sum of an infinite series.

Observe that the two embedded triangles are to be counted with opposite
signs (the differences in orientations at the two vertices of degree 1
cancel each other, while the non-triviality of the spin structures and the
complementarity of the sides result in a total of three sign changes, see
\S \ref{ss:antisym}); hence, in the ``symmetric'' case where the
six triangular regions delimited by $L_0\cup L_1\cup L_2$ have equal areas,
these two contributions cancel each other. The same is true of the other
(immersed) triangles with vertices $x_0,x_1,\bar{x}$, which arise in
similarly cancelling pairs. Hence, in the symmetric situation, we end up
having $m_2(x_0,x_1)=0$ as in \S \ref{sec:fsp2}; however in the general case
$m_2(x_0,x_1)$ can still be a non-zero multiple of $\bar{x}$. There are
similar statements for $m_2(y_0,y_1)$ and $m_2(z_0,z_1)$, which are
multiples of $\bar{y}$ and $\bar{z}$ respectively (and also vanish in the
symmetric case).

\subsection{HMS for products}\label{ss:product}

Let $W_1:X_1\to\C$ and $W_2:X_2\to\C$ be two Lefschetz fibrations, with
critical points respectively $p_i$, $1\le i\le r$ and $q_j$, $1\le j\le s$,
and associated critical values $\lambda_i=W_1(p_i)$ and $\mu_j=W_2(q_j)$.
Then $W=W_1+W_2:X_1\times X_2\to\C$ is a Lefschetz fibration with $rs$
critical points $(p_i,q_j)$, and associated critical values
$W(p_i,q_j)=\lambda_i+\mu_j$ (we will assume that these are pairwise
distinct and nonzero).

For all $t\in\C$, the fiber $M_t=W^{-1}(t)\subset X_1\times X_2$
can be viewed as the total space of a fibration $\phi_t:M_t\to\C$ given by
$\phi_t(p,q)=W_1(p)$, with fiber
$\phi_t^{-1}(\lambda)=W_1^{-1}(\lambda)\times W_2^{-1}(t-\lambda)$.
The $r+s$ critical values of $\phi_t$ are $\lambda_1,\dots,\lambda_r$
and $t-\mu_1,\dots,t-\mu_s.$
If $t$ varies along an arc $\gamma$ joining $0$ to
$\lambda_i+\mu_j$, the critical value $t-\mu_j$ of $\phi_t$ converges to the
critical value $\lambda_i$ by following the arc $\gamma-\mu_j$. Hence, the
vanishing cycle $L_\gamma\subset M_0$ associated to the arc $\gamma$
is a {\it fibered\/} Lagrangian sphere, mapped by $\phi_0$ to the arc
$\tilde\gamma=\gamma-\mu_j$ joining the critical values $-\mu_j$ and
$\lambda_i$ of $\phi_0$.

More precisely, the fiber of $\phi_0$ above an interior point of
$\tilde\gamma$ is symplectomorphic to the product $\Sigma_1\times \Sigma_2$
of the smooth fibers of $W_1$ and $W_2$, and its intersection with the
vanishing cycle $L_\gamma$ is a product of two Lagrangian spheres
$S_i\times T_j\subset \Sigma_1\times \Sigma_2$, where $S_i$ and $T_j$
correspond to vanishing cycles of $W_1$ and $W_2$ associated to
the critical values $\lambda_i$ and $\mu_j$ respectively. Above the
end points of $\tilde\gamma$, the product $S_i\times T_j$ collapses to
either $\{p_i\}\times T_j$ (above $\tilde\gamma(1)=\lambda_i$) or
$S_i\times \{q_j\}$ (above $\tilde\gamma(0)=-\mu_j$). Denoting by $n_i$
the complex dimension of $X_i$, a model for the topology of the restriction
of $\phi_0$ to $L_\gamma$ is given by the map $\phi:S^{n_1+n_2-1}\to [0,1]$
defined over the unit sphere in $\R^{n_1+n_2}$ by
$(x_1,\dots,x_{n_1},x_{n_1+1},\dots,x_{n_1+n_2})\mapsto
x_1^2+\dots+x_{n_1}^2.$

\begin{figure}[ht]
\setlength{\unitlength}{8mm}
\begin{picture}(6.5,4)(-1.5,-3.1)
\multiput(0,0.5)(1,0){6}{\circle*{0.1}}
\multiput(-0.5,0)(0,-1){4}{\circle*{0.1}}
\put(-0.3,0.7){\small $\lambda_1$}
\put(1.7,0.7){\small $\lambda_{i}$}
\put(3.7,0.7){\small $\lambda_{i'}$}
\put(4.7,0.7){\small $\lambda_{r}$}
\put(-1.5,-0.05){\small $-\mu_{s}$}
\put(-1.5,-1.05){\small $-\mu_{j'}$}
\put(-1.5,-2.05){\small $-\mu_{j}$}
\put(-1.5,-3.05){\small $-\mu_{1}$}
\put(-0.5,-1){\line(5,3){2.5}}
\put(-0.5,-1){\line(3,1){4.5}}
\put(-0.5,-2){\line(1,1){2.5}}
\put(0,-1.7){\small $L_{ij}$}
\put(-0.1,-0.2){\small $L_{ij'}$}
\put(1.8,-0.6){\small $L_{i'j'}$}
\end{picture}
\caption{The vanishing cycles of $W=W_1+W_2:X_1\times X_2\to\C$}
\label{fig:productvcs}
\end{figure}
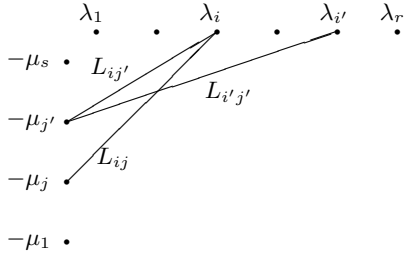

Up to a suitable isotopy we can assume that the critical values $\lambda_i$
all have the same imaginary part, and $0<\mathrm{Im}(\lambda_i)\ll
\mathrm{Re}(\lambda_1)\ll\dots\ll \mathrm{Re}(\lambda_r)$ (so that line
segments joining the origin to $\lambda_i$ form an ordered collection that
can be used to define $\FS(W_1)$). Similarly, assume that $\mu_j$ all have
the same real part, and $0<\mathrm{Re}(\mu_j)\ll\mathrm{Im}(\mu_s)\ll\dots
\ll\mathrm{Im}(\mu_1).$ Then there is a natural way to choose arcs
$\gamma_{ij},$ $1\le i\le r,$ $1\le j\le s,$ joining the origin to
$\lambda_i+\mu_j$, with both real and imaginary parts monotonically
increasing, in such a way that the lexicographic ordering of
the labels $ij$ coincides with the clockwise ordering of the arcs
$\gamma_{ij}$ around the origin. The arcs $\tilde{\gamma}_{ij}$ to which
the vanishing cycles $L_{ij}\subset M_0$ project under $\phi_0$ are then
as shown in Figure \ref{fig:productvcs}.

In this situation, we have the following result, which gives supporting
evidence for Conjecture \ref{conj:product}:

\begin{prop}
The vanishing cycles $L_{ij}$ of $W$ are in one-to-one
correspondence with pairs of vanishing cycles $(S_i,T_j)$ of $W_1$ and
$W_2$, and $$\mathrm{Hom}_{\FS(W_1+W_2)}(L_{ij},L_{i'j'})\simeq
\mathrm{Hom}_{\FS(W_1)}(S_i,S_{i'})\otimes \mathrm{Hom}_{\FS(W_2)}(T_j,T_{j'}).$$
\end{prop}

\proof[Sketch of proof]
For $i<i'$ and $j<j'$, the intersections between $L_{ij}$ and
$L_{i'j'}$ localize into a single smooth fiber of $\phi_0$, whose
intersection with $L_{ij}$ is $S_i\times T_j$ while the intersection with
$L_{i'j'}$ is $S_{i'}\times T_{j'}$ (up to isotopy in general, but by
suitably modifying the fibrations $W_1$ and $W_2$ to make them trivial
over large open subsets and by choosing the arcs $\gamma_{ij}$ carefully
we can make this hold strictly). Therefore, in this case intersections
points between $L_{ij}$ and $L_{i'j'}$ correspond to pairs of intersections
between $S_i$ and $S_{i'}$ and between $T_j$ and $T_{j'}$, so
$\mathrm{Hom}(L_{ij},L_{i'j'})\simeq \mathrm{Hom}(S_i,S_{i'})\otimes
\mathrm{Hom}(T_j,T_{j'}).$ After choosing suitable trivializations of the
canonical bundles (so that the phase of $L_{ij}$ at an
intersection point can easily be compared with the sums of the phases of
$S_i$ and $T_j$), it becomes
easy to check that this isomorphism is compatible with gradings.

When $i=i'$ and $j<j'$
the intersections between $L_{ij}$ and $L_{ij'}$ lie in a singular fiber of
$\phi_0$ (of the form $W_1^{-1}(\lambda_i)\times \Sigma_2$), inside which
$L_{ij}$ and $L_{ij'}$ identify with $\{p_i\}\times S_j$ and $\{p_i\}\times
S_{j'}$ respectively (see Figure \ref{fig:productvcs}); recalling that
$\mathrm{Hom}(S_i,S_i)=\C$ by definition, we obtain the desired formula.
Similarly for $L_{ij}\cap L_{i'j}$ when $i<i'$ and $j=j'.$ Finally, the
case $i=i'$ and $j=j'$ is trivial.

In all other cases, there are no morphisms from $L_{ij}$ to $L_{i'j'}$.
Indeed, if either $i>i'$ or $i=i'$ and $j>j'$ then $(i,j)$ follows $(i',j')$
in the lexicographic ordering, so there are no morphisms from $L_{ij}$ to
$L_{i'j'}$. The only remaining case is when $i<i'$ and $j>j'$; in that
case the triviality of $\mathrm{Hom}(L_{ij},L_{i'j'})$ follows from the fact
$L_{ij}\cap L_{i'j'}=\emptyset$ (because the projections $\tilde\gamma_{ij}$
and $\tilde\gamma_{i'j'}$ are disjoint).
\endproof

In order to prove Conjecture \ref{conj:product}, one needs to achieve a
better understanding of pseudo-holomorphic discs in $M_0$ with boundary
in $\bigcup L_{ij}$. This is most easily done in the case of low-dimensional
examples such as the mirror to $\CP^1\times\CP^1$ (already studied in a
different manner in \S \ref{ss:f0f1}), or more generally any situation
where the fibers are $0$-dimensional, because the description then becomes
purely combinatorial. Another piece of supporting evidence is the following

\begin{lemma} When $i<i'<i''$ and $j<j'<j''$, the composition
$m_2:\mathrm{Hom}(L_{ij},L_{i'j'})\otimes
\mathrm{Hom}(L_{i'j'},L_{i''j''})\to\mathrm{Hom}(L_{ij},L_{i''j''})$ is
expressed (up to homotopy) in terms of compositions in $\FS(W_1)$ and
$\FS(W_2)$ by the formula $m_2(s\otimes t,s'\otimes t')=m_2(s,s')\otimes
m_2(t,t')$.
\end{lemma}

\proof[Sketch of proof] After deforming the fibrations $W_1$ and $W_2$ and
the arcs $\gamma_{ij},\,\gamma_{i'j'},\,\gamma_{i''j''}$ (hence
``up to homotopy'' in the statement), we can
assume that all intersections between $L_{ij}$, $L_{i'j'}$ and $L_{i''j''}$
occur in a portion of $M_0$ where the fibration $\phi_0$ is trivial.
Choose an almost-complex structure which is locally a product in
$\phi_0^{-1}(U)\simeq U\times\Sigma_1\times \Sigma_2\subset M_0.$ Then every
pseudo-holomorphic disc with boundary in $L_{ij}\cup L_{i'j'}\cup
L_{i''j''}$ contributing to $m_2$ projects under $\phi_0$ to the same
triangular region in $U$ (the unique triangular region with boundary in
$\tilde\gamma_{ij}\cup\tilde\gamma_{i'j'}\cup\tilde\gamma_{i''j''}$, which
we can assume to be arbitrarily small),
while the projections to the factors $\Sigma_1$
and $\Sigma_2$ are exactly those pseudo-holomorphic discs which contribute
to $m_2:\mathrm{Hom}(S_i,S_{i'})\otimes \mathrm{Hom}(S_{i'},S_{i''})\to
\mathrm{Hom}(S_i,S_{i''})$ and $m_2:\mathrm{Hom}(T_j,T_{j'})\otimes
\mathrm{Hom}(T_{j'},T_{j''})\to\mathrm{Hom}(T_j,T_{j''}).$
\endproof

Other parts of Conjecture \ref{conj:product} are also accessible to
similar methods. However, the general situation is quite subtle, partly
because the definition of higher compositions in a product of two
$A_\infty$-categories is more complicated than one might think, but also
because one has to deal with more complicated moduli spaces of
pseudo-holomorphic discs.

\end{document}